\documentclass[times,final]{elsarticle}

\usepackage{latexsym}
\usepackage{amsmath,amssymb,amsthm}
\usepackage{algpseudocode}
\usepackage{algorithm, algorithmicx, algcompatible}
\include{definitions}

\usepackage{url}
\usepackage{xcolor}
\definecolor{newcolor}{rgb}{.8,.349,.1}

\begin{document}


\title{Multiscale Sampling for the Inverse Modeling of Partial Differential Equations}

\author[1]{Alsadig Ali}
\author[2]{Abdullah Al-Mamun}

\author[1,*]{Felipe Pereira} 

\author[3]{Arunasalam Rahunanthan}

\address[1]{Department of Mathematical Sciences, The University of Texas at Dallas, Richardson, TX, 75080, USA}
\address[2]{Institute of Natural Sciences, United International University, Dhaka-1212, Bangladesh}
\address[3]{Department of Mathematics and Computer Science, Central State University, Wilberforce, OH, 45384, USA }
\address[*]{Corresponding author: Felipe Pereira, luisfelipe.pereira@utdallas.edu }

\begin{abstract}
We are concerned with a novel Bayesian statistical framework for the characterization of natural subsurface formations, a very challenging task. Because of the large dimension of the stochastic space of the prior distribution in the framework, typically a dimensional reduction method, such as a Karhunen-Leove expansion (KLE), needs to be applied to the prior
distribution to make the characterization computationally tractable. Due to the large variability of properties of subsurface formations 
(such as permeability and porosity) it may be of value to localize the sampling strategy so that it can better adapt to large local
variability of rock properties.

In this paper, we introduce the concept of multiscale sampling to localize the search in the stochastic space. We combine the simplicity of a preconditioned Markov Chain Monte Carlo method with a new algorithm to 
decompose the stochastic space into orthogonal subspaces, through a one-to-one mapping of the subspaces to subdomains 
of a non-overlapping domain decomposition of the region of interest. The localization of the search is performed by a 
multiscale blocking strategy within Gibbs sampling: we apply a KL expansion locally, at the subdomain level. Within each subdomain,
blocking is applied again, for the sampling of the KLE random coefficients.

The effectiveness of the proposed framework is tested in the solution of inverse problems related to elliptic partial differential 
equations arising in porous media flows. We use multi-chain studies in a multi-GPU cluster to show that the new algorithm clearly 
improves the convergence rate of the preconditioned MCMC method. Moreover, we illustrate the importance of a few conditioning 
points to further improve the convergence of the proposed method.
\end{abstract}

\maketitle 

\section{Introduction}
Markov chain Monte Carlo (MCMC) methods have important applications and have experienced enormous developments since the original work of Metropolis-Hastings \cite{metropolis_1949,metropolis_1953,hastings1970}.
These methods differ on the strategy used in their sampling stage. They can be divided in two large classes, depending on whether they use or not the gradient of the likelihood function. Methods that take advantage of the gradient information typically show improvement in the convergence rate to the equilibrium distribution. These methods include, among many others, the MALA-Gibbs \cite{mala_2020}, Hamiltonian Monte Carlo \cite{HMC_2011}, and the active subspace method \cite{Colorado_2016}. There are also 
Hessian-based procedures \cite{Omar_2012}. In the case of porous media flows gradient calculations can be computationally very expensive (see \cite{moraes_2017} and references therein).

Our motivation for this work is the gradient-free uncertainty quantification in inverse problems associated with porous media flow problems in the field scale. 
For this class of problems the likelihood calculation typically requires the numerical solution of  very large problems in
fine computational grids \cite{jaramillo_2021}. There are methods known as upscaling \cite{durlofsky_2005}, that aim at solving an approximate problem (with effective coefficients)  on
coarse grids, thus reducing drastically the cost of the simulations. The idea of upscaling was used in \cite{efendiev2006} to define a two-stage or
preconditioned MCMC (see also \cite{christenfox05} where a procedure of this type was proposed). In the preconditioned MCMC a sample is first
tested through a coarse grid numerical simulation with upscaled coefficients. If it passes this filter, then a full fine grid simulation has to be performed to determine if the sample is accepted. Computationally this is a quite effective procedure because they do not require gradient calculations and samples 
can be discarded with a coarse grid simulation. The methods of \cite{efendiev2006,christenfox05} have been further investigated over the years.
Some developments of these methods include their application to flows in fractured porous media \cite{ginting2011}, 
their multi-physics version \cite{Pereira_Rahu_2015}, and their parallelization in multi-core devices \cite{Pre_fetching_Pereira_Rahu_2014}. 
This procedure has also been successfully applied in the solution of inverse problems in geophysics \cite{Georgia_2019} and
more recently a multi-level version of \cite{christenfox05} has been introduced \cite{lykkegaard_2020}. We remark that although the preconditioned 
MCMC is computationally more competitive than the Metropolis-Hastings algorithm, it still shows slow convergence for large dimensional problems
\cite{Georgia_2019}.
Thus, the development of MCMC methods that show good convergence properties and do not require gradient calculations remains as an
important area for research.

In this work, the new Multiscale Sampling Method (MSM) is proposed. It is motivated by multiscale methods for the solution of second order 
elliptic equations that are based on a domain decomposition (see \cite{mumm_2014,MRCM_2018} and references therein) that can 
produce solution for large problems taking advantage of the solution of a family of smaller boundary value problems. In these methods
the domain of the equation is decomposed into non-overlapping subdomains, local multiscale basis functions are computed 
for each subdomain, and a global interface problem is constructed and solved to couple the local solutions and produce the 
global solution.  In the Multiscale Sampling Method the domain of the partial differential equation is also decomposed into non-overlapping subdomains and a
local truncated Karhunen-Lo\`eve expansion (KLE) \cite{loeve1977} is used for each subdomain. The final stage in the construction of one sample 
consists in applying a local averaging
procedure to remove discontinuities between adjacent subdomains. The localized sampling is performed by Gibbs sampling
\cite{Gibbs_1993}.

We perform several multi-chain MCMC studies in a multi-GPU cluster to compare the convergence of the 
preconditioned MCMC with and without multiscale sampling for high-dimensional problems. We find that multiscale sampling has a 
huge impact in improving convergence rates for all problems considered. Moreover, acceptance rates also increase when 
the multiscale sampling is used. One example is included to illustrate the improved convergence that results from combining the proposed method with a few points for conditioning the field of interest. 

This work is organized as follows. We begin by describing the governing equations for contaminant transport problems in Section \ref{chapter3_sec1}. In Section \ref{subsurface_3} we present a Bayesian framework for subsurface characterization and the KLE for the dimensional reduction. We also recall methods for convergence assessment of MCMC methods. In Section \ref{multiscale} we describe the proposed method. Kriging and conditioning that will be used in our numerical studies are discussed in Section \ref{kriging_condi}. Numerical results from our experiments appear in Section \ref{results_3}. Our conclusions appear in Section \ref{concl_3}.

\section{The Model Problem}
\label{chapter3_sec1}
\subsection{Motivation}
We consider a model for contaminant transport problems (or single-phase flow problems) in a subsurface aquifer $\Omega$ with a heterogeneous permeability field. In this model, the first equation of a system of governing equations is an elliptic equation
\begin{align}
\label{elliptic_3_1}
	\begin{cases}
		\boldsymbol{u}~~~~~ &= -k(\boldsymbol{x}) \nabla p ~~~~~ \text{ in } \Omega \\
		\nabla \cdot \boldsymbol{u} &= f\qquad\qquad~ \text{ in } \Omega,
	\end{cases}
\end{align}
where $\boldsymbol{u}$ and $p$ represent the Darcy velocity and the fluid pressure, respectively, $k(\boldsymbol{x})$ is known as the absolute permeability field of the rock (a positive definite tensor), and $f$ represents sources and sinks. The elliptic equation is coupled to a hyperbolic equation
\begin{align}
	\phi(\boldsymbol x) \frac{\partial s(\boldsymbol x,t)}{\partial t} + \nabla \cdot [s(\boldsymbol x,t)\boldsymbol u(\boldsymbol x ) ] = 0, \label{eq:trans_3}
\end{align}
where $s(\boldsymbol x,t)$ is the contaminant concentration in the water and $\phi (\boldsymbol x)$ is the porosity of the rock.

The aquifer may contain many monitoring and injection wells. Our goal is to characterize the permeability field of the domain of interest by using available fractional flow data defined by
\begin{equation*}
	\displaystyle F(t) = 1-   \frac{\int_{ {\partial \Omega}_{\text{out}}} u_n(\boldsymbol{x}) s(\boldsymbol{x},t)~dy}{\int_{{\partial \Omega}_{\text{out}}}u_n(\boldsymbol{x})~dy},
\end{equation*}
where ${\partial \Omega}_{\text{out}}$ and $u_n(\boldsymbol{x})$ are the well outflow boundary and the normal components of the velocity field, respectively. More details about single-phase flow problems can be found in  \cite{chen2006,guiraldello2020}. In this paper we illustrate the proposed method in terms of the elliptic equation \eqref{elliptic_3_1}. 
\subsection{Variational Formulation of the Pressure Equation}
In this work we consider 
$\Omega \subset \mathbb{R}^{2}$, a bounded domain with a Lipschitz boundary $\partial\Omega$. For problems in $\mathbb{R}^{3}$ a formulation
similar to the one we describe here is also applicable. The velocity-pressure system  is given by Eq. \eqref{elliptic_3_1}.
In porous media flow applications typical boundary conditions that occur are
Dirichlet (the pressure is given) and Neumann (the normal component of the velocity is specified), 
which are expressed as
\[
p= g_p\in H^{\frac{1}{2}}(\partial \Omega_p), \qquad \boldsymbol{u}\cdot\boldsymbol{\hat{n}}=g _{\boldsymbol{u}}\in H^{-\frac{1}{2}}(\partial \Omega_{\boldsymbol{u}}), 
\]
where $\partial\Omega = \overline{\partial\Omega_p} \cup \overline{\partial\Omega_{\boldsymbol{u}}}$, $\overline{\partial\Omega_p} \cap  \overline{\partial\Omega_{\boldsymbol{u}}} = \varnothing $ and 
$\boldsymbol{\hat{n}}$ is the outward unit normal vector. Moreover, we assume $f\in L^2(\Omega)$.

Our numerical approximation is derived from the weak formulation of the above velocity-pressure 
problem. In order to introduce the weak formulation, we first define the following spaces
\begin{align*}
	W(\Omega) &= L^{2}(\Omega), \\
	H(\div;\Omega) &= \{\boldsymbol{v}\in (L^{2}(\Omega))^2 \ |\ \nabla \cdot {{\boldsymbol{v}}} \in L^{2}(\Omega) \}, \\
\end{align*}
and the set
\begin{align*}
	V_{ g _{\boldsymbol{u}}}(\Omega) &= \{\boldsymbol{v} \in H(\div;\Omega)\ |\ \boldsymbol{v} \cdot \boldsymbol{\hat{n}} = g _{\boldsymbol{u}}\  \hbox{on} \  \partial\Omega_{\boldsymbol u}\},
\end{align*}
for some function $g _{\boldsymbol{u}}$.
The global weak form of the pressure-velocity system \eqref{elliptic_3_1}
is given by finding $\{p,\boldsymbol{u}\}\in W\times V_{g _{\boldsymbol{u}}}$ such that
\begin{eqnarray}\label{weakp1}
	(\nabla \cdot\boldsymbol{u},w)_{\Omega}  = (f,w)_{\Omega} , ~~ \forall w\in W,~~~~~~~~~~\\
	(k^{-1}({\boldsymbol x}) \boldsymbol{u},{\boldsymbol{v}})_{\Omega} - (p,\nabla \cdot  {\boldsymbol{v}})_{\Omega} 
	+ \langle g_b, {\boldsymbol v} \cdot \boldsymbol{\hat{n}} \rangle_{\partial \Omega_p} = 0,~~ \forall {\boldsymbol{v}} \in V_0,
	\label{weakp2}
\end{eqnarray}
where $(\cdot,\cdot)_\Omega$ is the  $L^2(\Omega)$ inner product
and $\langle \cdot, \cdot \rangle_{\partial \Omega }$ is the  $L^2(\partial \Omega)$ inner product involving line
integration over $\partial \Omega$.

The system (\ref{weakp1}-\ref{weakp2}) is approximated by the lowest order Raviart-Thomas space 
\cite{raviart_1977,Douglas_Pereira_MCMCAA_1997} that is equivalent to cell-centered finite differences for a uniform
partition of $\Omega$. The resulting problem for the pressure variable is symmetric positive definite. Thus,
it can be efficiently solved by a preconditioned gradient method \cite{saad2003}. We use the algebraic 
multigrid method as a preconditioner and our elliptic solver has been developed to run on GPUs \cite{Pre_fetching_Pereira_Rahu_2014,Rahu_Pereira_MATCOM_2011}. 

\section{Subsurface Characterization}
\label{subsurface_3}
\subsection{The Bayesian Framework}
\label{Bayes_3}
Our focus in this work is in the characterization of the permeability field conditioned on pressure data. 
The available pressure data comes in the form of a red-black chessboard pattern: we assume that pressure measurements are
available at all black cells (this type of problem has been investigated in \cite{mala_2020} and references therein). 
The (log of the) permeability field is denoted by $\pmb {\eta}$ and  $R_p$ refers to the reference pressure data.  
A Bayesian statistical approach consisting of a preconditioned MCMC method combined with a novel multiscale sampling strategy 
is used to solve the inverse problem for the permeability field. The posterior probability conditioned on the pressure data $R_p$ is given 
by Bayes' rule:  
\begin{equation}
\label{bayes_eqn_3}
P(\pmb{\eta}|R_p) \propto P(R_p|\pmb{\eta})P(\pmb{\eta}),
\end{equation} 
where $P(\pmb{\eta})$ denotes a prior distribution. 
The normalizing 
constant is not needed in an iterative search within MCMC methods. 
The (log of) permeability field $\pmb{\eta}(\pmb{\theta})$ 
is built by using a local KLE strategy and 
the chain $\pmb{\theta}$ is evolved by an MCMC method.
A Gaussian likelihood function is assumed (as in~\cite{efendiev2006}), and it is given by
\begin{equation}
\label{likelihood_fun_3}
P(R_p|\pmb{\eta}) \propto \exp\Big(-(R_p - R_{\pmb\eta})^\top\Sigma (R_p - R_{\pmb\eta})\Big),
\end{equation}
where $R_{\pmb\eta}$ refers to the simulated pressure data. We set the covariance matrix  $\Sigma$  to be $\Sigma =  \pmb{I}/2\sigma_R^2$, where $\pmb{I}$ and $\sigma_R^2$ refer to the identity matrix and the precision parameter, respectively.

An MCMC algorithm is used to sample from the posterior distribution~\eqref{bayes_eqn_3}.  In the MCMC algorithm an instrumental distribution $I(\pmb{\eta}_p|\pmb{\eta})$ is used to propose a sample $\pmb{\eta}_p =\pmb{\eta}(\pmb{\theta}_p)$ at each iteration,
where $\pmb{\eta}$ denotes the previously accepted sample. For a given permeability field, the system \eqref{weakp1}-\eqref{weakp2} is solved
numerically to give $R_{\pmb\eta}$, and the original Metropolis-Hastings \cite{metropolis_1949,metropolis_1953,hastings1970} acceptance probability of a proposed sample is given by   
\begin{equation}
\label{single_stage_prob_3}
{\alpha}(\pmb{\eta}, \pmb{\eta}_p) = \text{min}
\left(1,\frac{I(\pmb{\eta}|\pmb{\eta}_p)P(\pmb{\eta}_p|R_p)}{I(\pmb{\eta}_p|\pmb{\eta})P(\pmb{\eta}|R_p)}\right).
\end{equation} 
In this work we consider a preconditioned MCMC method that will be discussed in Section \ref{msm_section}.
\subsection{Dimensional Reduction}
\label{kle_subsection}
We use a Bayesian statistical framework along with MCMC methods where our numerical simulator requires a permeability value in each cell of a partition of the domain of interest. Therefore, we need to generate a large number of random permeability values (based on the grid size) in each iteration that is infeasible from a practical point of view. Thus, we need to reduce the dimension of the uncertainty parameter space describing the permeability field. We use KLE~\cite{loeve1977,ginting2013} to achieve the desired dimensional reduction of the parameter space. Next, we briefly discuss the KLE. 

We consider $\log\left[ k(\boldsymbol{x})\right]= Y^k(\boldsymbol{x})$ to be a sample of a Gaussian field, where $k(\boldsymbol{x})$ represents the permeability field and  $\boldsymbol{x}$ is a point in the domain $\Omega$. We also consider $Y^k(\boldsymbol{x})\in L^2(\Omega)$ with unit probability, i.e., $Y^k(\boldsymbol{x})$ is a second-order stochastic process. If we assume $E[(Y^k)^2] = 0$, then, the permeability field $Y^k(\boldsymbol{x})$ can be written for a given orthonormal basis $\left \{\varphi_i\right\}$ of $L^2(\Omega)$ as follows:
\begin{equation}
\label{kle1_3}
Y^k(\boldsymbol{x}) = \sum _{i=1}^{\infty} Y_i^k, \varphi_i(\boldsymbol{x}), \hspace*{0.3cm}
\end{equation}
where $Y_i^k = \int_{\Omega}Y^k(\boldsymbol{x})$ are random coefficients, and $\varphi_i(\boldsymbol{x})$ 
are eigenfunctions with the corresponding eigenvalues $\lambda_i = E[(Y_i^k)^2] > 0$. The pairs $(\lambda_i,\varphi_i(\boldsymbol{x}))$ satisfy the integral equation
\begin{equation}
\label{kle_efun_3}
\int_{\Omega}R(\boldsymbol{x}_1, \boldsymbol{x}_2)\varphi_i(\boldsymbol{x}_2)d\boldsymbol{x}_2 = \lambda_i \varphi_i(\boldsymbol{x}_1),~~ i = 1, 2, . . .
\end{equation}
for a given covariance function $R(\boldsymbol{x}_1,\boldsymbol{x}_2)$. Setting $\theta_i^k = Y_i^k/\sqrt{\lambda_i}$ in Eq. \eqref{kle1_3} we have
\begin{equation}
\label{kle_exp_3}
Y^k(\boldsymbol{x}) = \sum_{i=1}^{\infty} \sqrt{\lambda_i}{\theta}_i^k\varphi_i (\boldsymbol{x}),
\end{equation}
where the eigenfunctions $\varphi_i$ and the corresponding eigenvalues $\lambda_i$  satisfy Eq. \eqref{kle_efun_3}.  
The eigenvalues are assumed to be arranged in descending order. The series shown in the Eq. \eqref{kle_exp_3} is known as the Karhunen-Lo\`eve expansion. The first $N$ dominating eigenvalues are considered in the KLE so that the energy $E$ is above $95\%$ \cite{Laloy2014}, i.e., 
\begin{equation}
\label{energy_3}
E=\dfrac{\displaystyle\sum_{i=1}^N \lambda_i}{\displaystyle\sum_{i=1}^{\infty} \lambda_i}\geq 95\% .
\end{equation}
We, thus, can define the truncated KLE by
\begin{equation}
\label{kle_truncated_3}
Y^k_{N}(\boldsymbol{x})=\displaystyle \sum_{i=1}^{N}\sqrt{\lambda_i}\theta_i^k
\varphi_i(\boldsymbol{x}).
\end{equation}

\subsection{Convergence Assessment of MCMCs}
\label{mpsrf_3}
We consider a problem that consists of sampling the permeability field conditioned on pressure measurements. We use a Bayesian statistical approach (discussed in Section \ref{Bayes_3}) along with a preconditioned MCMC method \cite{efendiev2006,christenfox05} to characterize the permeability field of our domain of interest. Two critical issues, namely, where to begin (burn-in) and when  to terminate (convergence), need to be addressed when MCMC methods are used. We now discuss the convergence of MCMC methods that we use in our investigation to construct one of the rock properties (permeability field).

A number of convergence criteria \cite{brooks1998mcmc,polson1996,rosenthal1995} for MCMCs have been developed with a solid theoretical foundation. Several review papers, where authors used MCMC convergence diagnostics, are available in the literature~\cite{Roy2020,cowles1996,mengersen1999,Brooksgelman1998}. Note that in \cite{cowles1996} the authors discussed thirteen MCMC convergence diagnostics. The convergence diagnostics described in~\cite{Brooksgelman1998} are now widely used. In this work we use two popular diagnostic tools, namely, the Potential Scale Reduction Factor (PSRF) and the multivariate PSRF (MPSRF), to diagnose the convergence of MCMC algorithms. Between these two, the MPSRF method takes all the parameters into account for accessing convergence of the MCMC methods. Thus, the MPSRF is more restrictive than the PSRF.

The PSRF and MPSRF measures rely on multiple chains. Thus we are required to run $m>1$ independent chains in parallel  with different initial points drawn from an overdispersed distribution. The effect of starting at different initial points is made minimal by discarding the first few iterations as burn-in. Let us denote by $\boldsymbol\theta$ an $N$-dimensional parameter vector, and let $l$ represents the number of posterior draws for each of the $m$ chains. Furthermore, assume that $\boldsymbol {\theta}_j^{c}$ denote the value of the parameter vector $\boldsymbol {\theta}$ generated at iteration $c$ in $j$th chain of the MCMC algorithm. The posterior variance-covariance matrix is then estimated as \\
\begin{equation}
\mathbf{\widehat{V}} = \frac{l-1}{l}\mathbf{W} + \left( 1+ \frac{1}{m}\right)\frac{\mathbf{B}}{l}.
\end{equation}
The within- and between-sequence (chain) covariance matrix $\mathbf{W}$ and $\mathbf{B}$ are calculated as 
\begin{equation}
\mathbf{W} = \frac{1}{m(l-1)} \sum\limits_{j=1}^m \sum\limits_{c=1}^l \left(\boldsymbol {\theta}_j^{c} - \boldsymbol {\bar \theta}_{j.}\right) \left(\boldsymbol {\theta}_j^{c} - \boldsymbol {\bar \theta}_{j.}\right) ^{T},
\end{equation}
and 
\begin{equation}
\mathbf{B} = \frac{l}{m-1} \sum\limits_{j=1}^m \left(\boldsymbol {\bar \theta}_{j.}-\boldsymbol {\bar \theta}_{..}\right) \left(\boldsymbol {\bar \theta}_{j.}-\boldsymbol {\bar \theta}_{..}\right)^{T},
\end{equation} respectively.  $\boldsymbol {\bar \theta}_{j.}$ denote within chain mean and  $\boldsymbol {\bar \theta}_{..}$ represent the mean between $m$ combined chains, respectively. $T$ denotes the transpose of a matrix. The PSRFs are calculated using the two estimators $\mathbf{\widehat{V}}$ and  $\mathbf{W}$ defined by
\begin{equation}
\begin{aligned}
\text{PSRF}_\text{i} =\sqrt{\frac {\text{diag}(\mathbf{\widehat{V}})_i}{\text{diag}(\mathbf{W})_i}}, ~~~~\text{where} \, \,\, i=1,2,...,N.
\end{aligned} 
\end{equation}
A large PSRF$_{i}$ suggests that either the estimate of the between variance can be decreased by taking more samples into account or by taking further samples one could increase the within variance. It indicates that the simulated sequences have not yet traversed the parameter space completely. On the other hand, if the maximum of PSRF values is close to $1$, we can draw the conclusion that each of the $ m $ chains of $l$ simulated samples is close to the target distribution. The MPSRF is estimated by using the maximum root statistic. As in \cite{Brooksgelman1998} it is defined by

\begin{equation}
\begin{aligned}
\text{MPSRF} &= \sqrt{\max_{\bold a} \smallskip \frac{\bold a^T\bold{\widehat{V}} \bold a}{\bold a^T\bold{W}\bold a}}\\
&=\sqrt{\max_{\bold a} \smallskip \frac{\bold a^T\left[ \frac{l-1}{l}\bold{W} + \left( 1+ \frac{1}{m}\right)\frac{\bold{B}}{l} \right] \bold a }{\bold a^T \bold{W} \bold a}} \\
&=\sqrt{\frac{l-1}{l} + \left(\frac{m+1}{m}\right) \max_{\bold a} \frac{\bold a^T\frac{\bold{B}}{l}\bold a}{\bold a^T\bold{W} \bold a} \notag}\\ 
&=\sqrt{\frac{l-1}{l} + \left(\frac{m+1}{m}\right) \lambda},
\end{aligned} 
\end{equation} 
\noindent where $\bold a \in \mathbb{R}^N$ is an arbitrary vector, and $\lambda$ is the greatest eigenvalue of the positive definite matrix $\bold{W}^{-1} \bold{B}/l$. If the means of between chains are equal, the between chain covariance matrix $  \mathbf{B}$ 
becomes zero. In this case,  the chains mix well and $\lambda\rightarrow 0$.  Thus, as the MPSRF approaches to $1$, it guarantees a convergence for sufficiently large sample size. 

\section{Multiscale Sampling}
\label{multiscale}
\subsection{The Multiscale Prior Distribution}\label{prior}
\begin{figure}[ht]
	\centering
	{
		\includegraphics[scale=1.0]{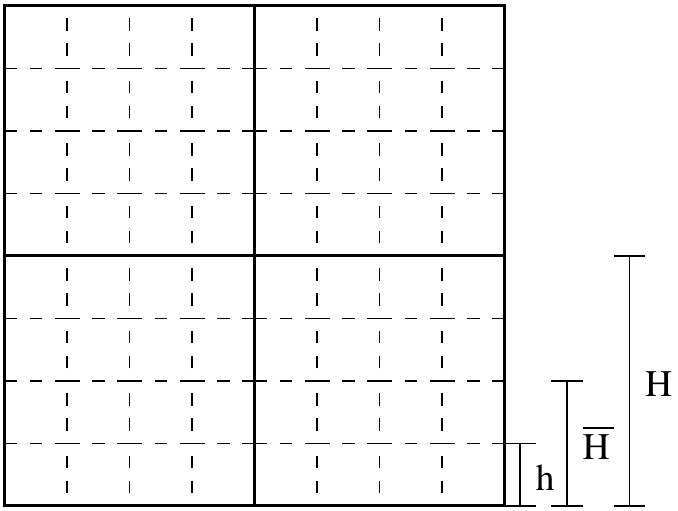}
	}
	\caption{The fine $\Omega^f$ (dashed lines) and coarse $\Omega^c$ (solid lines)
		partitions of $\O$ along with the three spatial scales used in
		the definition of the new multiscale procedure.}
	\label{fig1}
\end{figure}

We begin with the description of a decomposition of the domain $\Omega$. Our multiscale sampling strategy is based 
on two non-overlapping partitions of the domain $\Omega$: the first is a uniform fine Cartesian mesh 
$\Omega^f$ where the values of the absolute permeability field are piecewise constant. 
This is also the mesh used for the numerical solution of the system \eqref{weakp1}-\eqref{weakp2}. The second is a coarse 
Cartesian mesh $\Omega^c$ constructed as sets of elements in $\Omega^f$ (see Figure  \ref{fig1}) where a KLE will be 
applied for local dimensional reduction. The proposed method is based on 
partitions $\Omega^\gamma$ into rectangles
$\{\Omega_i^\gamma,\ i=1,\dots,M_\gamma\}$ (see Figure  \ref{fig1}), such that
\[
\bar{\Omega^\gamma}=\bigcup^{M_\gamma}_{i=1}\bar{\Omega_i^\gamma}; \quad \Omega_i^\gamma \cap \Omega_k^\gamma = 
\varnothing, \quad i \neq k, \quad \gamma = c,f.
\]
Define $\Gamma=\partial\Omega$ and, for $i=1,\dots,M_\gamma$:
\[
\quad\Gamma_{ik}^\gamma=\Gamma_{ki}^\gamma=
\partial\Omega_i^\gamma\cap\partial\Omega_k^\gamma, \quad \gamma = c,f.
\]

For each element of the coarse partition $\{\Omega_i^\gamma,\ i=1,\dots,M_c\}$, we define the set 
\[
\mathcal{S}_i = \{j: \Omega_j^f  \subset  \Omega_i^c \}.
\]
As indicated in Figure  \ref{fig1}, we refer to two length scales in the description of the new multiscale procedure:
$H$, the mesh size for the coarse partition and $h$, the mesh size of an underlying fine 
grid.

\begin{figure}[H]
	\centering
	{
		\includegraphics[scale=0.5]{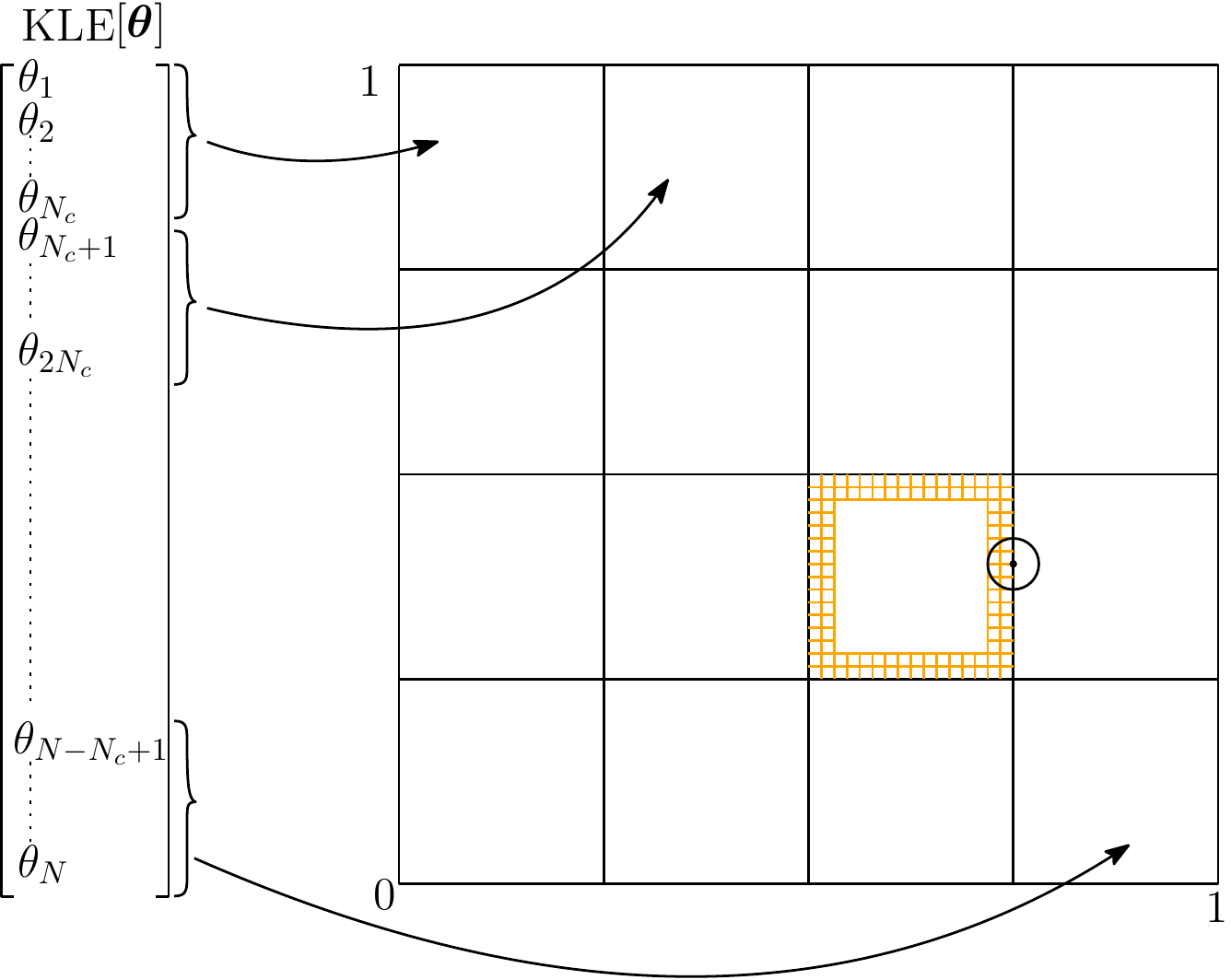}
	}
	\caption{The mapping between blocks of theta variables and elements of the coarse $\Omega^c$ (solid lines)
		partition of $\Omega$. Each block of thetas is used to generate a local Gaussian field within its corresponding subdomain.}
	\label{fig_map}
\end{figure}

We will consider a blocking strategy \cite{Brooksgelman1998} for Gibbs sampling within a Metropolis-Hastings algorithm. In order to define it we decompose the $\boldsymbol\theta$ vector in Eq. \eqref{kle_truncated_3} into  orthogonal subspaces corresponding to blocks with the same number of
components, that are denoted by $\boldsymbol\theta^i$, for $i=1,\dots,M_c$.
Each block of thetas is used to generate a local Gaussian field within its corresponding subdomain, as illustrated in
Fig. \ref{fig_map}.
The update of each  $\boldsymbol\theta^i$ block is based on the random walk sampler (RWS) of ~\cite{Cotter_2013}. It is given,
for  $i=1,\dots,M_c$, by 
\begin{equation}
\label{RW_sampler}
{\boldsymbol\theta_p}^i = \sqrt{1-\beta^2}\, {\boldsymbol\theta}^i + \beta\,{\boldsymbol\epsilon}^i, 
\end{equation} 
where the current sample is denoted by ${\boldsymbol\theta_p}^i$ and the previously accepted sample by $\boldsymbol\theta^i$. 
The algorithmic parameter $\beta$ is used for tuning the sampler and ${\boldsymbol\epsilon}^i$ represents a $\mathcal N(0,1)$-random 
vector. Not all components of ${\boldsymbol\theta_p}^i $ are updated simultaneously. Blocking is
used again so that only a subset of the components is modified in one MCMC iteration. We view ${\boldsymbol\theta_p}^i $ as a
column vector and we are going to refer to the {\it local blocking number} as the number of contiguous components of this column vector
that are updated simultaneously. For the purpose of sampling, we consider the ${\boldsymbol\theta_p}^i$ components ordered 
by their corresponding eigenvalues in the local KLE, from the largest towards the smallest one.

\begin{figure}[H]
	\centering
	{
		\includegraphics[scale = 0.45]{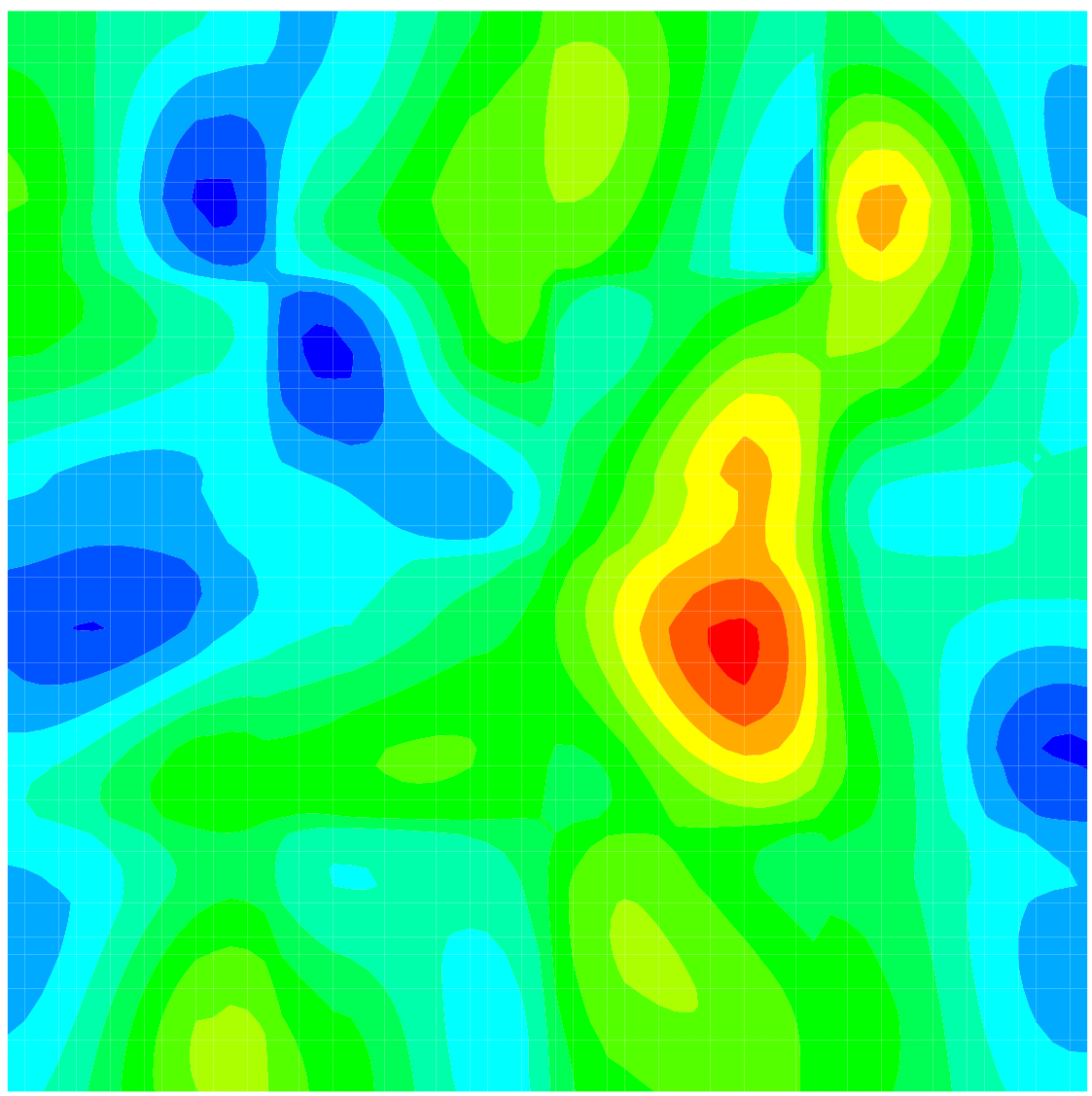} 
		\includegraphics[scale = 0.45]{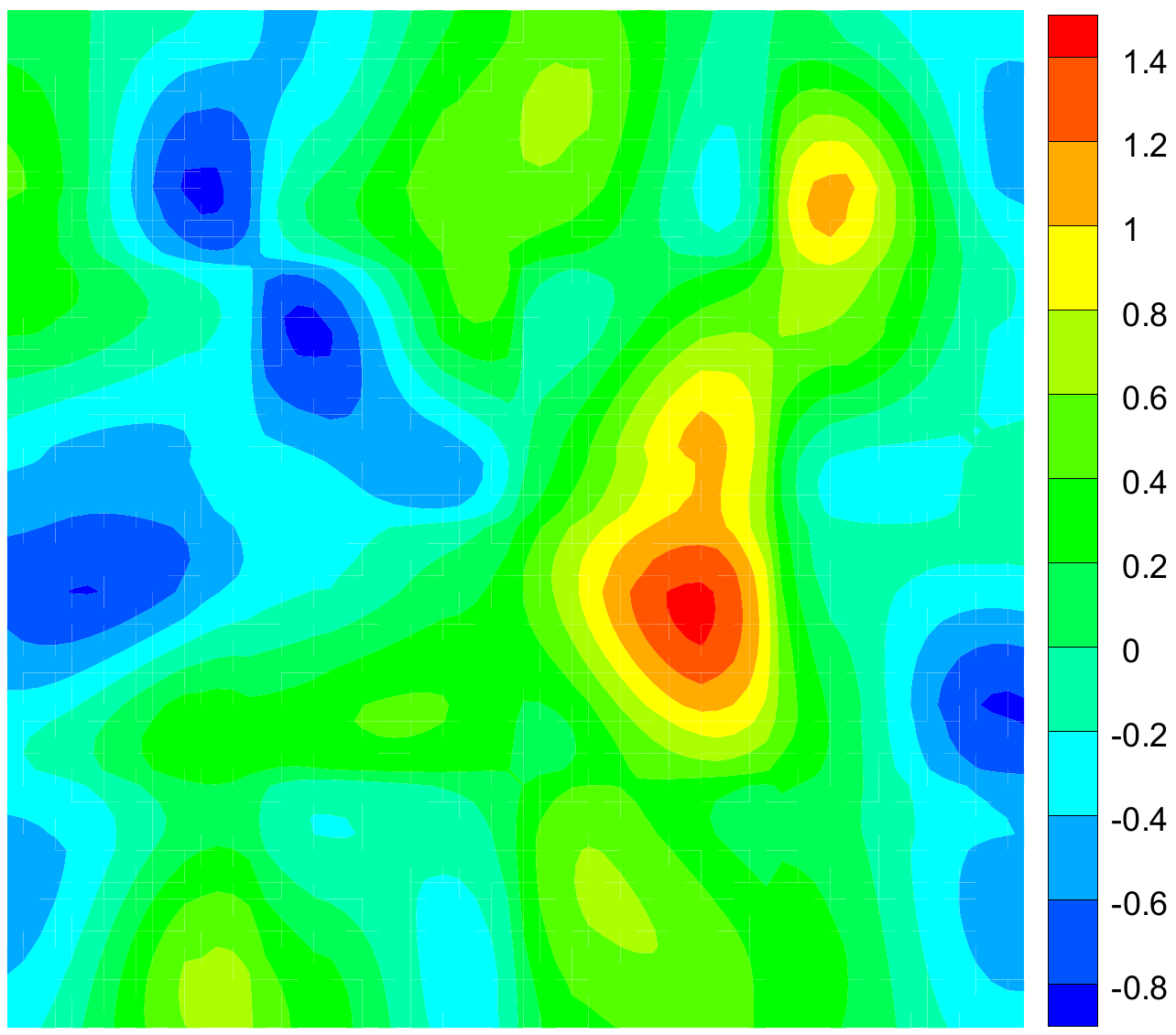}
	}
	\caption{Sample of a permeability field before (left) and after (right) a local averaging is applied at the top right subdomain of a domain decomposition with $H = 0.25$.}
	\label{fig_jump}
\end{figure}

The samples produced by the local sampling strategy discussed above produces Gaussian samples that show discontinuities in 
$\Gamma_{ik}^c$. Motivated by downscaling strategies developed for multiscale methods (that aim at removing flux discontinuities at
subdomain boundaries \cite{guiraldello2020}) in order to complete the construction of one sample from our multiscale prior 
distribution an averaging method is used to condition each sample on the available data at nearest neighbor subdomains. The
averaging procedure is illustrated in Figure  \ref{fig_map}. 
A length scale $\overline H$ is set (a fraction of the correlation length 
that enters in the construction of the prior distribution) and, for each $i$, all cells of $\mathcal{S}_i$ that are at a distance of $\overline H$ (or less)
to $\Gamma_{ik}^c$ have their current values replaced by local averages (that preserve both their mean value and variance - if they were
uncorrelated). Figure  \ref{fig_jump} illustrates a sample before and after this averaging procedure. In this figure, $H = 0.25$ and the averaging is applied on the boundary of the top right subdomain. Note that if the correlation lengths are not equal, the circle for the averaging in Figure  \ref{fig_map} should be replaced by an ellipse.

In conclusion, the multiscale prior distribution requires three user-specified parameters:
\begin{itemize}
	\item The value of $H$: the subdomain size;
	\item The value of $\overline H$: the length scale local averages are taken;
	\item The local blocking number.
\end{itemize}
We remark that the number of blocks for the localized Gibbs sampling in each subdomain is given by $N_{\text{local}} = N/(M_c * N_{lb})$, where
$N_{lb}$  denotes the local blocking number. We refer to these blocks as $B_k$, $k=1,\dots, N_{\text{local}}$, and we also define the local stochastic dimension to be $N_c=N/M_c$. We conjecture that improved convergence for MCMCs should be observed for $H>L$, where $L$ denotes the largest correlation length in the definition of the prior distribution. Further studies are needed to check the validity of this conjecture.

\subsection{The Multiscale Sampling Method}\label{msm_section}
We now provide a detailed algorithm of the new Multiscale Sampling Method (MSM) that consists of a preconditioned MCMC with a multiscale prior distribution. 
If the number of subdomains $M_c = 1$, then the proposed algorithm reduces to the classical
preconditioned MCMC ~\cite{efendiev2006,christenfox05} with Gibbs sampling associated with the local blocking number.

We first discuss the algorithm of the preconditioned MCMC method. The filtering step of this method is based on a coarse-scale model approximation of the governing system \eqref{weakp1}-\eqref{weakp2}. The coarse-scale discretization is similar to the fine-scale discretization and the permeability field $\pmb{\eta}(\pmb{\theta})$ is projected on the coarse-scale. An upscaling procedure~~\cite{durlofsky1991} is used to set an effective permeability field that provides a similar average response as that of the underlying fine-scale problem. The numerical simulator is run on the coarse-scale model and produces the coarse-grid pressure field $R_c$. The coarse-scale and fine-scale acceptance probabilities are estimated as \begin{equation}
\begin{aligned}
{\alpha}_c(\pmb{\eta}, \pmb{\eta}_p) &= \text{min}\left(1,\dfrac{I(\pmb{\eta}|\pmb{\eta}_p)P_c(\pmb{\eta}_p|R_p)}{I(\pmb{\eta}_p|\pmb{\eta})P_c(\pmb{\eta}|R_p)}\right),\text{and}\\
{\alpha}_f(\pmb{\eta},\pmb{\eta}_p) &= \text{min}\left(1,\dfrac{P_f(\pmb{\eta}_p|R_p)P_c(\pmb{\eta}|R_p)}{P_f(\pmb{\eta}|R_p)P_c(\pmb{\eta}_p|R_p)}\right),
\end{aligned}		
\end{equation}
where $P_c$ and $P_f$ are the posterior probabilities calculated at coarse- and fine-scale, respectively.
In MSM we first construct a local permeability field $\pmb\eta(\pmb\theta^i)$ for each subdomain $\Omega^c_i, i=1,\dots M_c$ using Eq. \eqref{kle_truncated_3}. To do so, we generate local KLE data in subdomains $\Omega^c_i, i=1,\dots M_c$ with size $\dfrac{1}{2^n}$, $n=1,2,\dots $. Then we construct the global permeability field by taking local averages. The MSM algorithm is presented in Algorithm \ref{alg_multiscale_a}.
\begin{algorithm}[H]
	\caption{The Multiscale Sampling Method (MSM)}
	\label{alg_multiscale_a}
	\begin{algorithmic}[1]
		\STATE For a given covariance function $R$ solve Eq. \eqref{kle_efun_3} to get a KLE in Eq. \eqref{kle_truncated_3}, which is used in all the subdomains, $\Omega_i^c, i=1,\dots,M_c$.
		\FOR{$j=1$ to $M_{\text{mcmc}}$}
		\FOR{$i=1$ to $M_c$}
		\FOR{$k=1$ to $N_{\text{local}}$}								
		\STATE  Generate  i.i.d., $\mathcal N(0,1)$ Gaussian variables to construct $\pmb{\theta}_p$ using Eq. \eqref{RW_sampler} for block $B_k$ in $\Omega_i^c$.
		\STATE Construct a local Gaussian sample (in each subdomain) using the KLE to set a preliminary 
		value for the Gaussian sample at the $\mathcal{S}_i$ cells.
		\STATE Run the local averaging algorithm to remove discontinuities.
		\STATE Compute the upscaled permeability on the coarse-scale using $\pmb{\eta}_p$.
		\STATE Solve the forward problem on the coarse-scale to get $R_c$.
		\STATE Compute the coarse-scale acceptance probability ${\alpha}_c(\pmb{\eta}, \pmb{\eta}_p)$.
		\IF{ $\pmb{\eta}_p$ is accepted}
		\STATE Use $\pmb{\eta}_p$ in the fine-scale simulation to get $R_f$.
		\STATE Compute the fine-scale acceptance probability ${\alpha}_f(\pmb{\eta},\pmb{\eta}_p)$.
		\IF{ $\pmb{\eta}_p$ is accepted}
		$\pmb{\eta} = \pmb{\eta}_p$. 
		\ENDIF		
		\ENDIF		
		\STATE $j=j+1.$
		\ENDFOR
		\ENDFOR			
		\ENDFOR
	\end{algorithmic}
\end{algorithm}

\section{Kriging and Conditioning}
\label{kriging_condi}
In this section, we combine the multiscale sampling method with conditioning by projection for the sampling. The conditioning by projection method has been discussed in detail in \cite{ali2021conditioning}. This method consists of two steps. In the first step, for given permeability values at sparse locations in the domain, generate a kriged field for the domain. In the second step, project the i.i.d $\mathcal N(0,1)$-random vector onto the nullspace of a data matrix defined in terms of KLE to calculate the linear combination in Eq. \eqref{kle_truncated_3}. The final permeability field is obtained by adding the fields defined in both steps (after taking exponential). In the following subsections we briefly discuss the kriging interpolation and projection method for conditioning.

\subsection{Kriging Interpolation}
Kriging is an interpolation method that is derived from a regionalized variable theory \cite{kriging_1978,mining_1976}. It employs a limited set of sampled data points to compute the value of a variable throughout a continuous spatial field. Kriging uses the spatial correlation between sampled points to interpolate the values in the spatial field. This interpolation gives the exact values of the field at the known locations.

\subsection{The Projection Method for Conditioning}
Following the discussion in \cite{ali2021conditioning}, we need to extract a data matrix that is defined in terms of KLE. We assume the Gaussian field $Y(\boldsymbol{x})$ defined in Eq. \eqref{kle_truncated_3} is a Gaussian perturbation on top of a kriged field $\widehat{Y}(\boldsymbol{x})$, thus we can write
\begin{equation}
\begin{aligned}
\label{condi_perm1}
Y(\boldsymbol{x}) - \widehat{Y}(\boldsymbol{x}) = \sum_{i=1}^{N}\sqrt{\lambda_i}\varphi_i(\boldsymbol{x})\theta_i
  =\pmb{\phi}^T(\boldsymbol{x})\sqrt{D}\pmb{\theta},
\end{aligned}  
\end{equation}
where, for each $\boldsymbol{x}, \pmb{\phi}(\boldsymbol{x})\in \mathbb{R}^N$, and $D$ is a diagonal matrix containing $N$ dominating eigenvalues. If we have $M$ measured data values at sparse locations, then we can define the following homogeneous linear system of equations
\begin{equation*}
\begin{aligned}
\label{system}
A\pmb{\theta}= \boldsymbol{0},
\end{aligned}  
\end{equation*}
where $A=\pmb{\phi}^T(\boldsymbol{\hat{x}})\sqrt{D}\in \mathbb{R}^{M\times N}$ is the desired data matrix. Finally, we project the vector $\pmb{\theta}$ onto the nullspace of the matrix $A$ to get the closest vector to $\pmb{\theta}$ in the nullspace of the data matrix $A$. i.e. 
\begin{equation*}
\label{projection}
\pmb{\hat{\theta}} = P\pmb{\theta}, 
\end{equation*}
 where $P$ is a projection matrix \cite{LA_strang_2019}. Therefore, we can write 
\begin{equation*}
\begin{aligned}
\label{condi_perm3}
Y(\boldsymbol{x}) = \widehat{Y}(\boldsymbol{x})+ \sum_{i=1}^{N}\sqrt{\lambda_i}
\varphi_i(\boldsymbol{x})\hat{\theta_i}, \quad \pmb{\hat{\theta}} = (\hat{\theta}_1, \dots , \hat{\theta}_N).
\end{aligned}  
\end{equation*}

\section{Numerical Results}
\label{results_3}
In this section, we describe the simulation study for the problem of interest. We test the proposed multiscale sampling method in four examples. In each example, we  numerically solve the system containing Eqs. \eqref{weakp1}-\eqref{weakp2} on the domain $\Omega =[0,1]\times [0,1]$. We present a comparative study between the preconditioned MCMC method with and without multiscale sampling in the first three examples. In the last example, we analyze MSM with and without conditioning for a problem of higher dimensional stochastic space. In MSM, we apply KLE to construct a permeability field for each subdomain, and then construct the global permeability field. In KLE, we use the following covariance function:

\begin{equation}
\label{kle_3}
\begin{aligned}
R(\boldsymbol{x}_1,\boldsymbol{x}_2) = \sigma_Y^2\, \text{exp}\left(-\frac{|x_1 - x_2|^2}{2L_x^2} - \frac{|y_1 - y_2|^2}{2 L_y^2}\right),
\end{aligned}
\end{equation}
where $L_x$ and $L_y$ are the correlation lengths and $\sigma_Y^2= \text{Var}[(Y^k)^2]$. We take $\sigma_Y^2=1$ in all the four examples. Moreover, we set the source term $f=0$ and impose Dirichlet boundary conditions, $p=1$ and $p=0$, on the left and right boundaries, respectively. We also set a no-flow (Neumann-type boundary condition) condition  on the other two boundaries. We run four MCMC chains for each method. In order to remove the discontinuities between subdomains in our numerical studies, we set the length scale $\overline H$ to be $\overline H = \min\{\frac{L_x}{2},\frac{L_y}{2}\}$. Below we discuss the numerical results.

\subsection{Example 1} 
\begin{figure}[H]
	\centering
	\includegraphics[scale = 0.7]{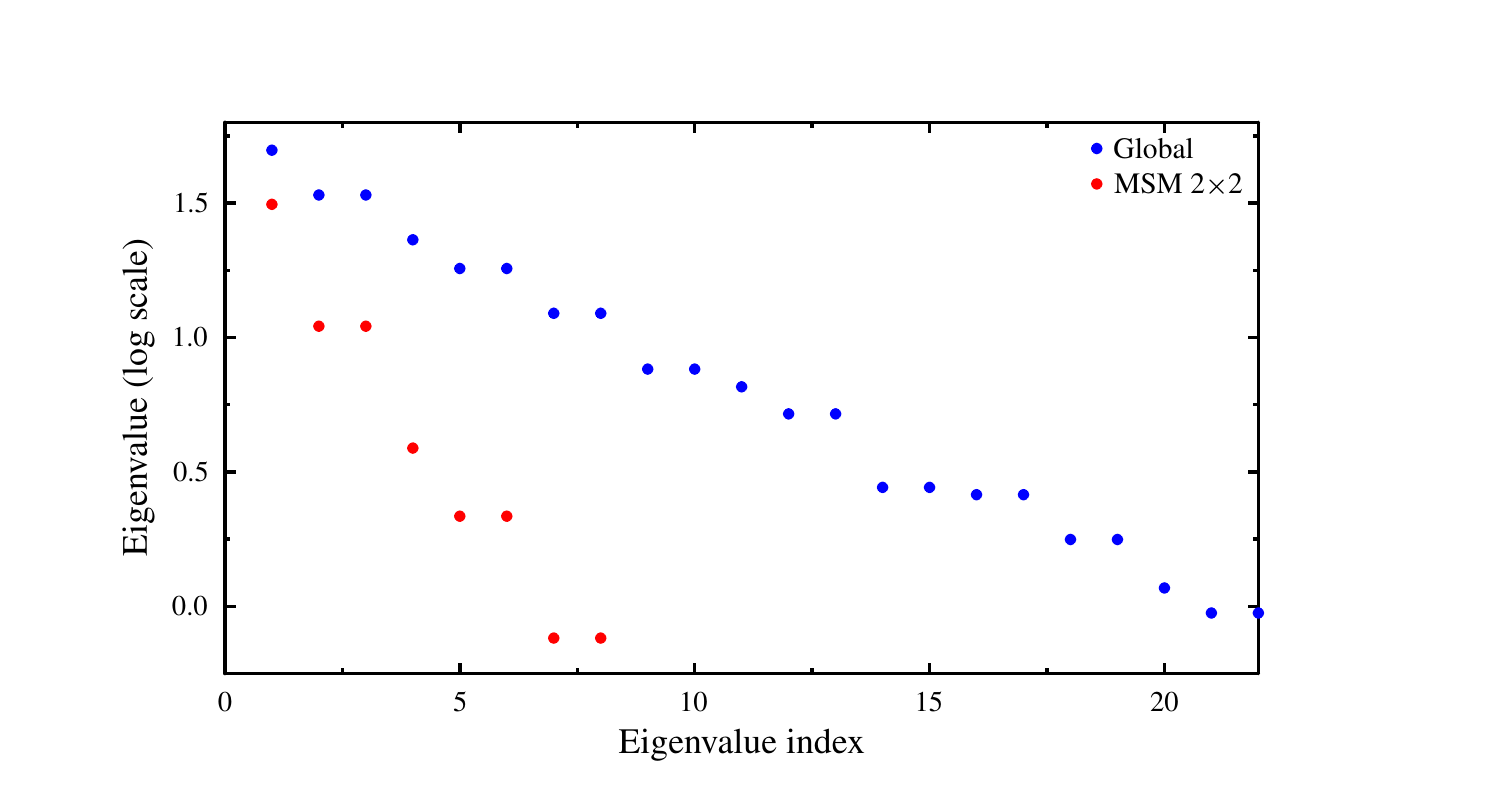}
	\caption{Decay of eigenvalues for the global and multiscale samplings in the first example.}
	\label{eigen_16x16_ch3}
\end{figure}

\begin{figure}[H]
	\centering
	\includegraphics[width= 2.2in]{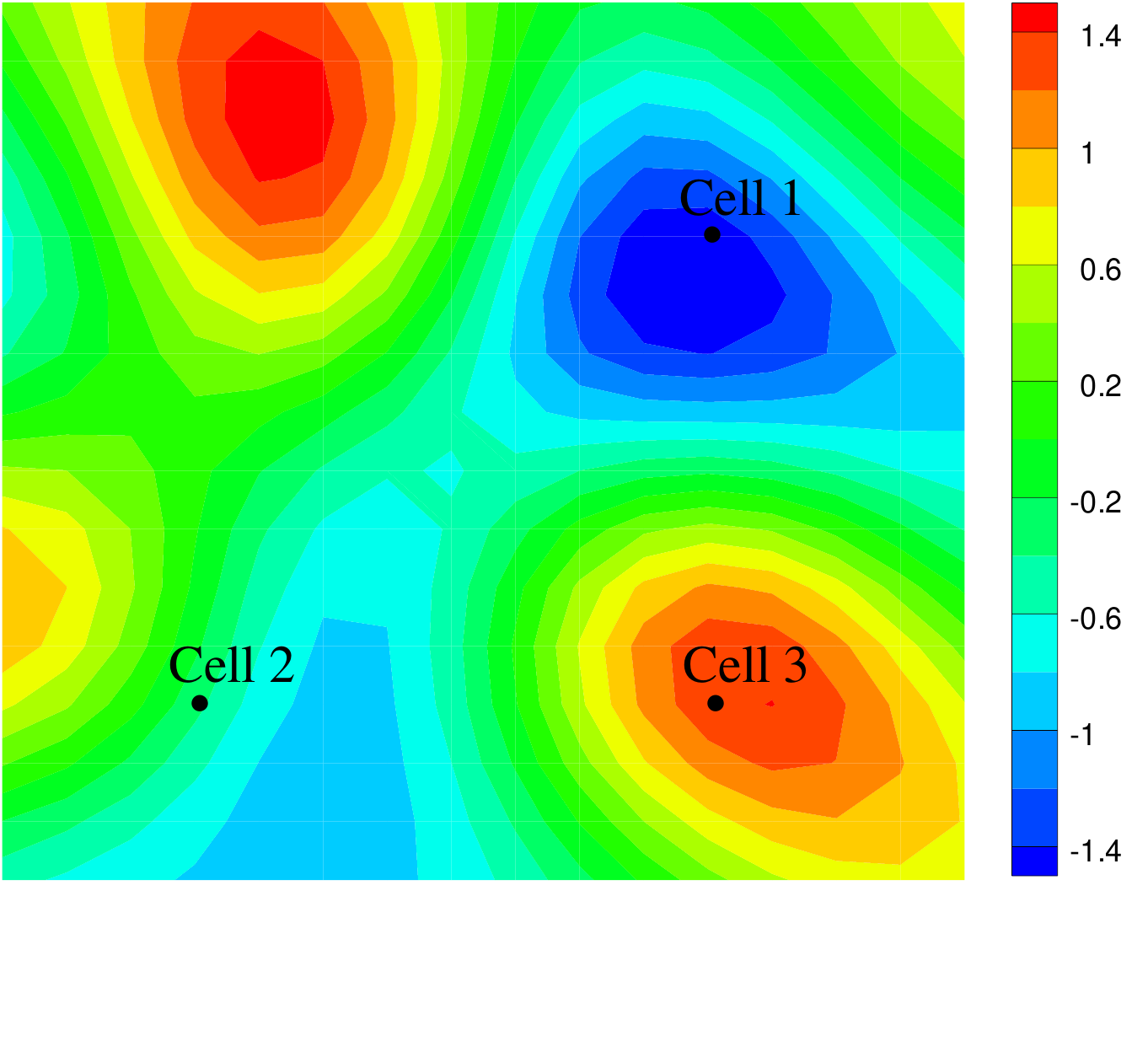}
	\hspace{5mm}
	\includegraphics[width= 2.2in]{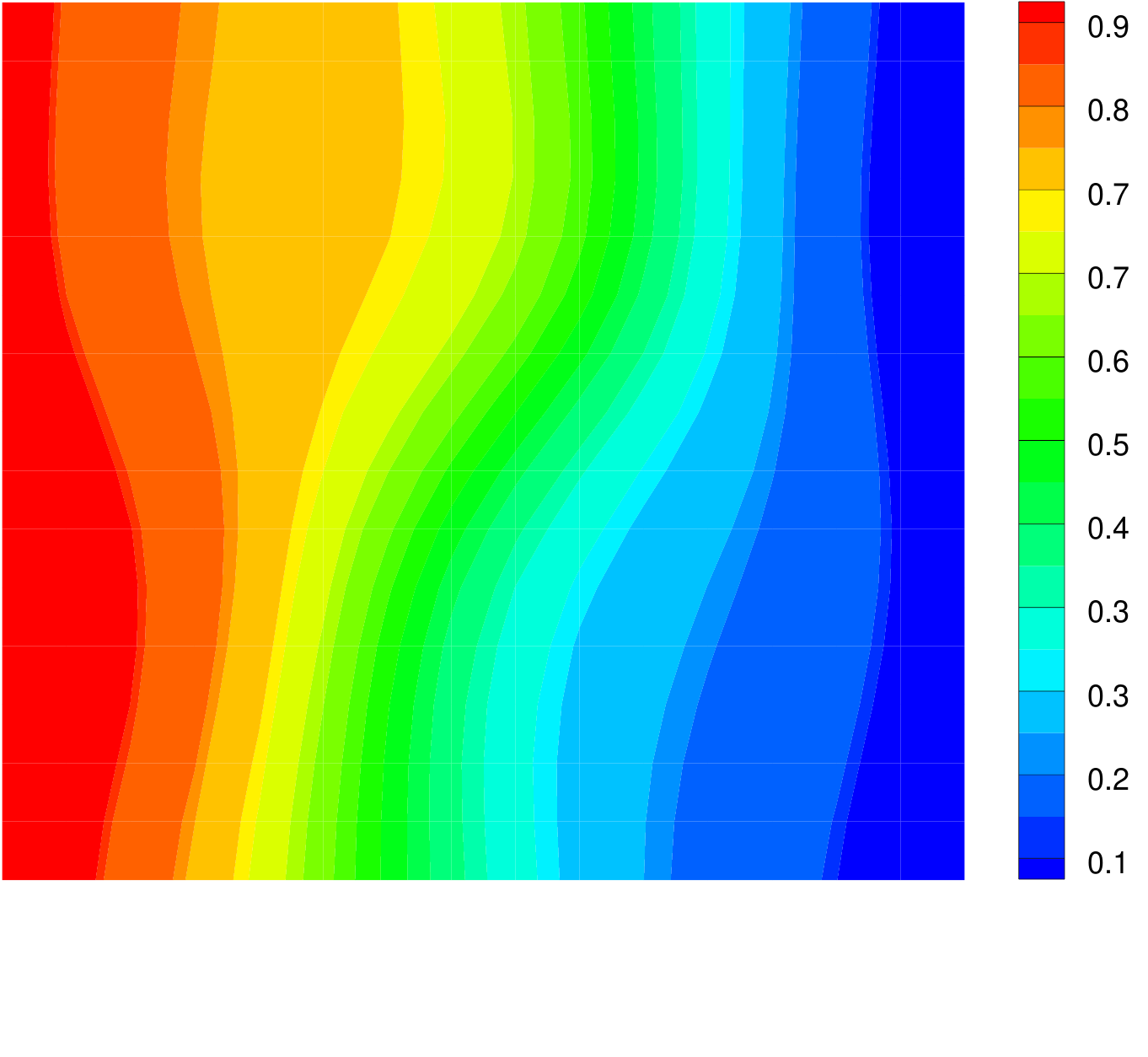}
	\caption{Reference log permeability field (left) and the corresponding reference pressure field (right) for the first example.}
	\label{ref_perm_1}
\end{figure}

In the first example, we consider $L_x=L_y=0.2$ in Eq. \eqref{kle_3}. We then generate KLEs for the global and MSM $2\times 2$ samplings. In MSM $2\times 2$ sampling, we use $H=0.5$. Figure  \ref{eigen_16x16_ch3} illustrates the decay of the eigenvalues (in log scale) for both samplings. Note that the relationship between the eigenvalues in the global sampling and the eigenvalues in multiscale sampling can be obtained directly by a change of variables in Eq. \eqref{kle_efun_3}. 
We take the first $20$ eigenvalues that preserve more than $97\%$ of the total energy for the global sampling. Five eigenvalues are used for each subdomain in the multiscale sampling. We generate a reference synthetic permeability field on a computational fine mesh of size $16\times 16$, and then run our numerical simulator to generate the corresponding reference pressure field. We run  the MCMC algorithms conditioned on this pressure field. Figure  \ref{ref_perm_1} shows these reference fields. Furthermore, we use a coarse mesh of size $8\times 8$ as a filtering step in the preconditioned MCMC. We let the local blocking number $N_{lb}=1$ and $\beta = 0.5$ in Eq. \eqref{RW_sampler}.

\begin{figure}[H]
	\centering
	\includegraphics[scale = 0.55]{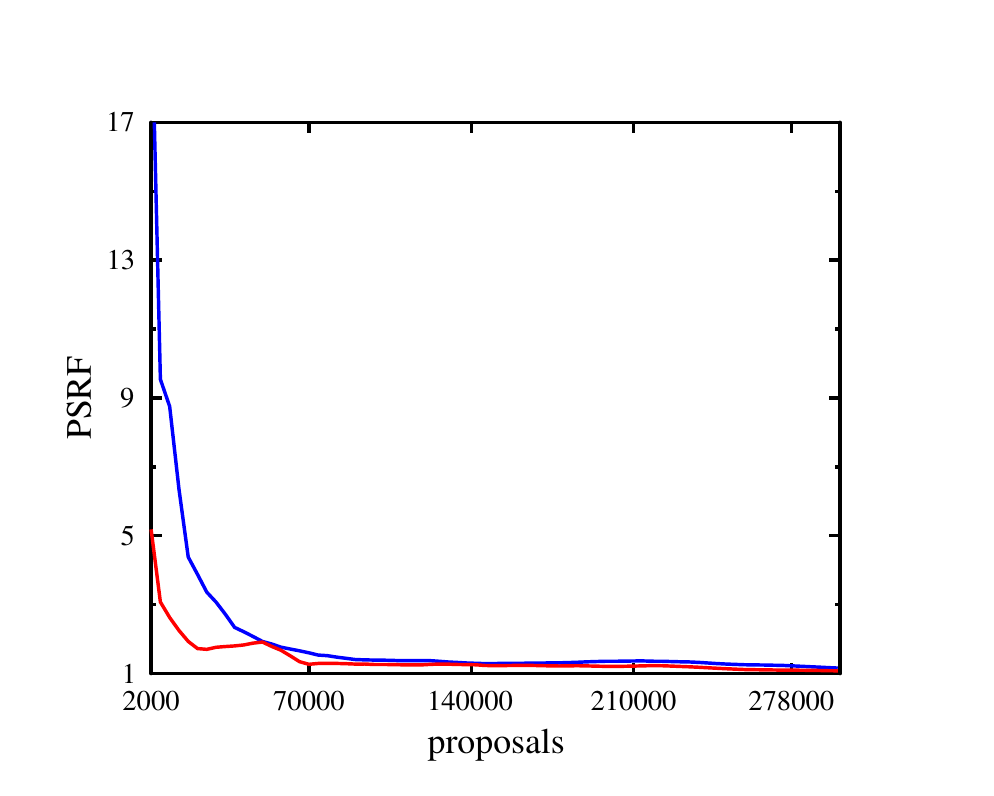}
	\includegraphics[scale= 0.55]{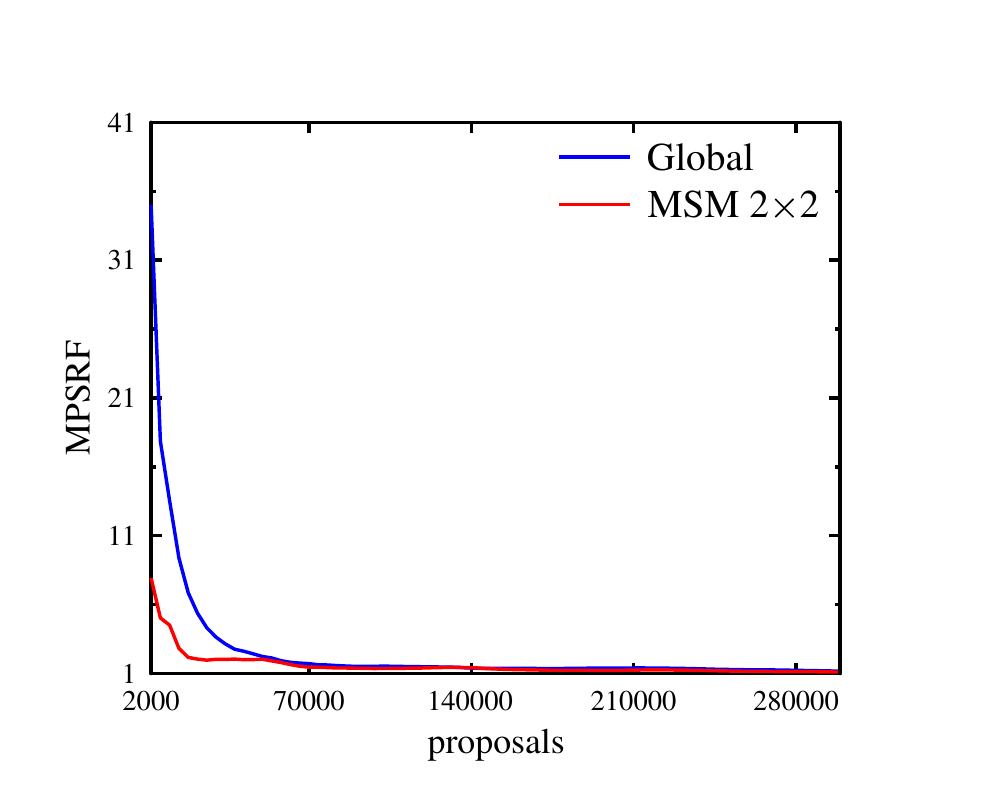}\\
	\includegraphics[scale = 0.55]{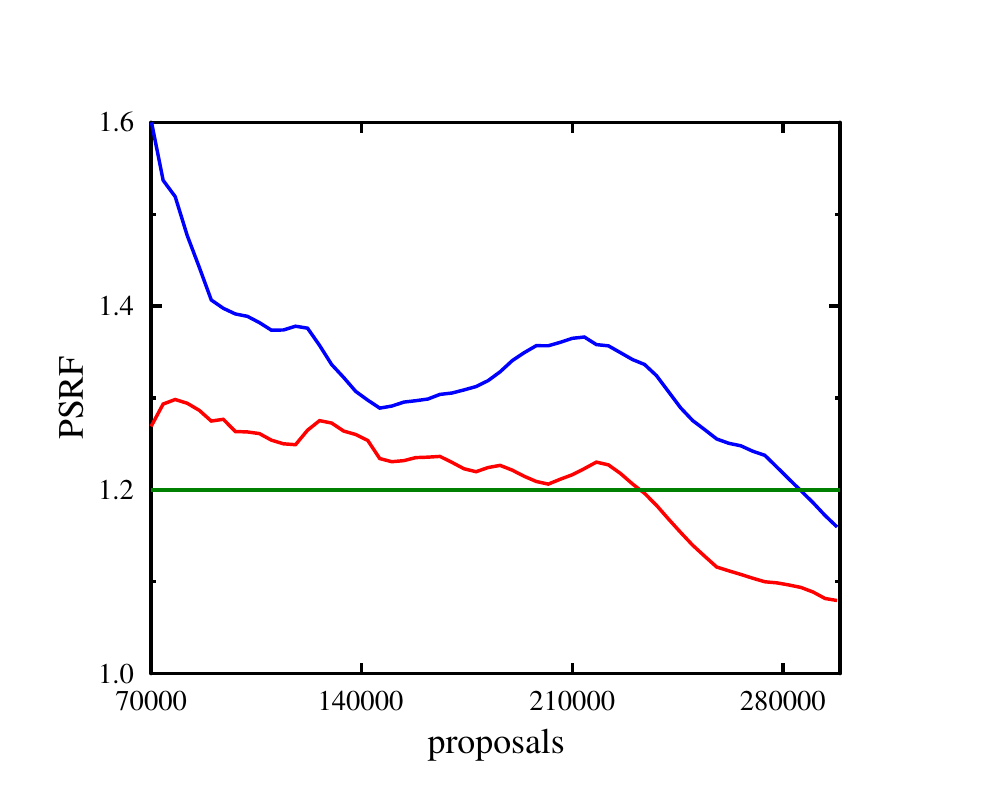}
	\includegraphics[scale= 0.55]{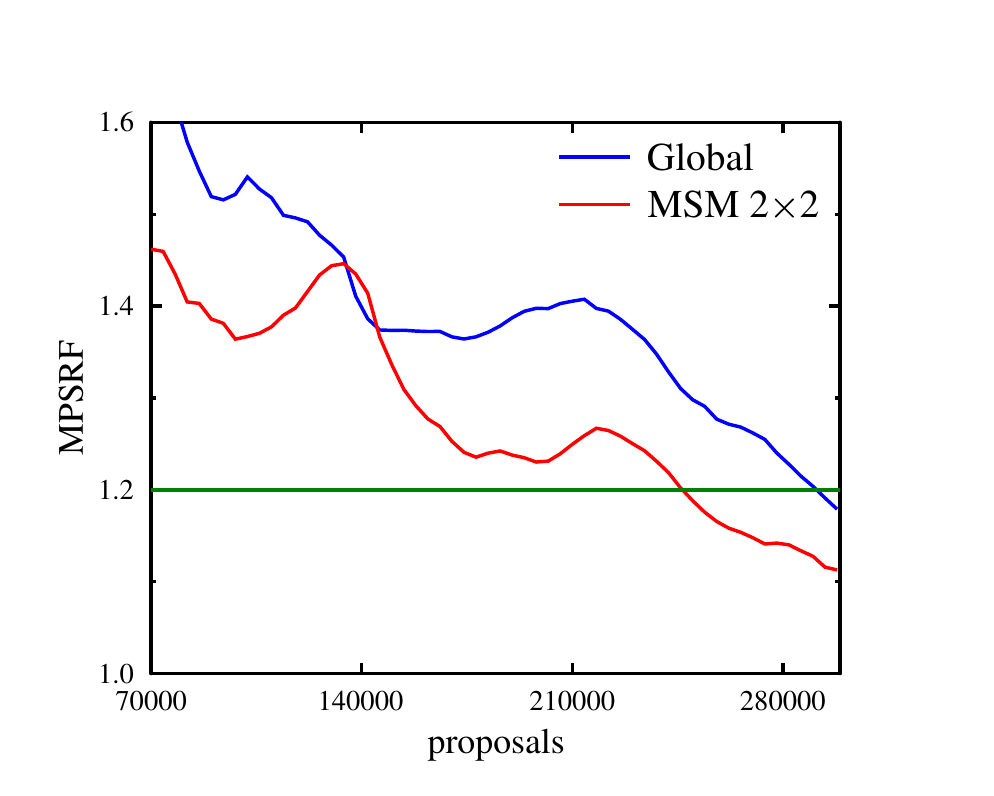}
	\caption{Top: The maximum of PSRFs and MPSRF for the MCMC method with and without multiscale sampling in the first example.
	Bottom: Tails of the maximum of PSRFs and MPSRF curves.}
	\label{MPSRF_16x16_ch3}
\end{figure}

\begin{figure}[H]
	\centering
	\includegraphics[scale = 0.9]{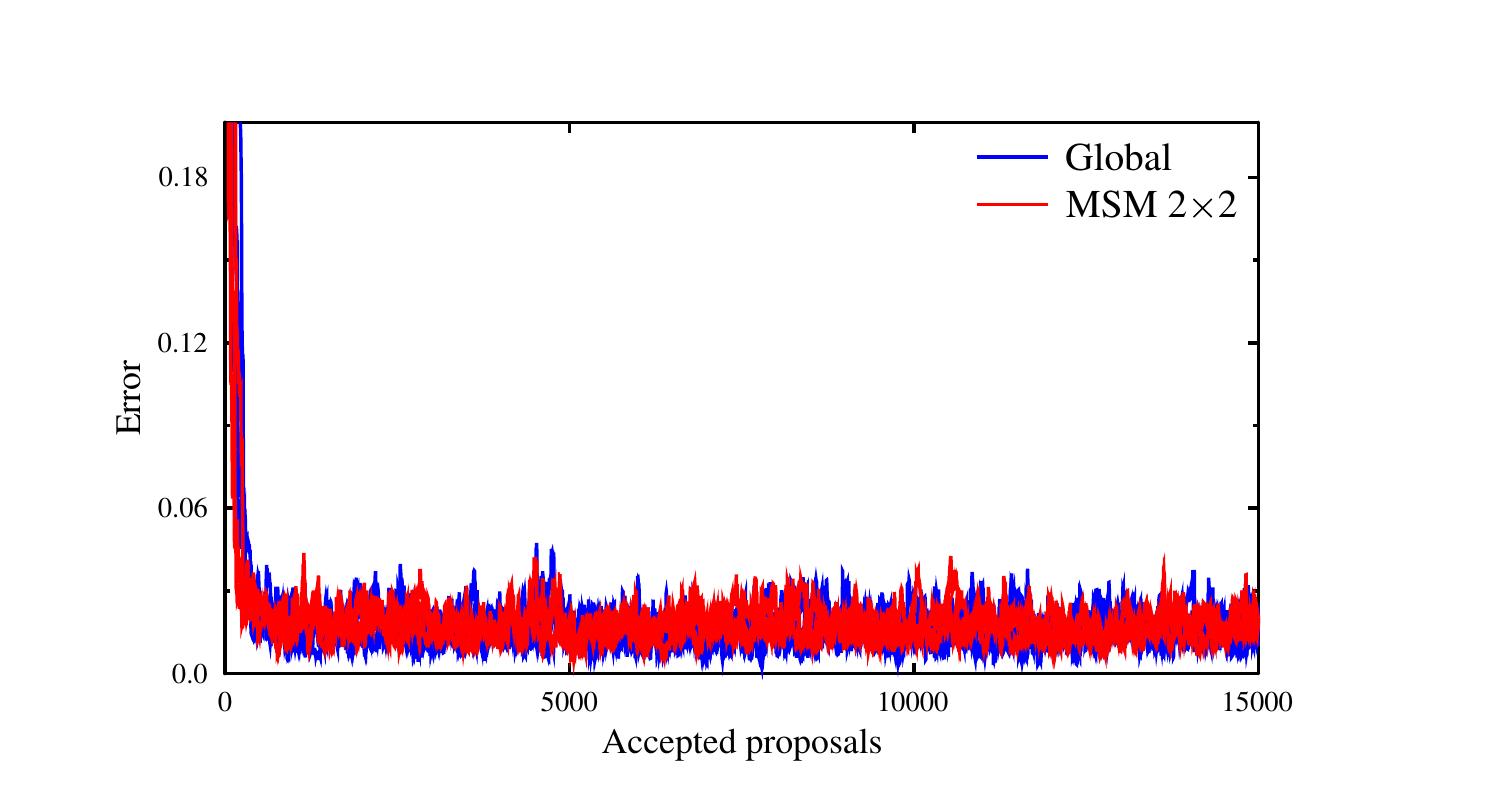}
	\caption{Error curves of the preconditioned MCMC with and without multiscale sampling for the first example.}
	\label{error_16x16_ch3}
\end{figure}

As we discussed in subsection \ref{mpsrf_3}, we analyze the convergence of the MCMCs using PSRFs and MPSRF. An MCMC method converges to the stationary distribution if both MPSRF and the maximum of PSRFs get closer to 1. In \cite{Brian2007} the author considered a value of $1.2$ for these parameters to confirm the convergence of the chains. In line with that, we decide to stop the simulation once these parameters reach $1.2$. Figure  \ref{MPSRF_16x16_ch3} shows that the preconditioned MCMC methods with and without multiscale sampling converge. However, the plots at the bottom in Figure  \ref{MPSRF_16x16_ch3} show that the preconditioned MCMC method with multiscale sampling converges to the stationary distribution earlier than the method without multiscale sampling. Table \ref{Lx_Ly_1} shows the acceptance rates for both methods as well as the precision parameters for coarse- and fine-grid simulations. The acceptance rate increases slightly when we use MSM. The errors between the reference and simulated pressure data, which are used in the likelihood function, for both methods are shown in Figure  \ref{error_16x16_ch3}. Both methods produce similar error curves.

\begin{table}[H] 
	\caption{A comparison of acceptance rates for the MCMC with and without MSM for the first example.}
	\center
	\begin{tabular}{|cccc|}
		\hline
		&  \quad  MCMC with global sampling & \quad MCMC with multiscale sampling&  \\
		\hline
		$\sigma_F^2$ &   $10^{-3}$  & $10^{-3}$  & \\               
		$\sigma_C^2$    & $5\times10^{-3}$ & $5\times10^{-3}$&   \\  
		acc. rate                     & $53\%$   &$55\%$&  \\
		\hline
	\end{tabular}
	\label{Lx_Ly_1}     
\end{table}

\begin{figure}
	\centering
	\begin{tabular}[b]{c}
		\includegraphics[scale=.4]{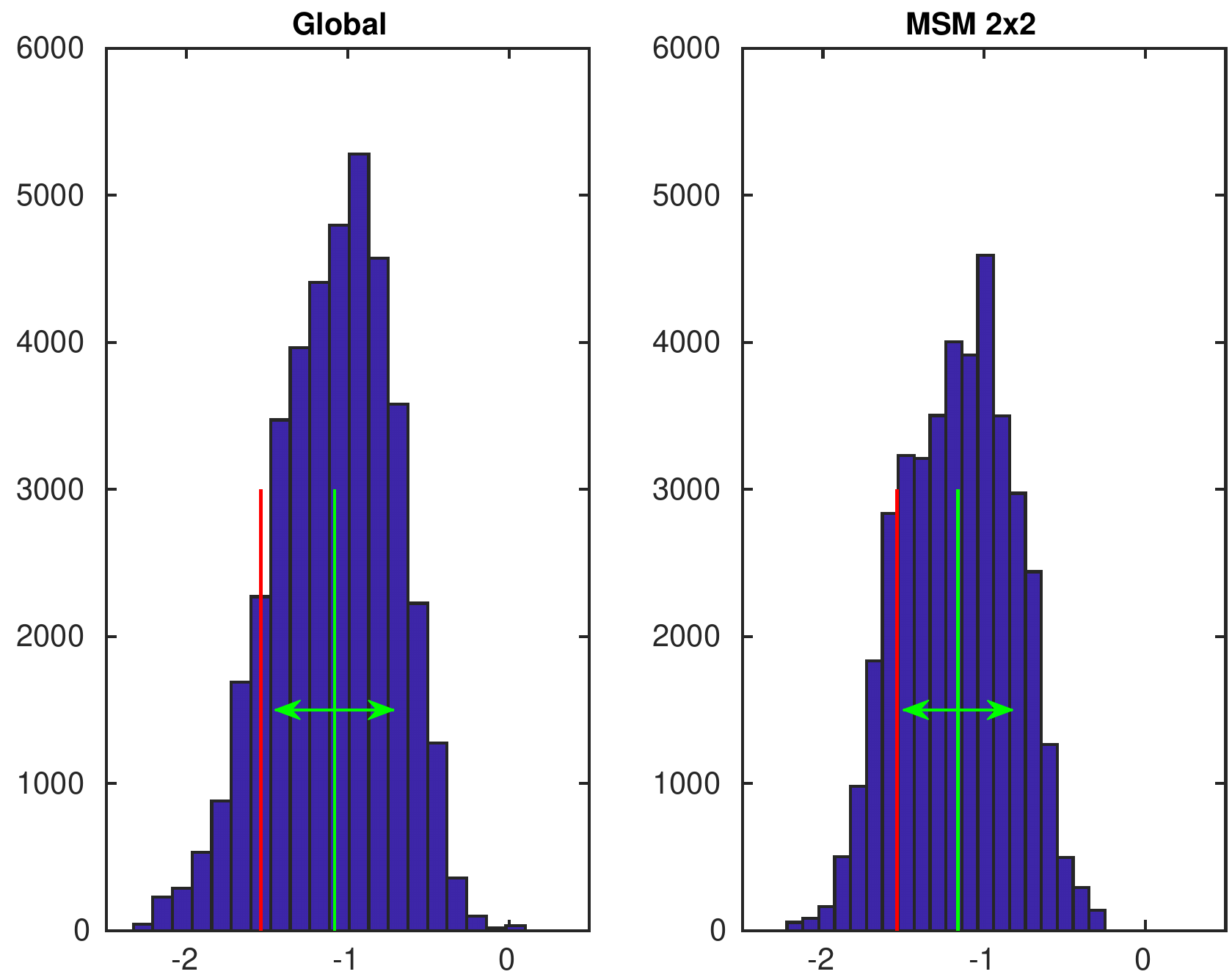}\\
		Cell 1
	\end{tabular}
	\begin{tabular}[b]{c}
		\includegraphics[scale=.4]{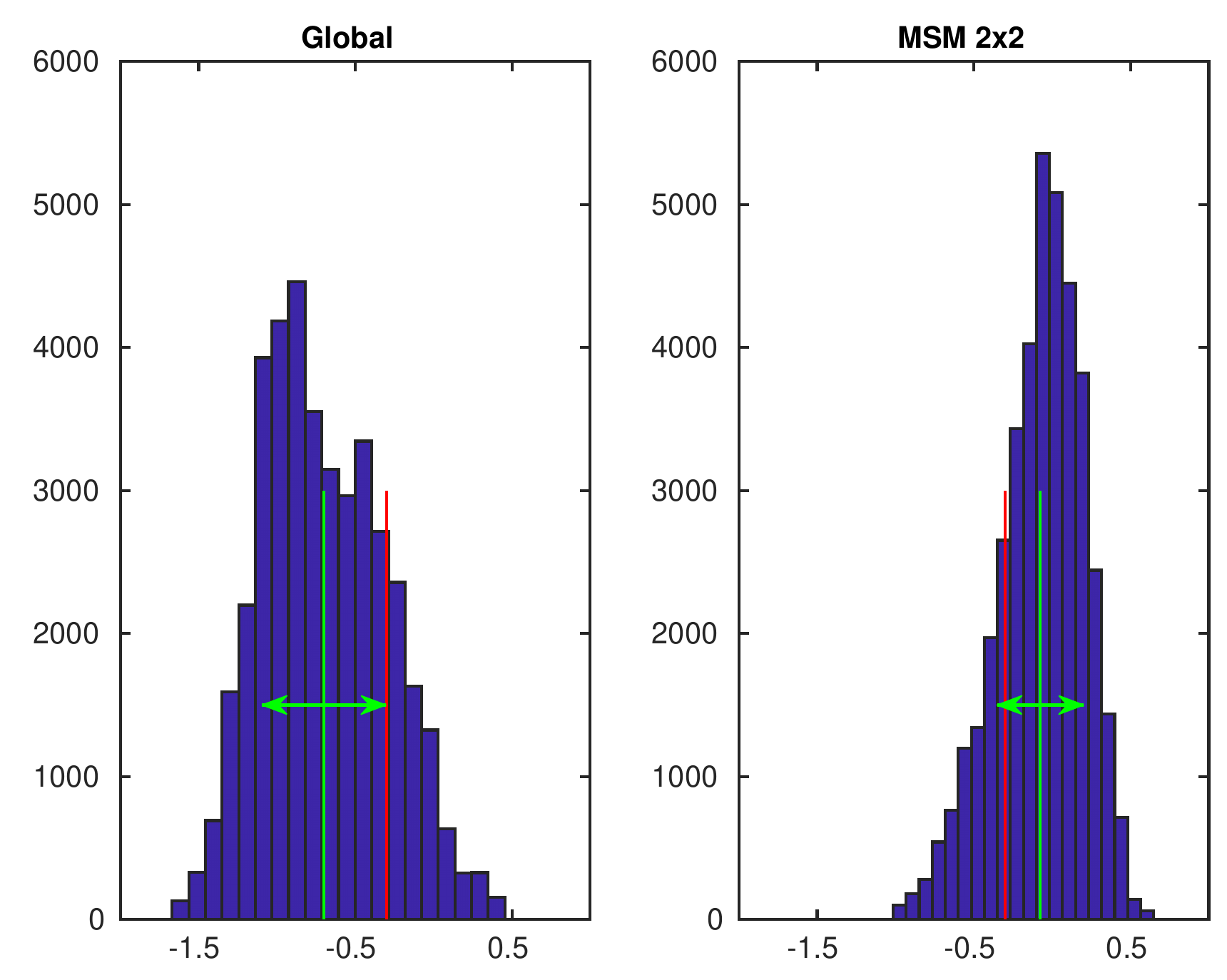}\\
		Cell 2
	\end{tabular}
	\begin{tabular}[b]{c}
		\includegraphics[scale=.4]{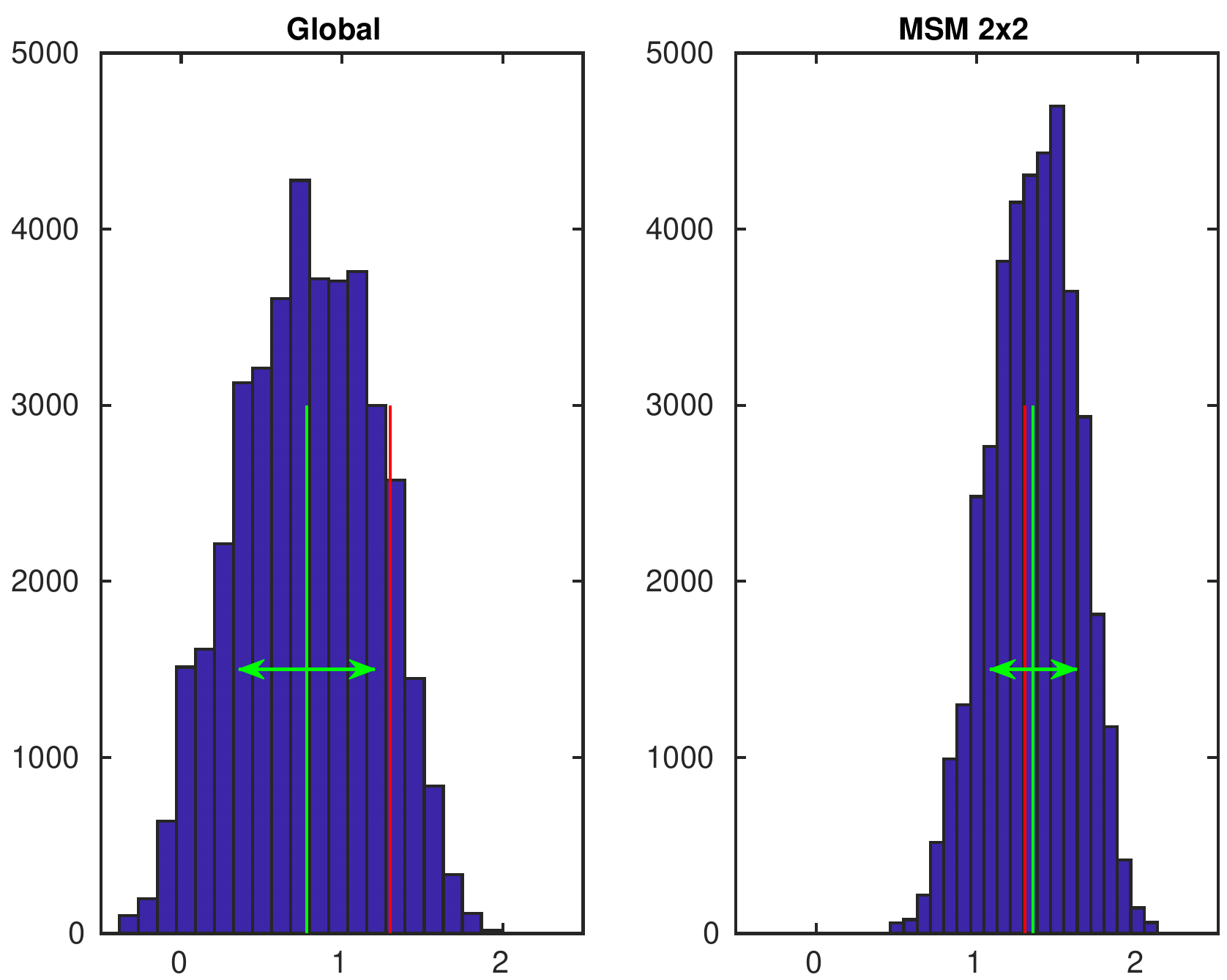}\\
		Cell 3
	\end{tabular}
	\vspace*{-0.2cm}
	\caption{Posterior histograms for three selected cells.}
	\label{posterior}
\end{figure}

After the convergence of both MCMC methods, we take 10000 log permeability values from each chain and draw the posterior histograms for three cells with high, medium and low permeability values in the computational domain. See Figure  \ref{ref_perm_1} for those three cells. Figure  \ref{posterior} shows the posterior histograms with the true (red vertical line) and mean (green vertical line) values of the log permeability for the cells.
Table \ref{cells} shows these values and the corresponding standard deviation. We observe in Figure  \ref{posterior} that when we use the multiscale sampling method, the mean values are within one standard deviation (green horizontal line) for all the three cells. Also, in MSM, the mean of the posterior histogram is almost the same as the true value for the cell 3. We do not observe a similar behavior in the posterior histograms in the global sampling method.

\begin{center}
	\begin{table}[H]
		\caption{True and mean values of log permeability with the corresponding standard deviation for three cells.}\vspace{0.1cm}
		\setlength{\tabcolsep}{17pt}
		\centering
		\scalebox{0.85}{
		\begin{tabular}{|l c c c c c c |}
			\hline 
			& \multicolumn{2}{c}{Cell 1} & \multicolumn{2}{c}{Cell 2} & \multicolumn{2}{c|}{{Cell 3}}  \\
			& Global & MSM & Global & MSM & Global & MSM \\ 
			\hline
			\\\\[-3.95\medskipamount ]
			True & $-1.54$ & $ - $       & $-0.3$  & $ - $        & $1.3$ &$  - $  \\
			Mean & $-1.08$ & $-1.16$ & $-0.69$ & $-0.07$ & $0.78$ & $1.35$  \\
			SD & $ 0.36$ & $0.34 $ & $0.39$ & $0.27$ & $0.42$ & $0.26$ \\
			\hline 
		\end{tabular}}
		\label{cells}
	\end{table}
\end{center} 

\begin{figure}[H]
	\centering
		\includegraphics[width= 1.5in]{figures/ex_1_perms/ref_perm_16x16_crop}\\
\includegraphics[width=1.5in]{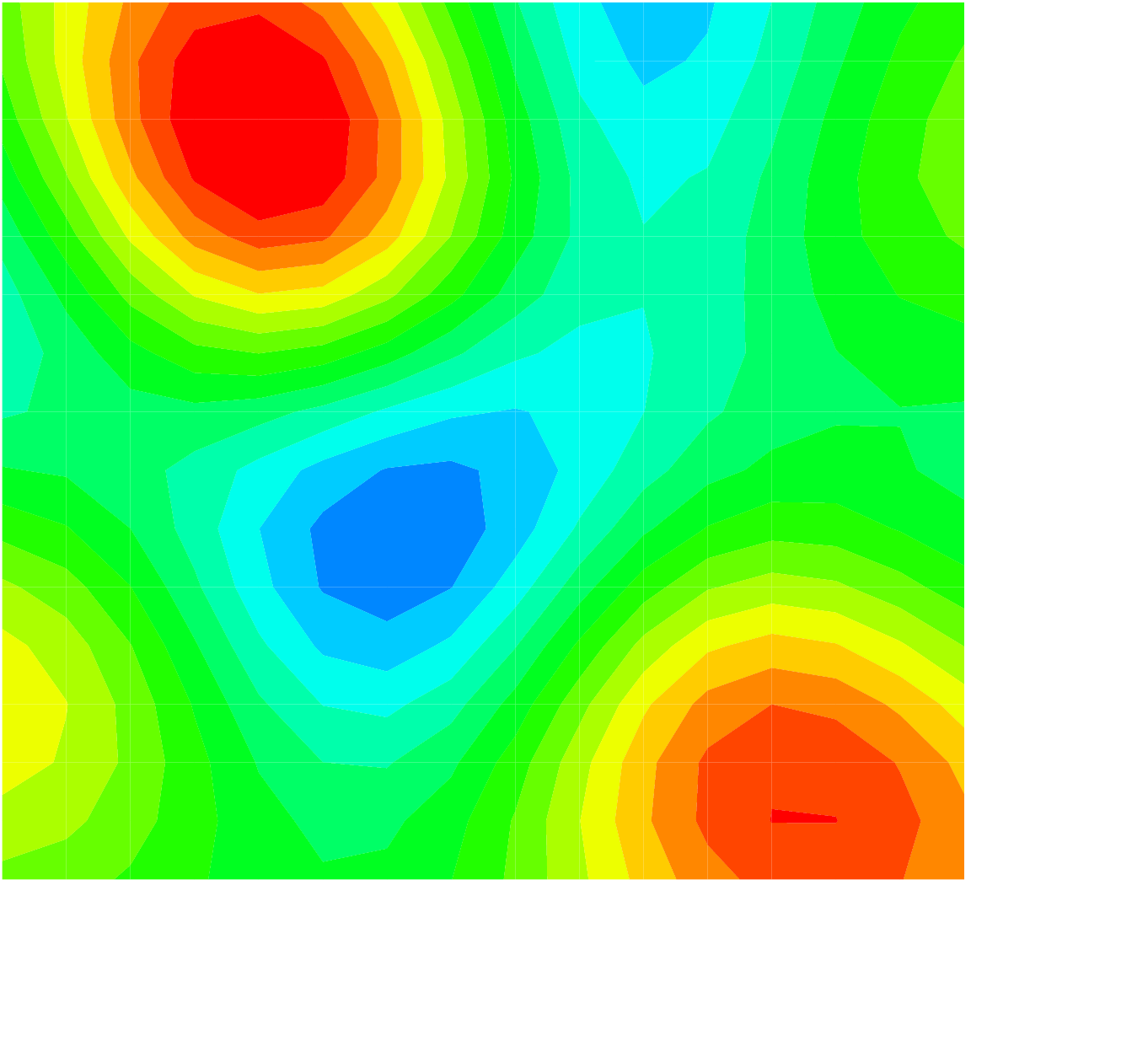}
\includegraphics[width=1.5in]{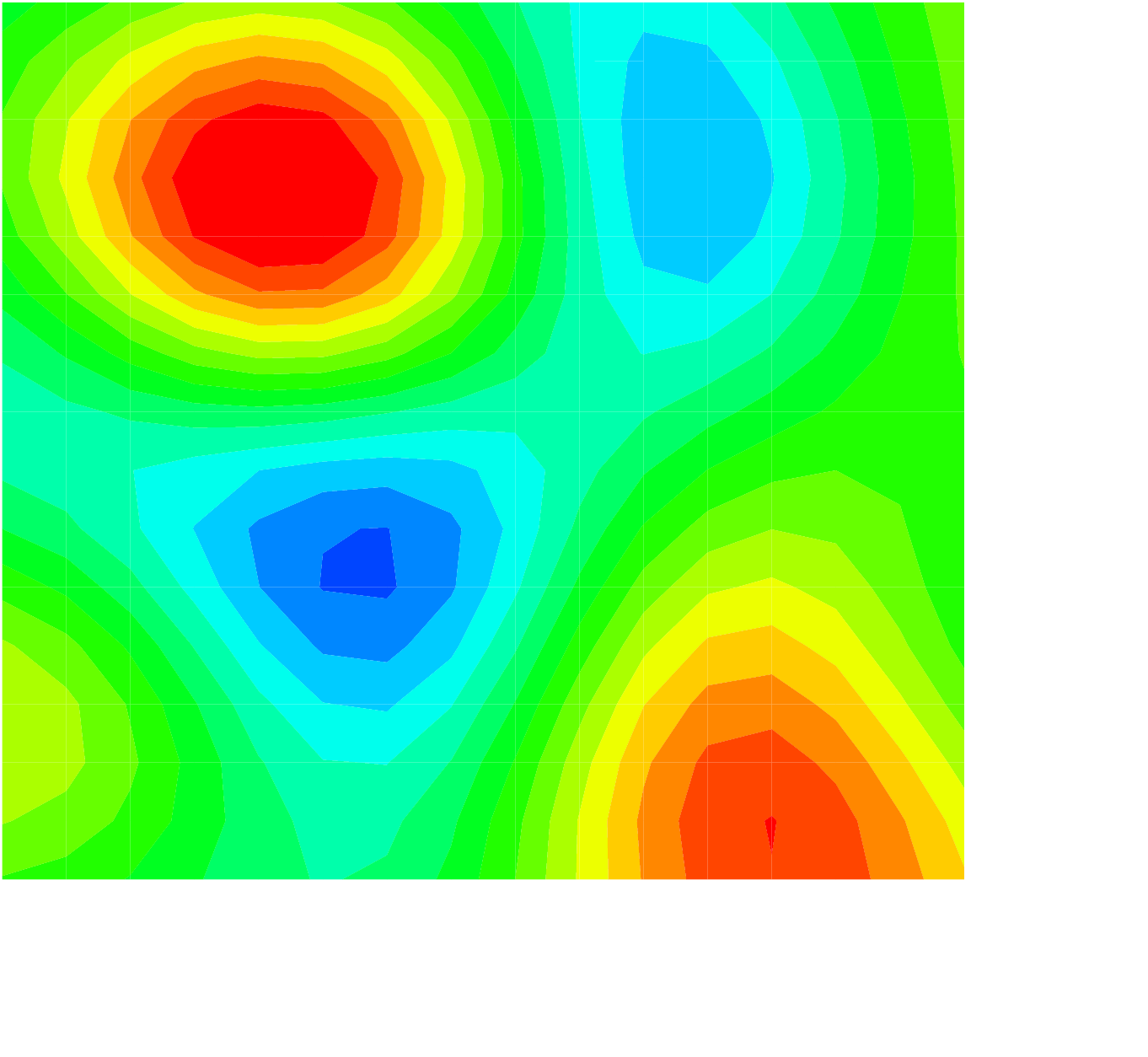}
\includegraphics[width=1.5in]{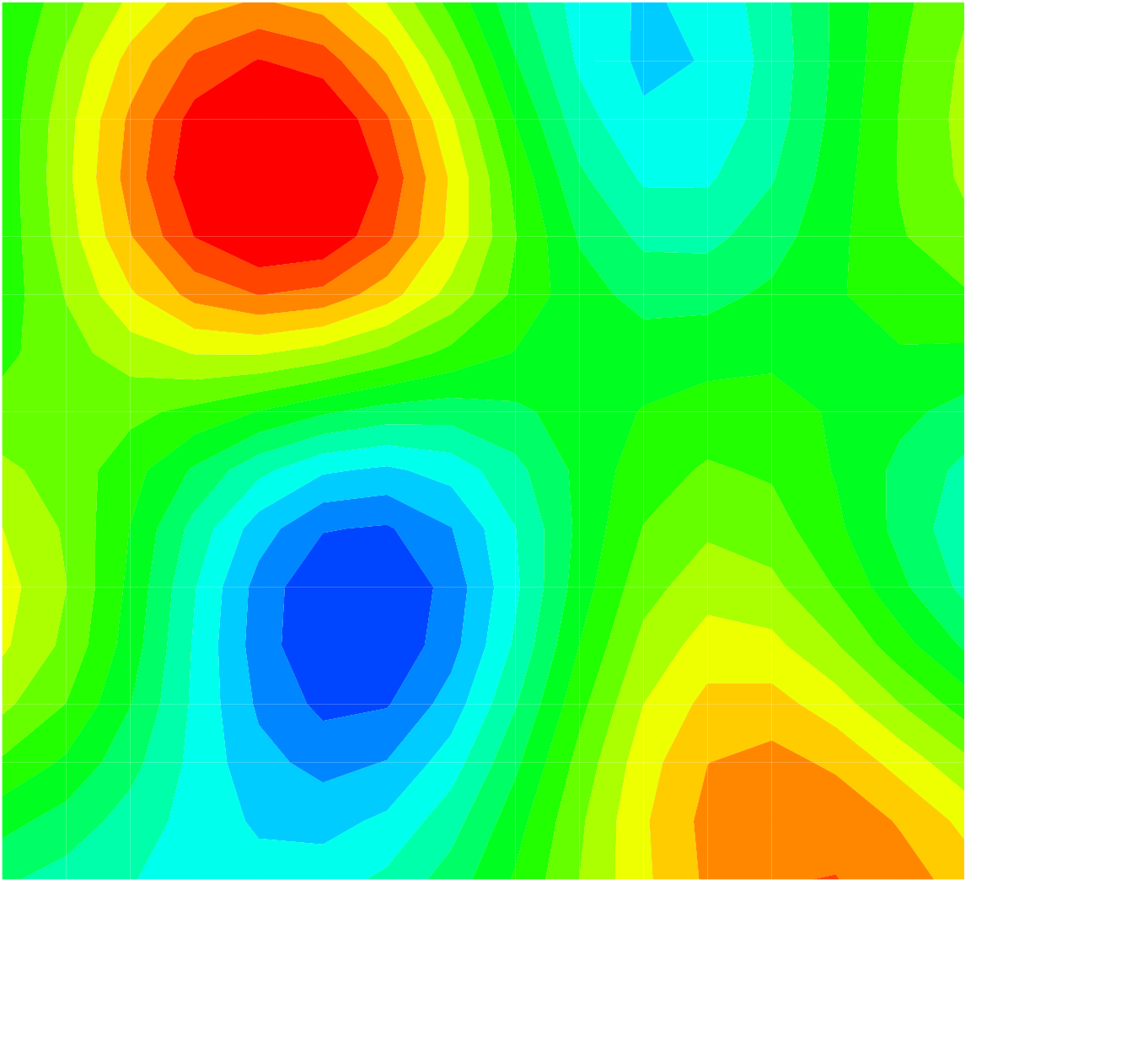}
\\
	\includegraphics[width=1.5in]{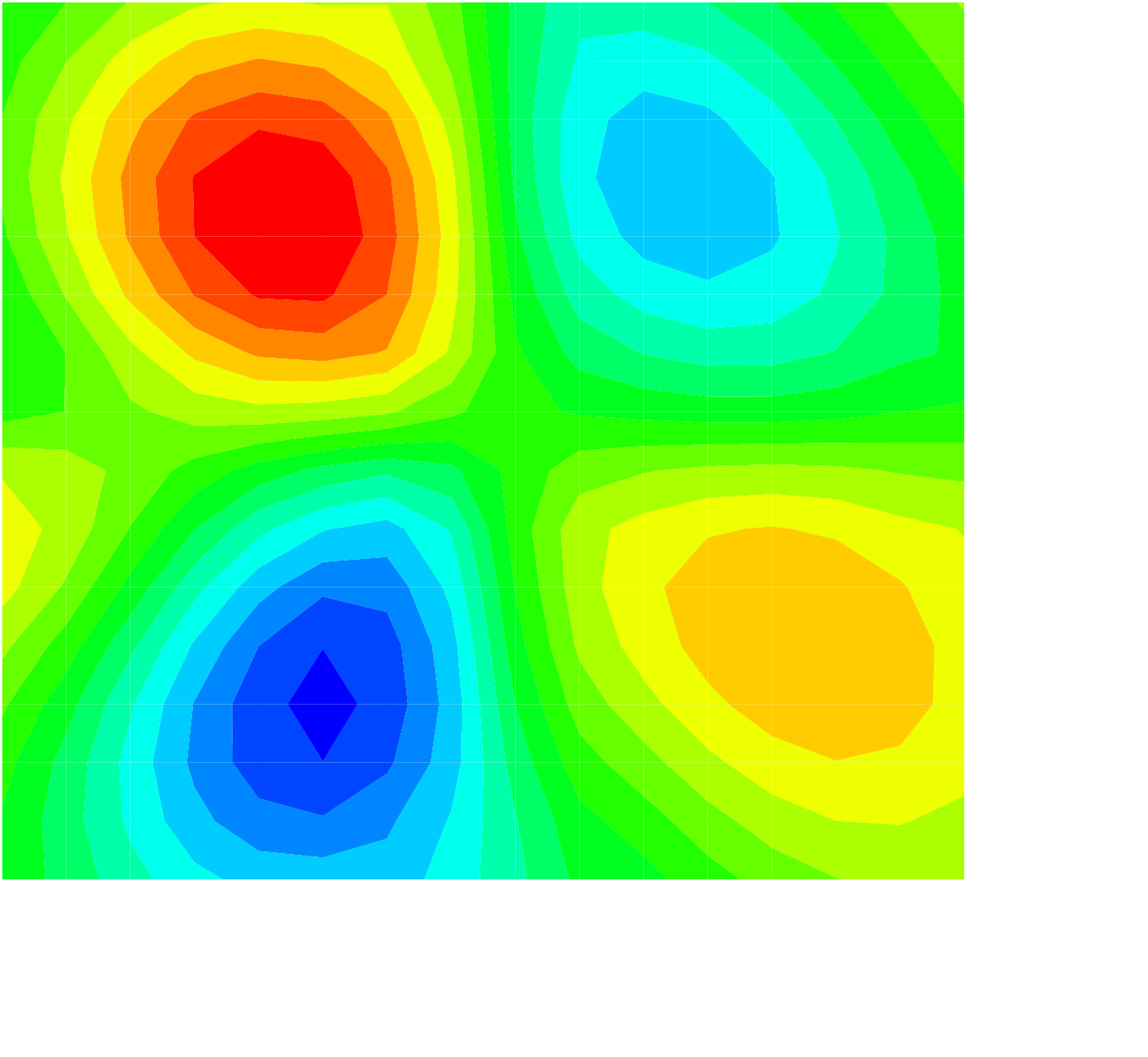}
	\includegraphics[width=1.5in]{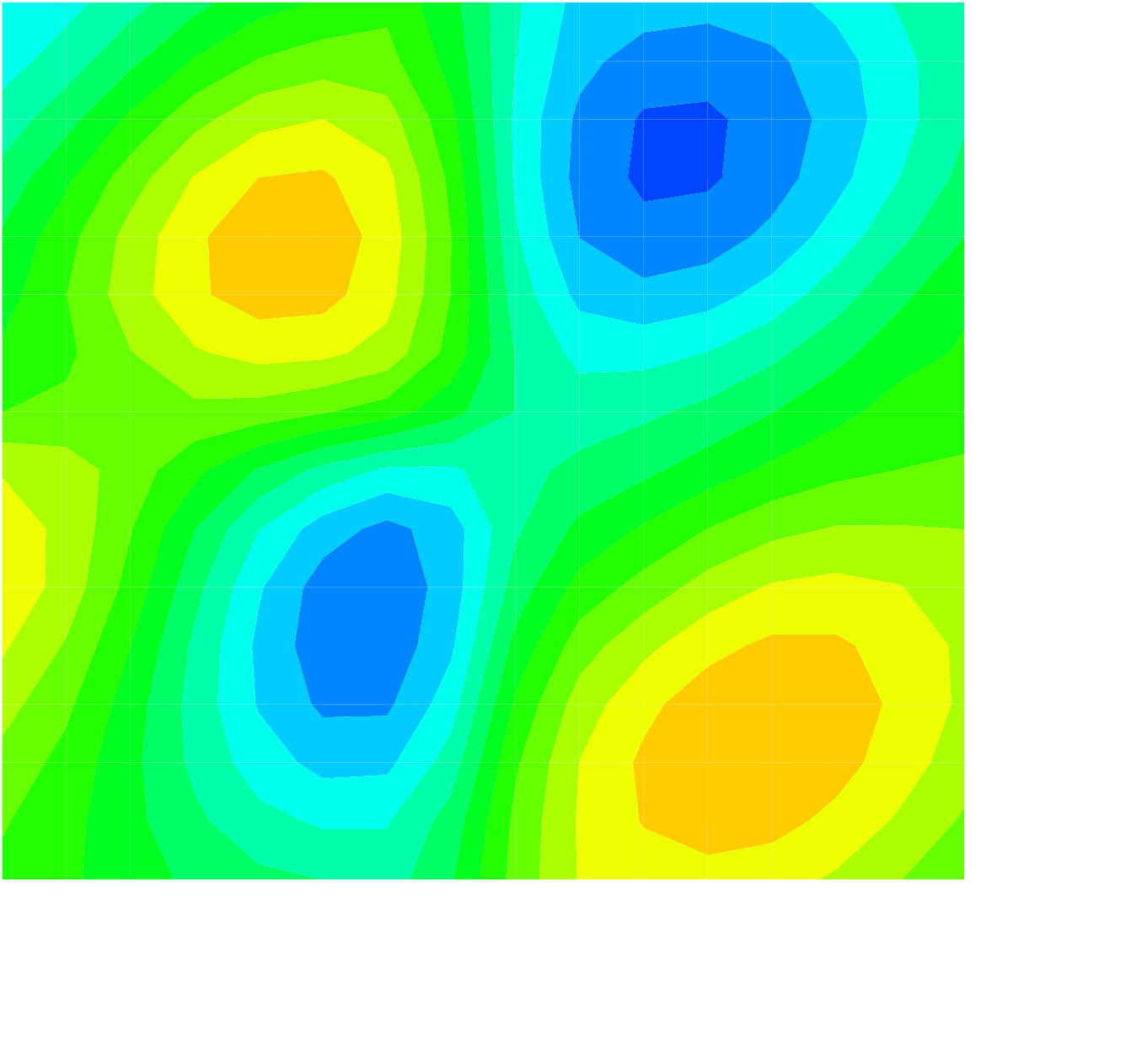}
	\includegraphics[width=1.5in]{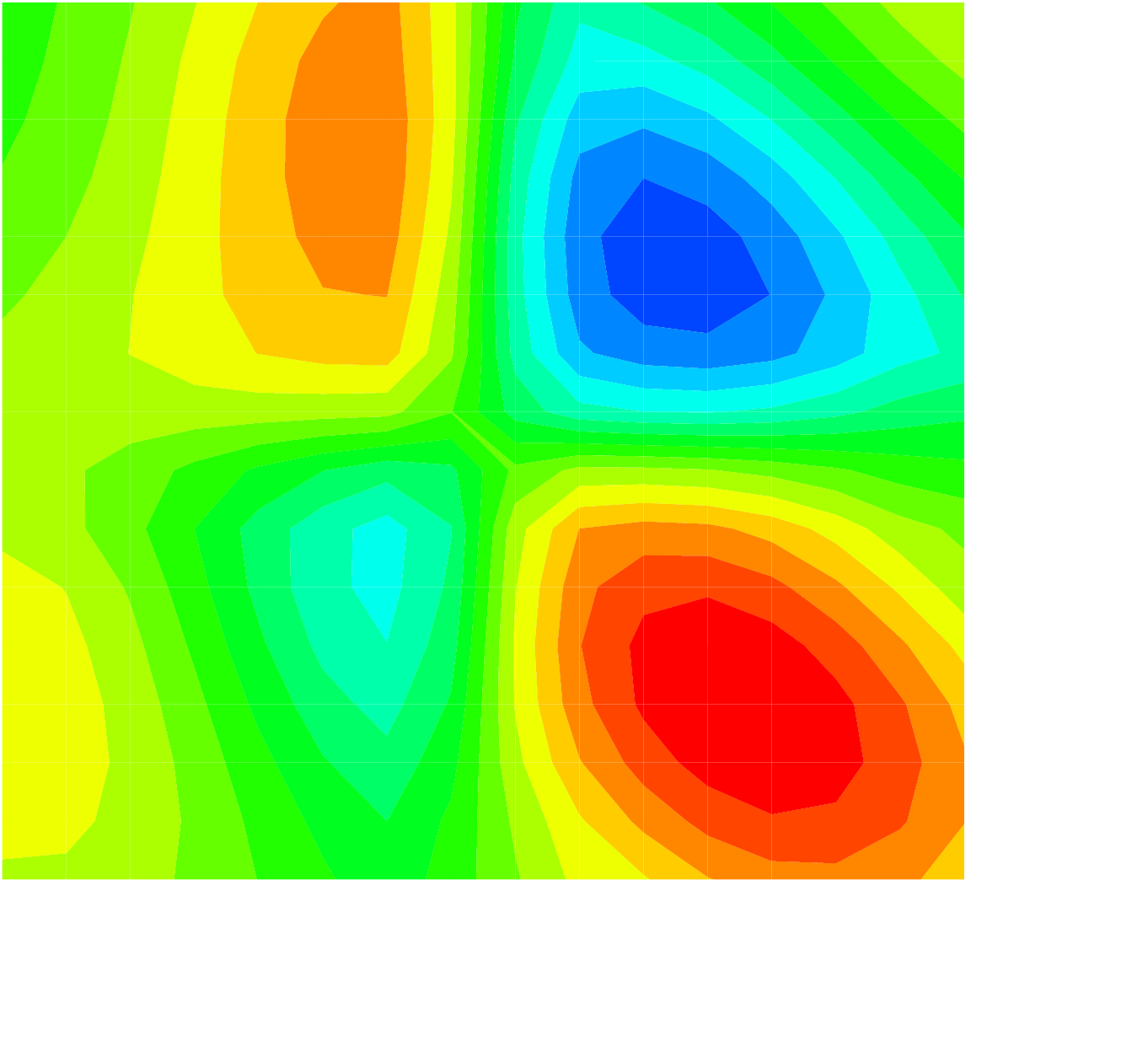}
\\
	\caption{First row: Reference log permeability filed. Second row: Accepted permeability fields in the global sampling method. Third row: Accepted permeability fields in MSM $2\times 2$.
	From left to right, log permeability fields at 20000, 40000 and 60000 iterations, respectively, from chain 1 in the first example.}
	\label{perm_1_1}
\end{figure}

\begin{figure}[H]
	\centering
	\includegraphics[width= 1.5in]{figures/ex_1_perms/ref_perm_16x16_crop.pdf}\\
\includegraphics[width=1.5in]{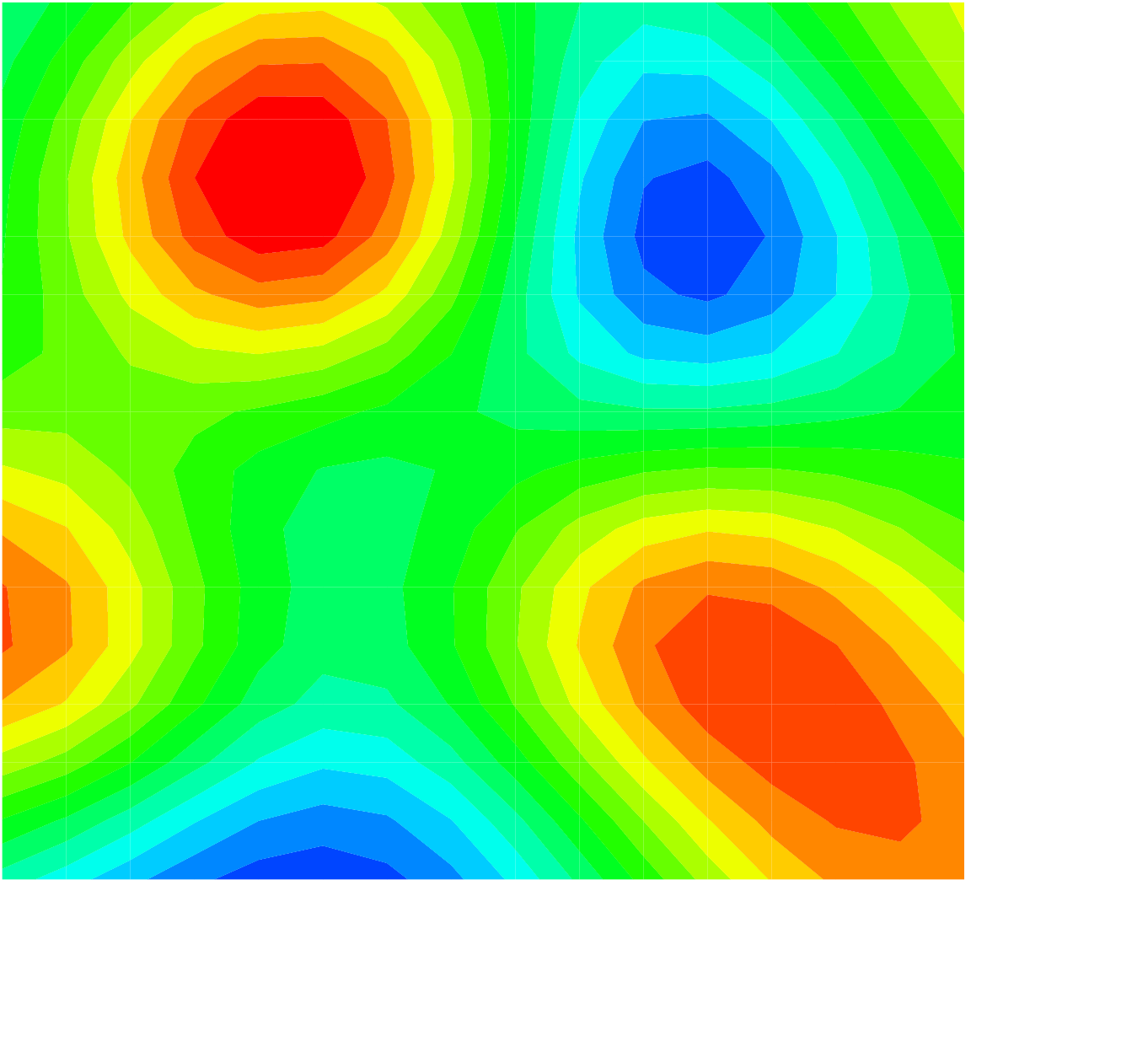}
\includegraphics[width=1.5in]{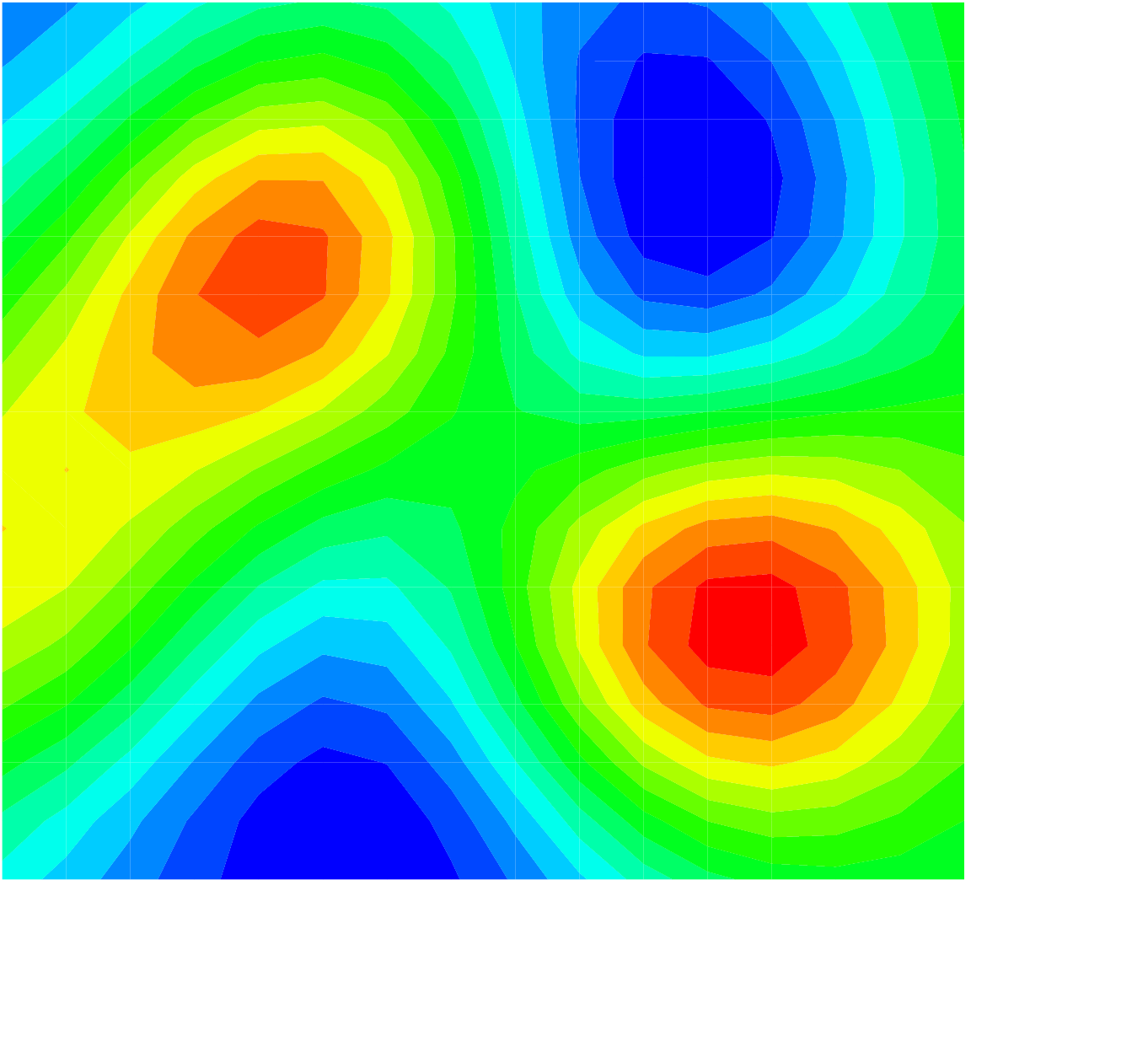}
\includegraphics[width=1.5in]{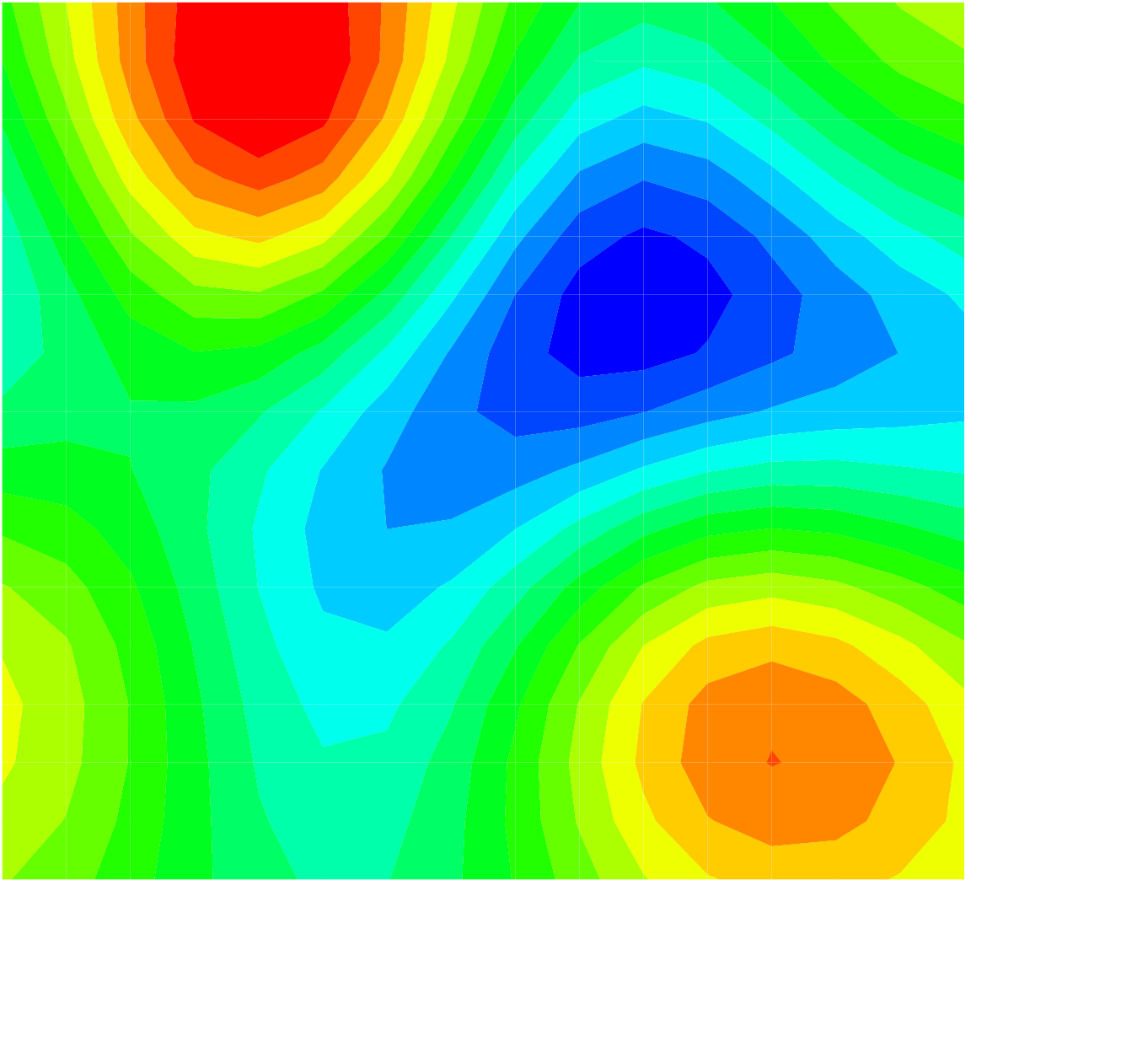}\\

\includegraphics[width=1.5in]{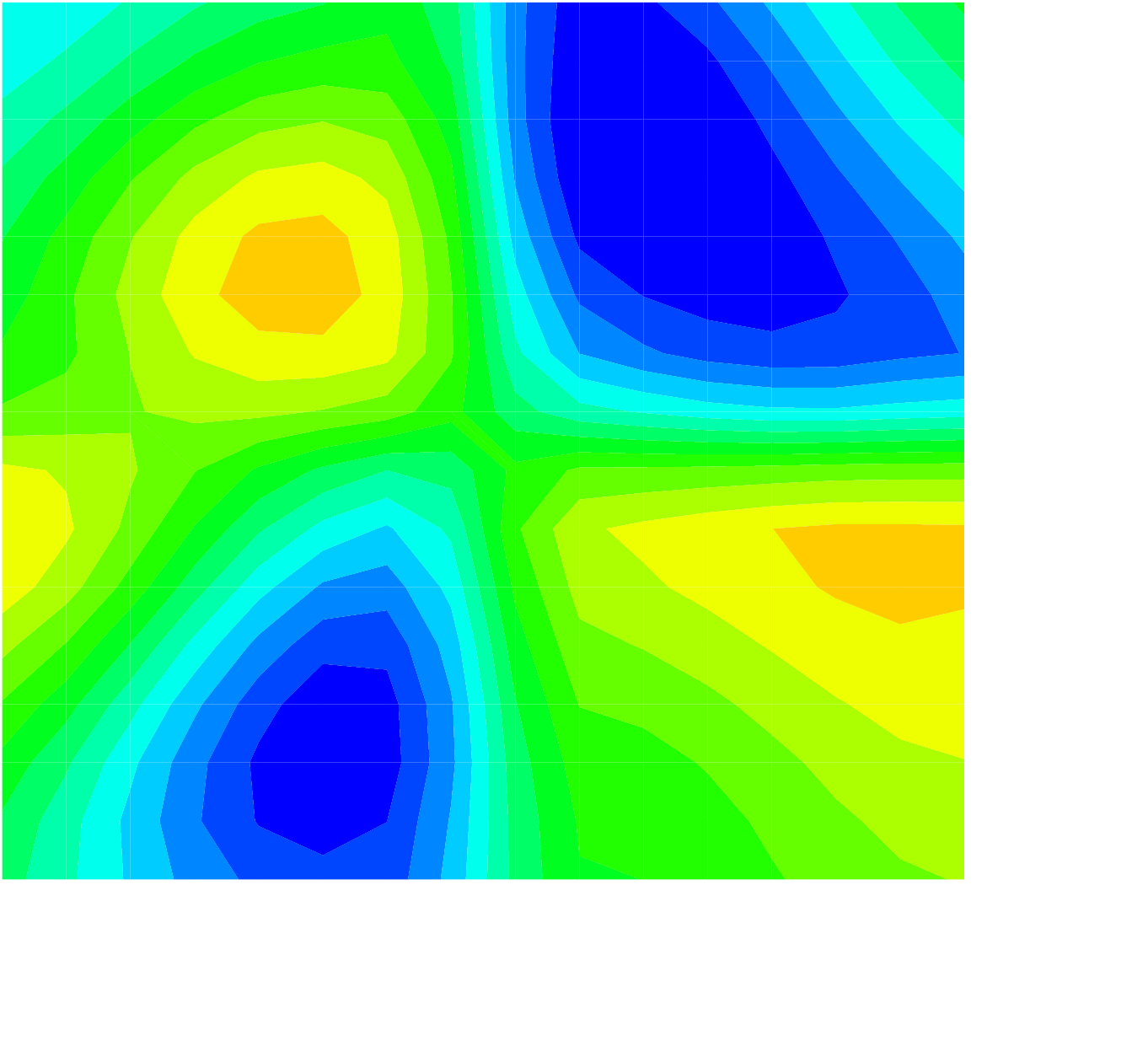}
\includegraphics[width=1.5in]{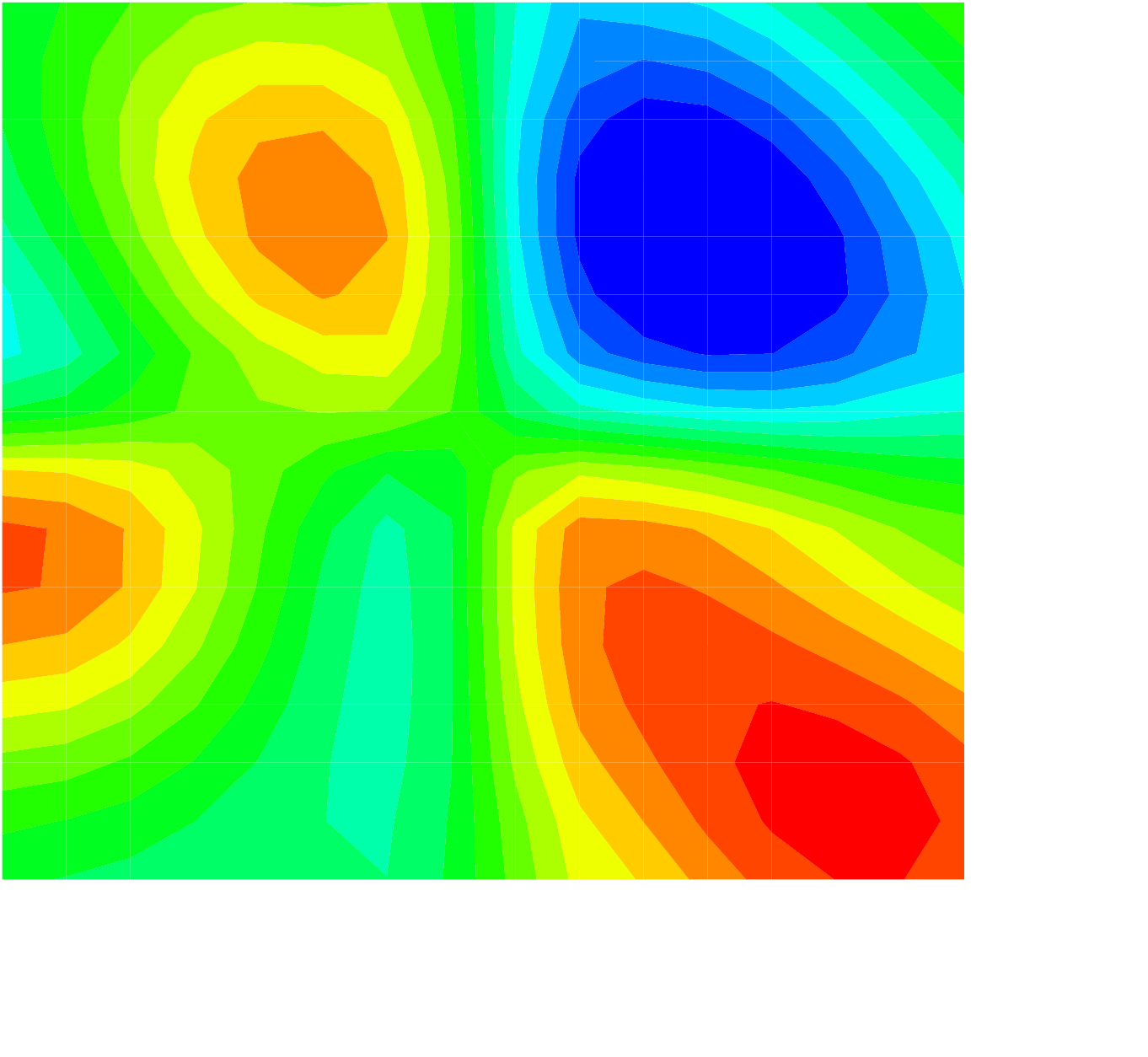}
\includegraphics[width=1.5in]{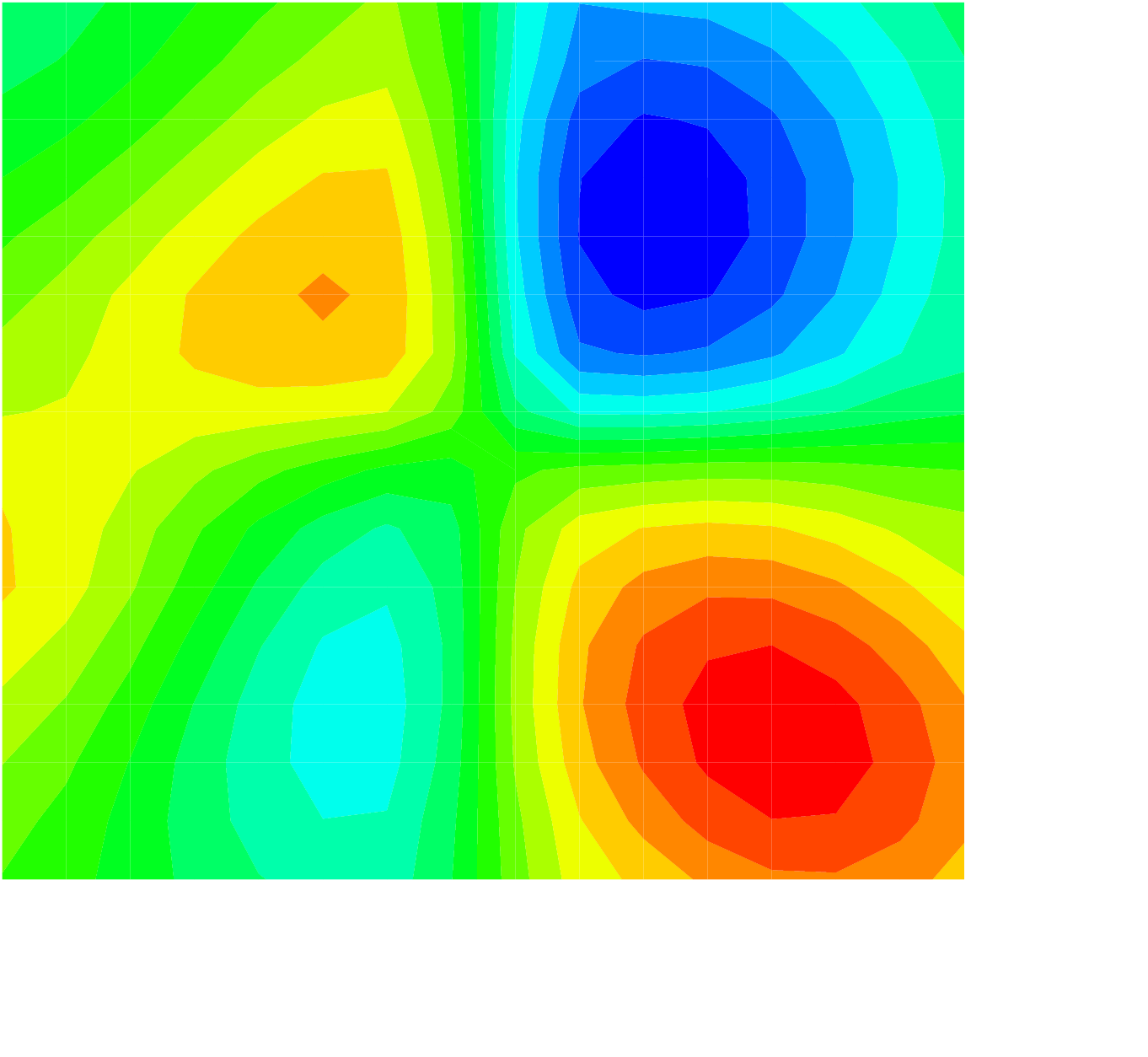}
\\

	\caption{First row: Reference log permeability filed. Second row: Accepted permeability fields in the global sampling method. Third row: Accepted permeability fields in MSM $2\times 2$.
	From left to right, log permeability fields at 20000, 40000 and 60000 iterations, respectively, from chain 2 in the first example.}
	\label{perm_1_2}
\end{figure}

We now compare the reference field with some of the simulated permeability fields from two selected chains. See Figures \ref{perm_1_1} and \ref{perm_1_2}. Other chains also show a similar behavior. Although both MCMC methods converged, we observe that at iteration 60000 the permeability field obtained from the multiscale sampling is closer to the reference permeability field than the field obtained without the multiscale sampling.  


\subsection{Example 2}
\begin{figure}[H]
	\centering
	\includegraphics[scale = 0.7]{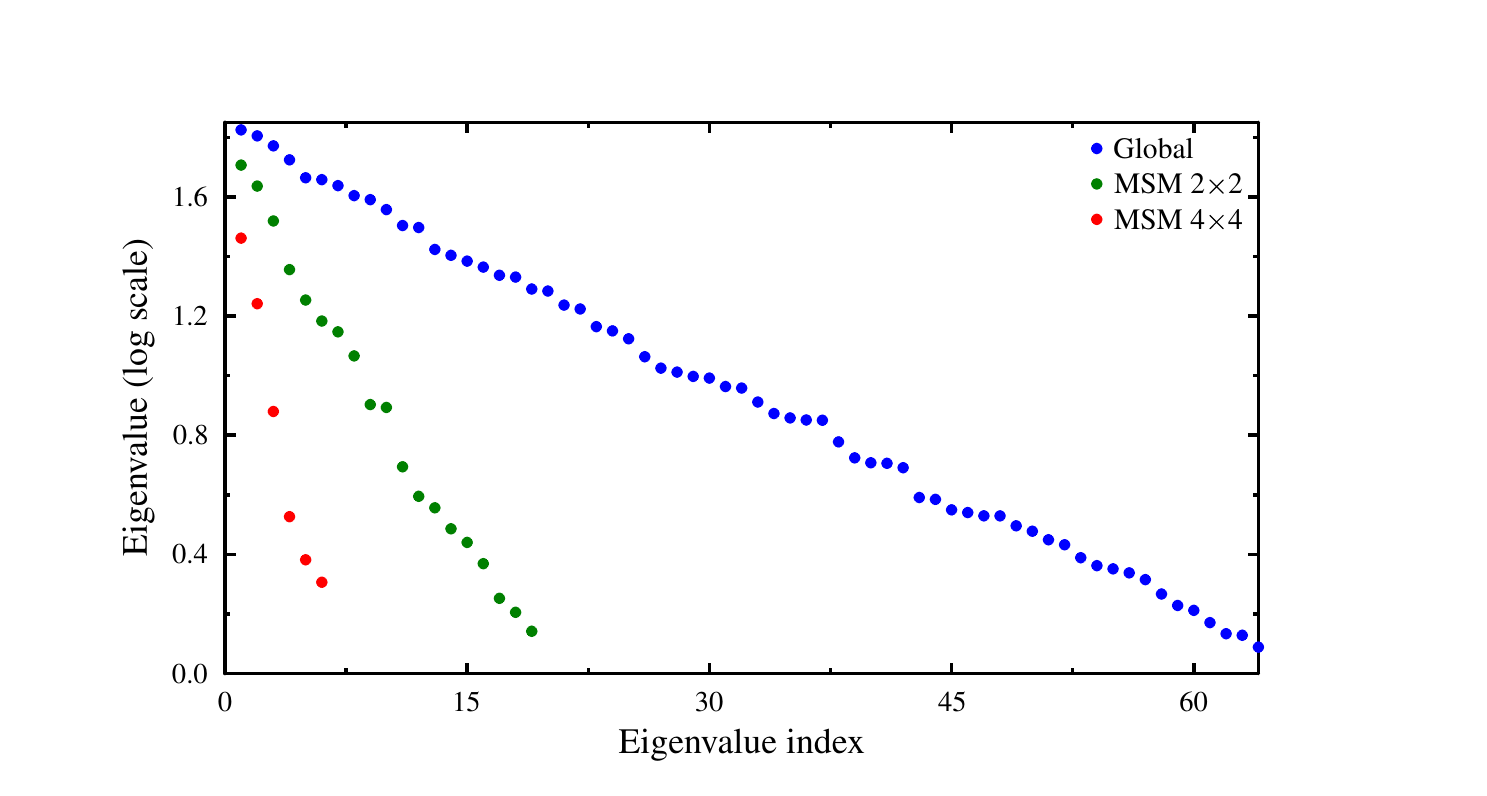}
	\caption{Decay of eigenvalues for the global and multiscale sampling in the second example.}
	\label{eigen_32x32_ch3}
\end{figure}

\begin{figure}[H]
	\centering
		\includegraphics[width=2.2in]{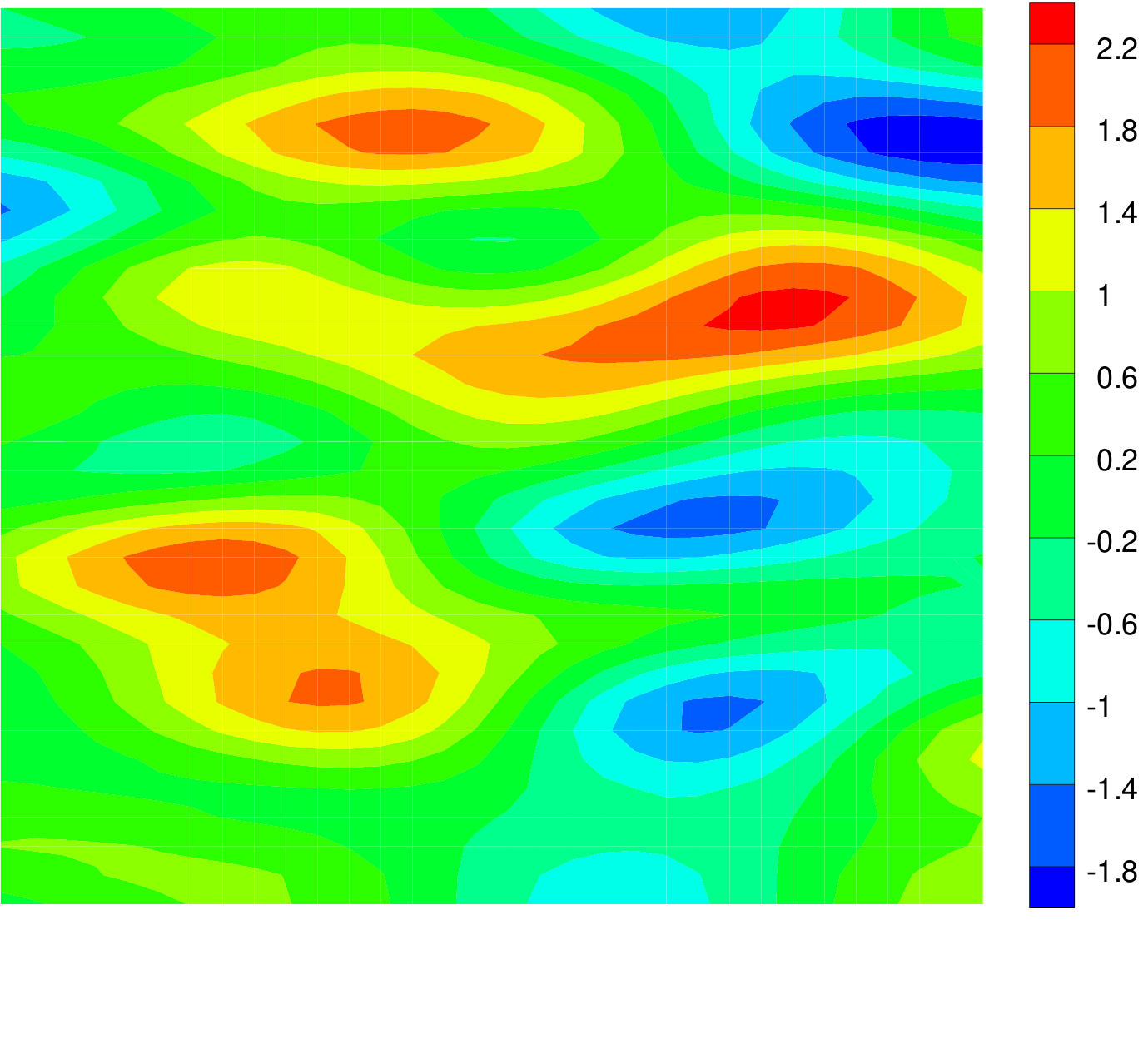}
	\hspace{5mm}			\includegraphics[width=2.2in]{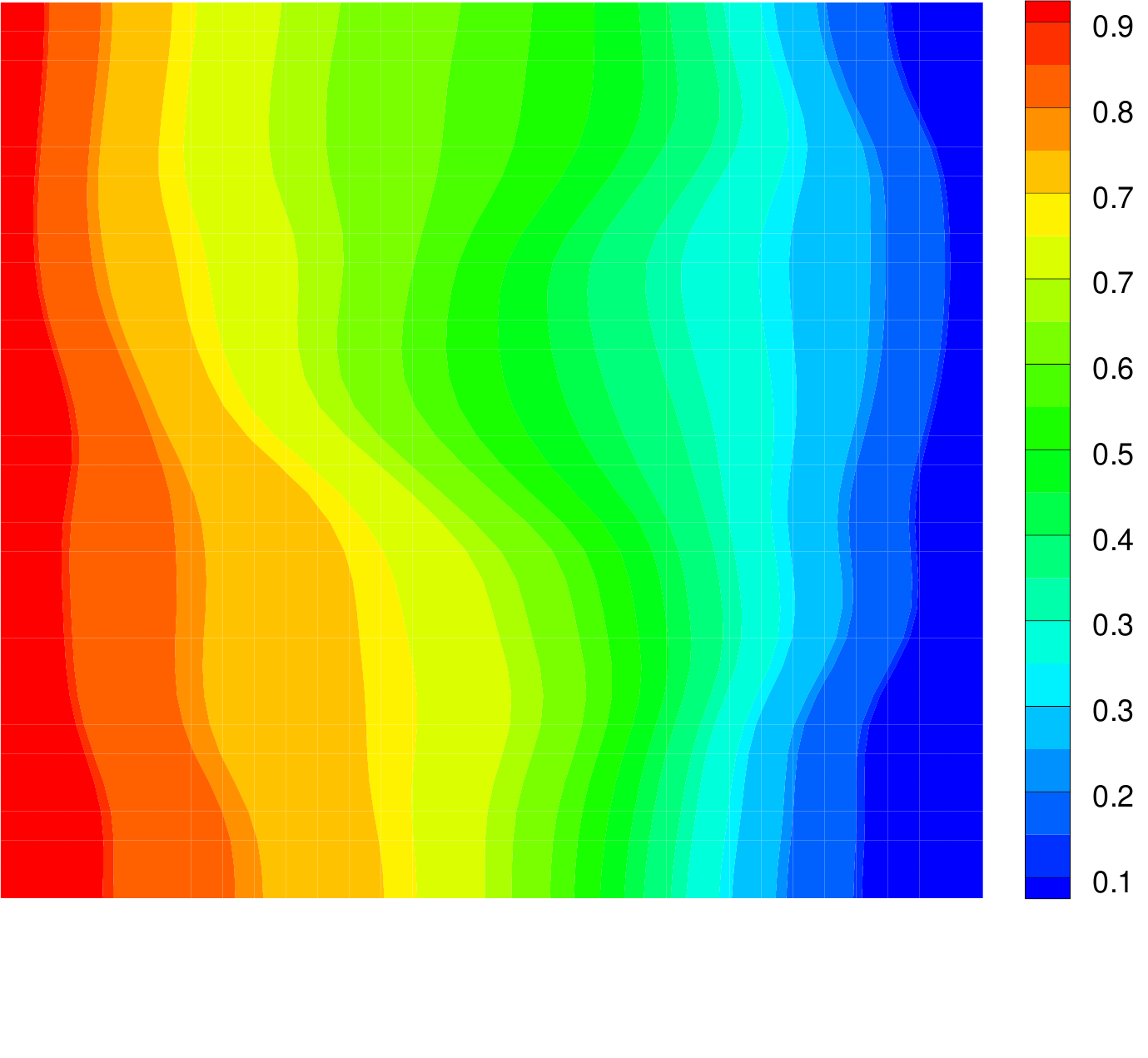}

	\caption{Reference log permeability field (left) and the corresponding reference pressure field (right) for the second example.}
	\label{ref_perm_2}
\end{figure}

In the second example we consider the case where the correlation lengths are not equal, i.e., $L_x= 0.2 $ and $L_y=0.06$ in Eq. \eqref{kle_3}. Figure  \ref{eigen_32x32_ch3} shows the decay of the eigenvalues (in log scale) for the global sampling, MSM $2\times 2$, and MSM $4\times 4$ methods.
In KLE, we consider the first $64$ eigenvalues, which preserve $97.87\%$ of the total energy, in the global sampling method. For MSM $2\times 2$ and MSM $4\times 4$, we take $16$ and $4$ eigenvalues, respectively. We generate a reference synthetic permeability field on a computational grid of size $32\times 32$ and then run the numerical simulator to generate the corresponding reference pressure field. See Figure  \ref{ref_perm_2} for the reference fields.	
We use the same coarse mesh of size $8\times 8$ as in the first example. We also use the same local blocking number $N_{lb}=1$. We set the tuning parameter $\beta= 0.75$ in Eq. \eqref{RW_sampler}. 

\begin{figure}[H]
	\centering
	\includegraphics[scale = 0.55]{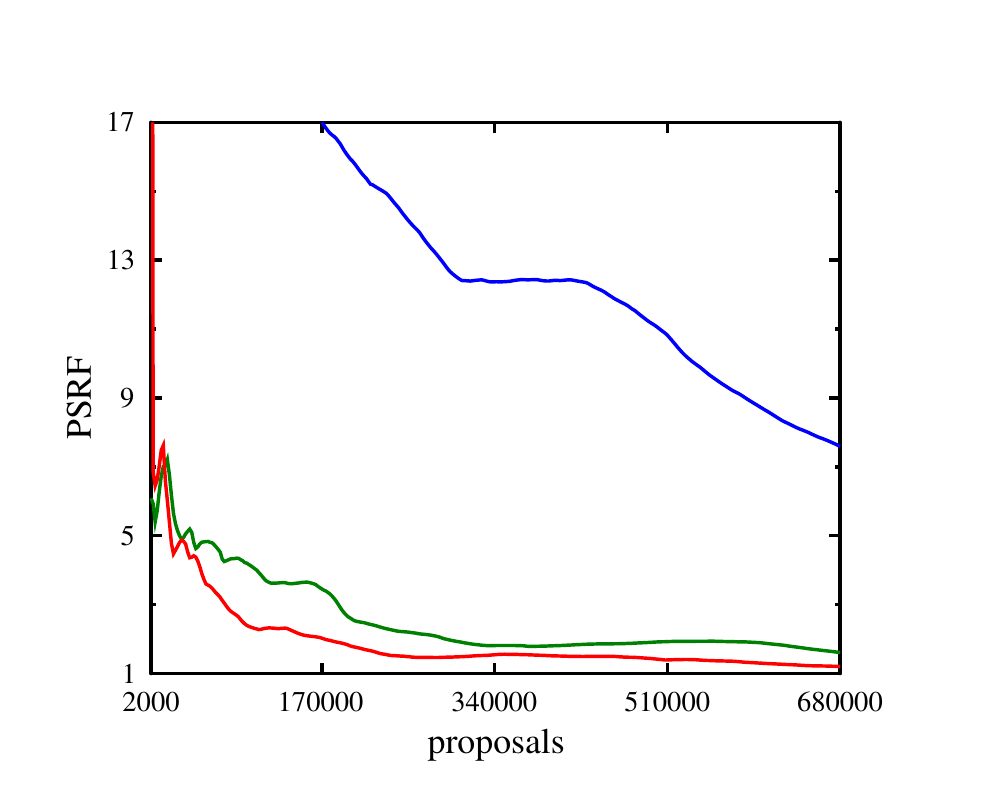}
	\includegraphics[scale= 0.55]{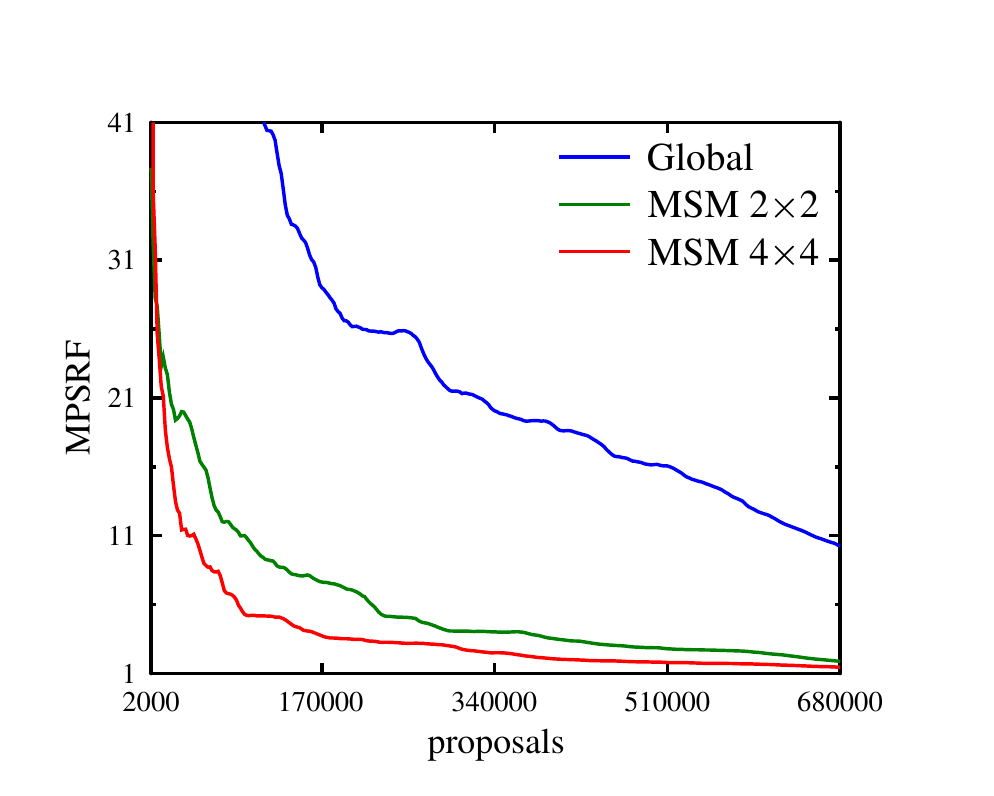}
	\caption{The maximum of PSRFs and MPSRF for the MCMC method with and without multiscale sampling for the second example.}
	\label{MPSRF_32x32_ch3}
\end{figure}

\begin{figure}[H]
	\centering
	\includegraphics[scale = 0.9]{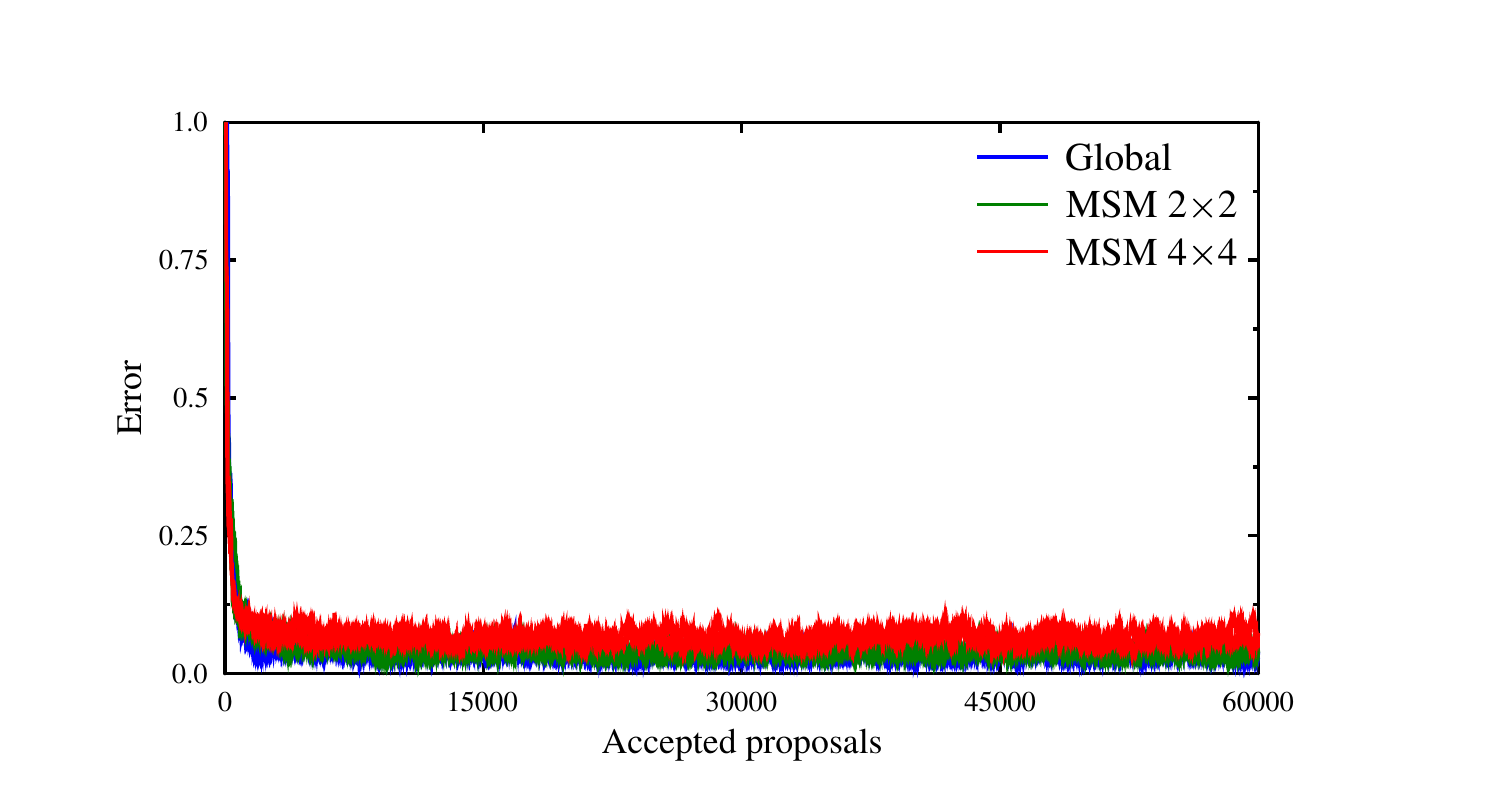}
	\caption{Error curves of the preconditioned MCMC with and without multiscale sampling for the second example.}
	\label{error_32x32_ch3}
\end{figure}

\begin{figure}[H]
	\centering
	\includegraphics[width=1.5in]{figures/ex_2_perms/ref_32x32_crop.pdf}\\
	
	\includegraphics[width=1.5in]{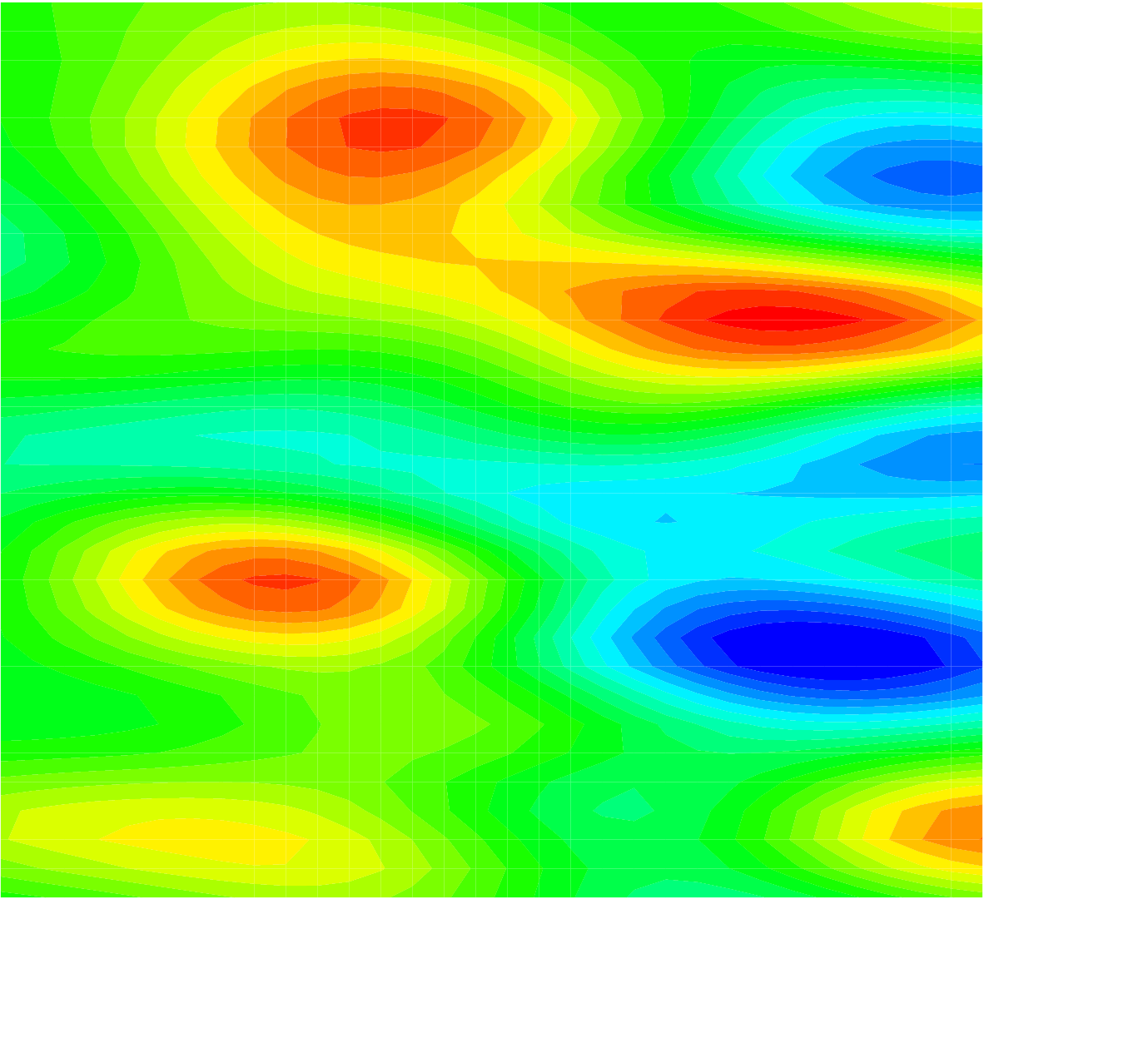}				\includegraphics[width=1.5in]{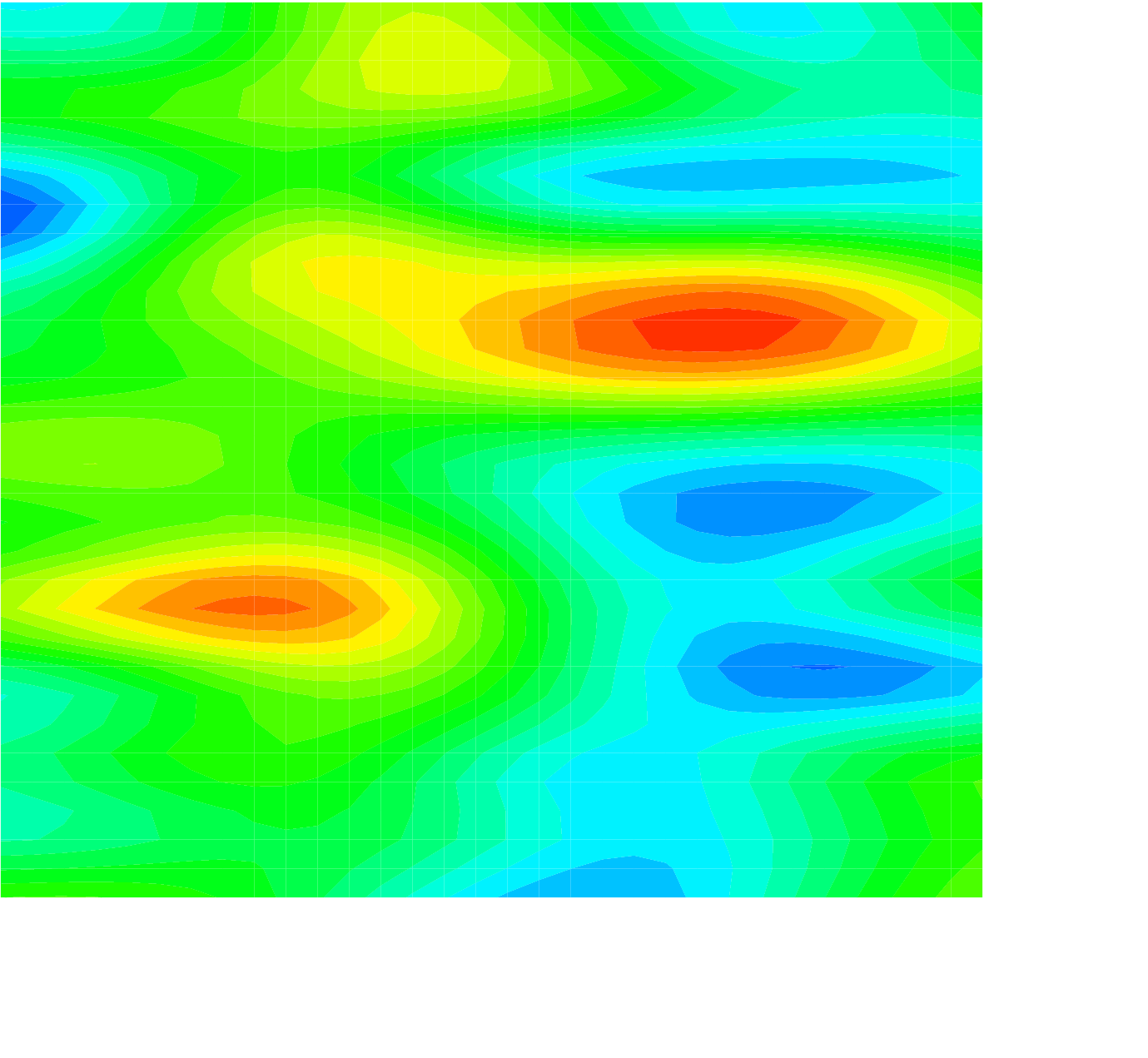}						\includegraphics[width=1.5in]{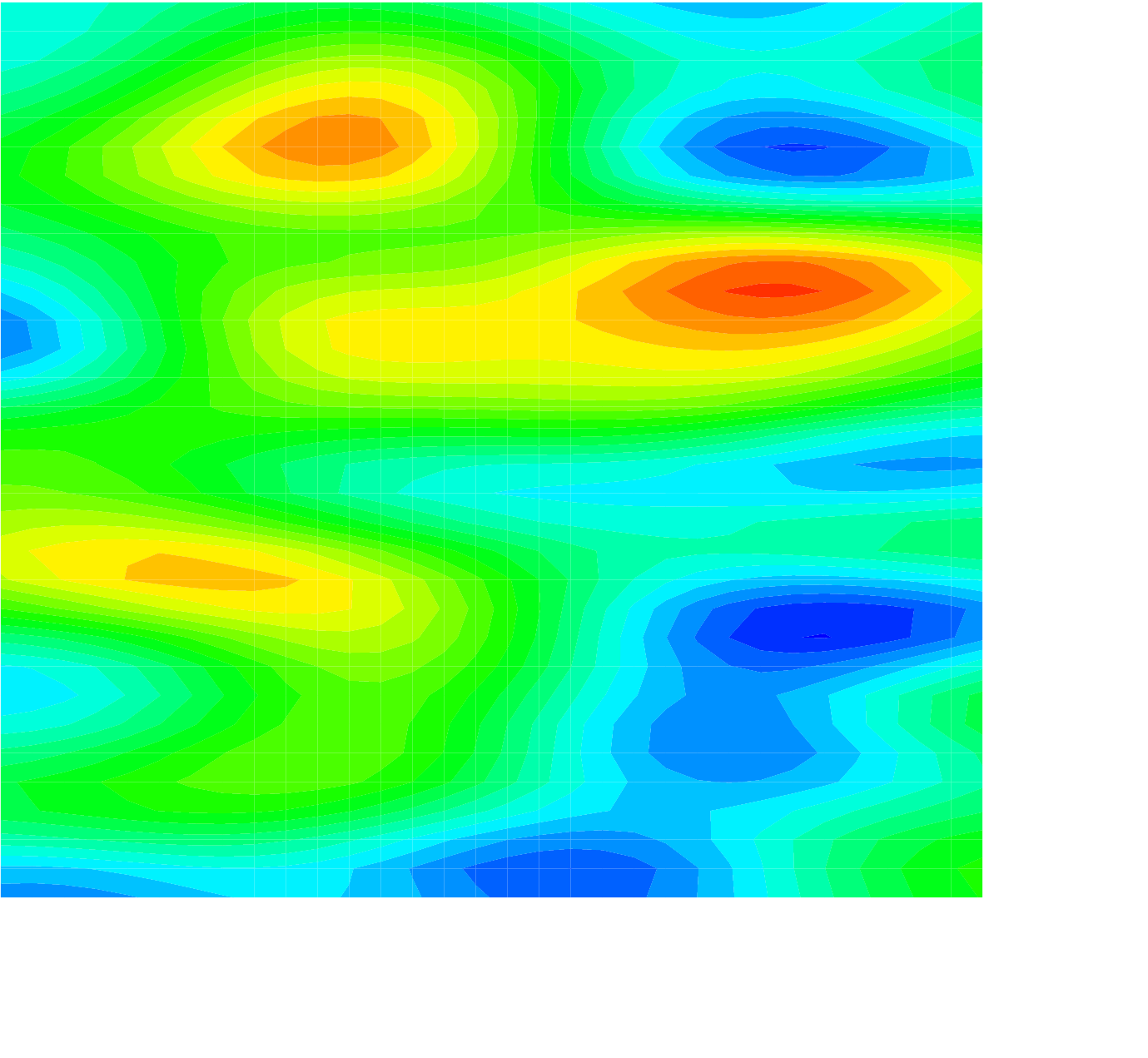}\\								
	\includegraphics[width=1.5in]{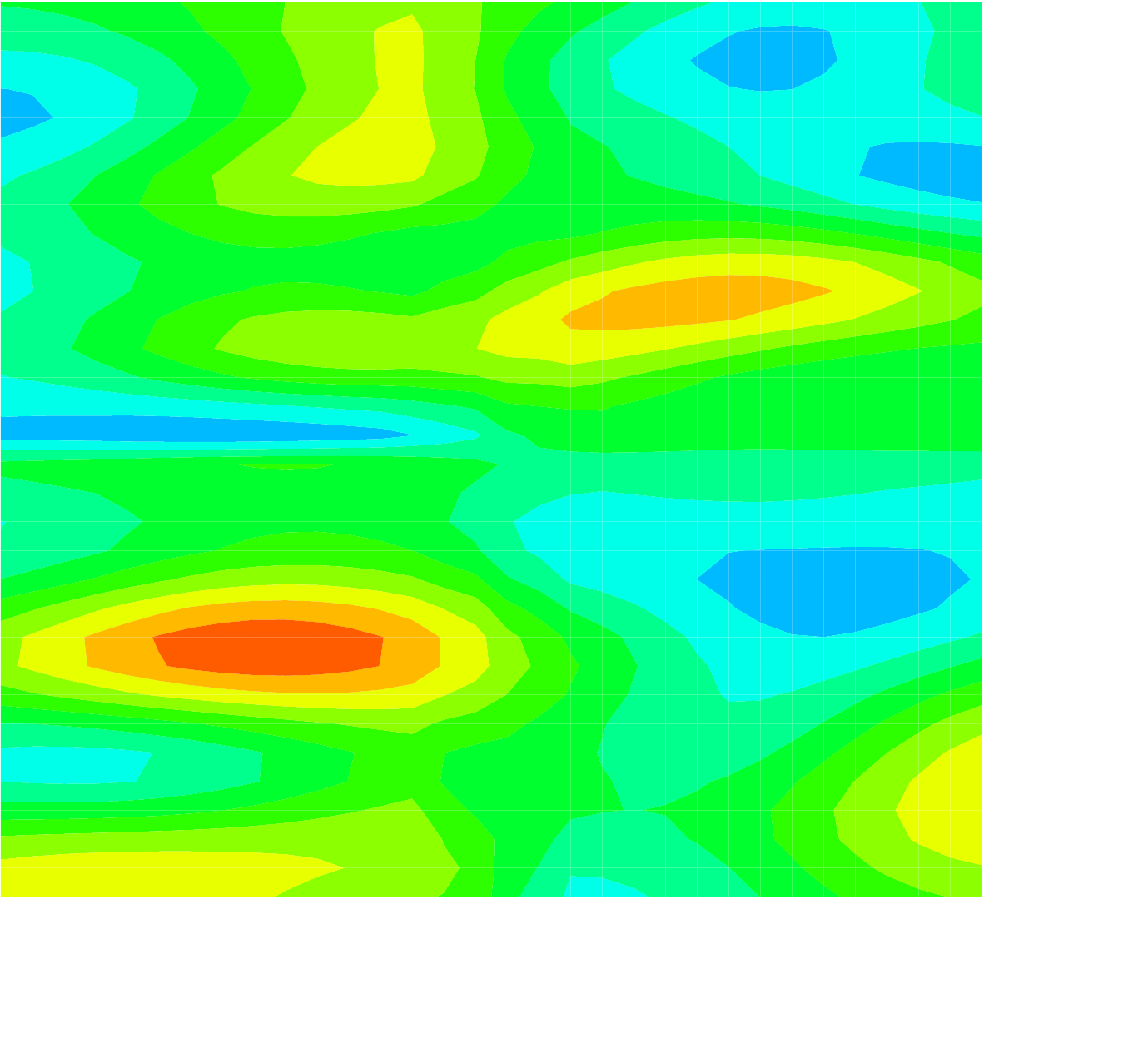}								\includegraphics[width=1.5in]{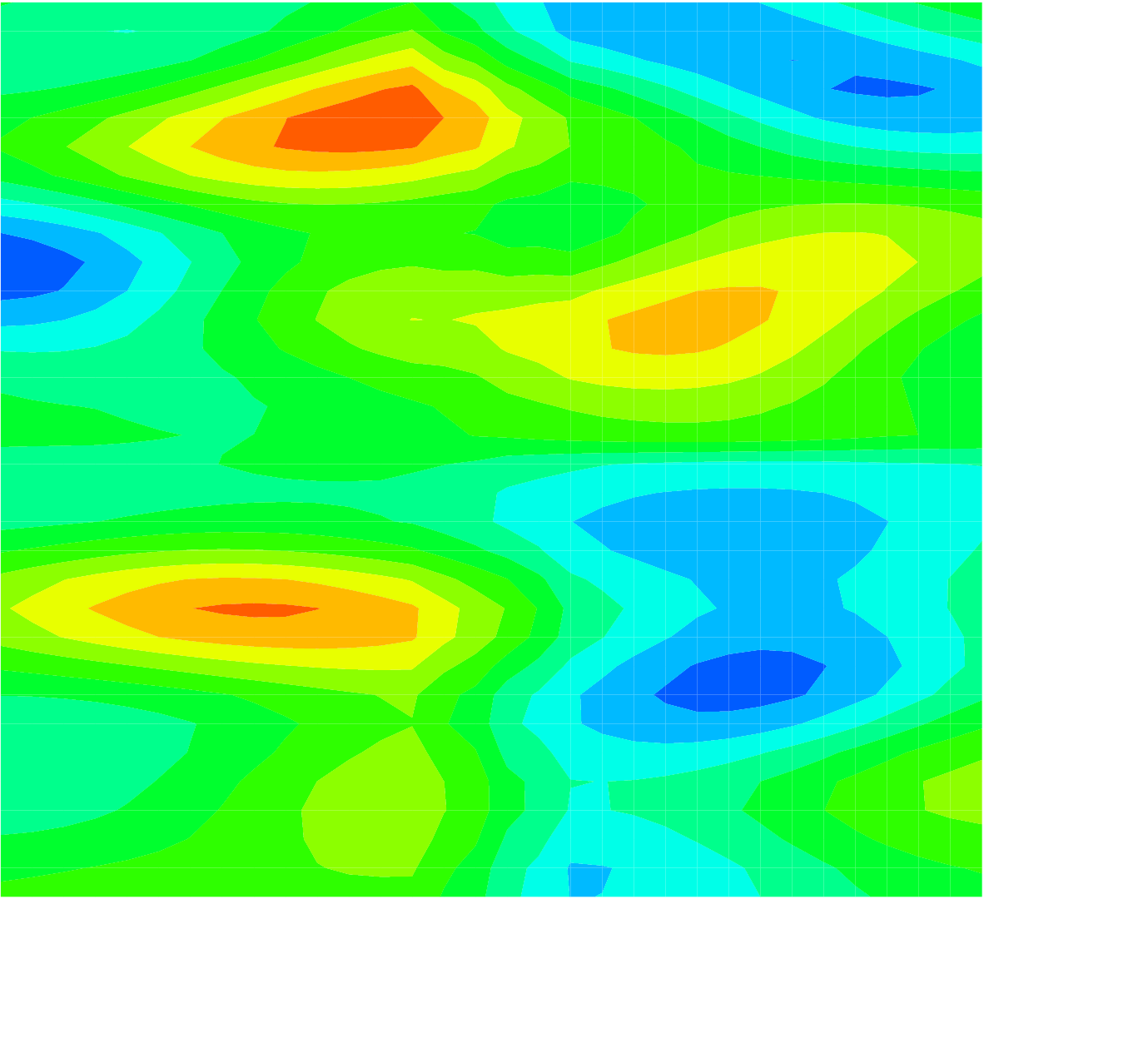}							\includegraphics[width=1.5in]{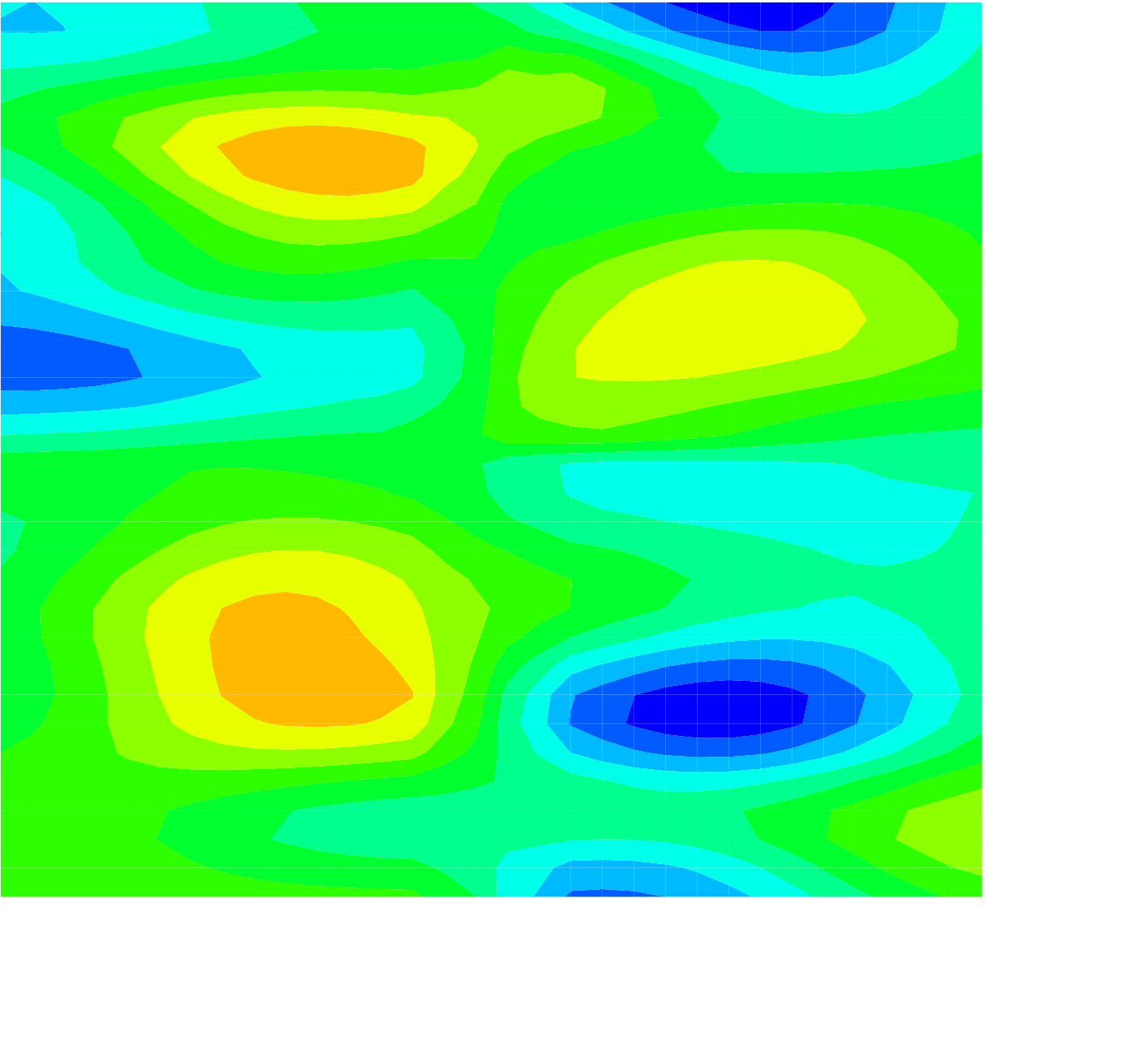}\\
		
	\includegraphics[width=1.5in]{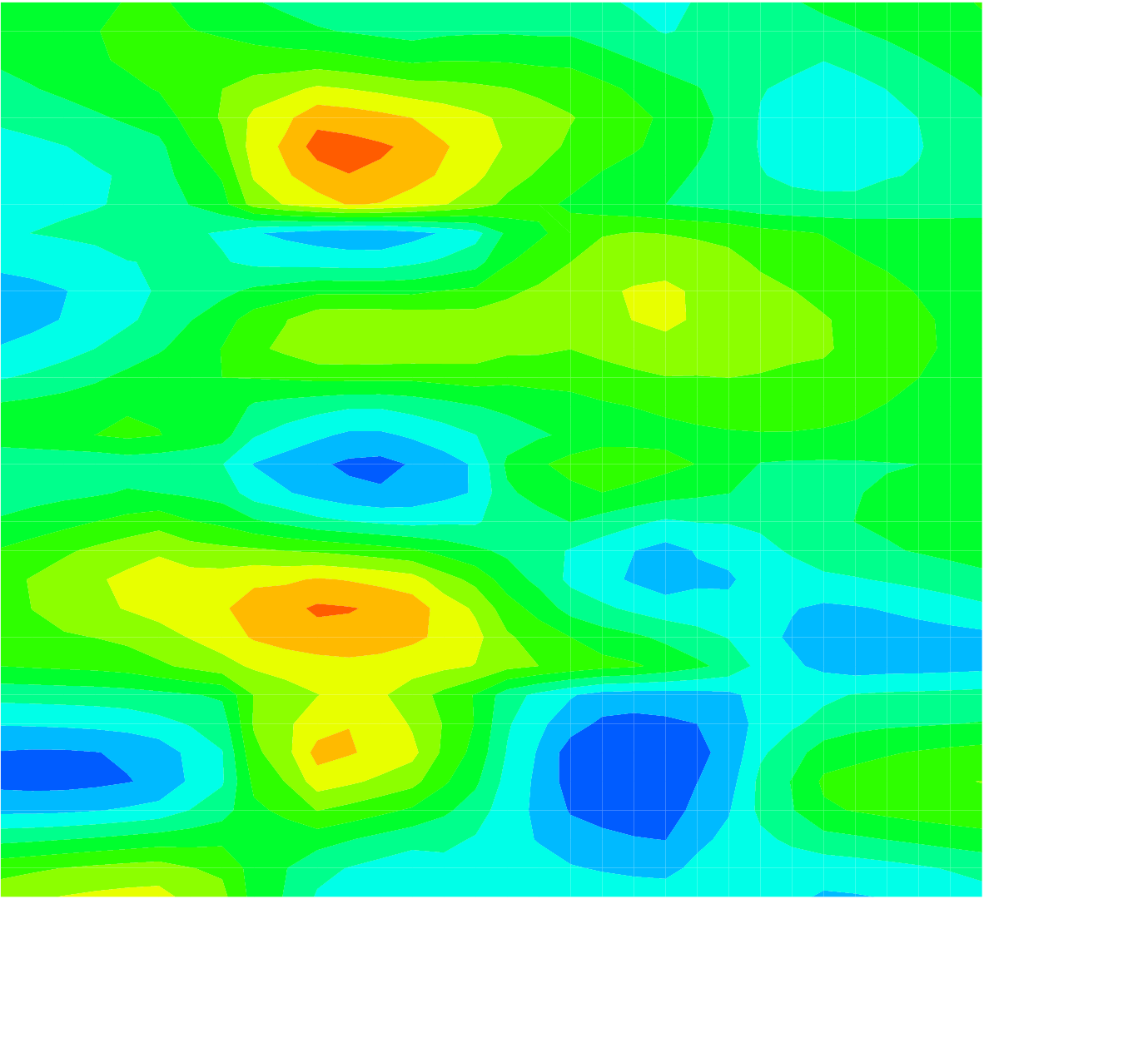}
	\includegraphics[width=1.5in]{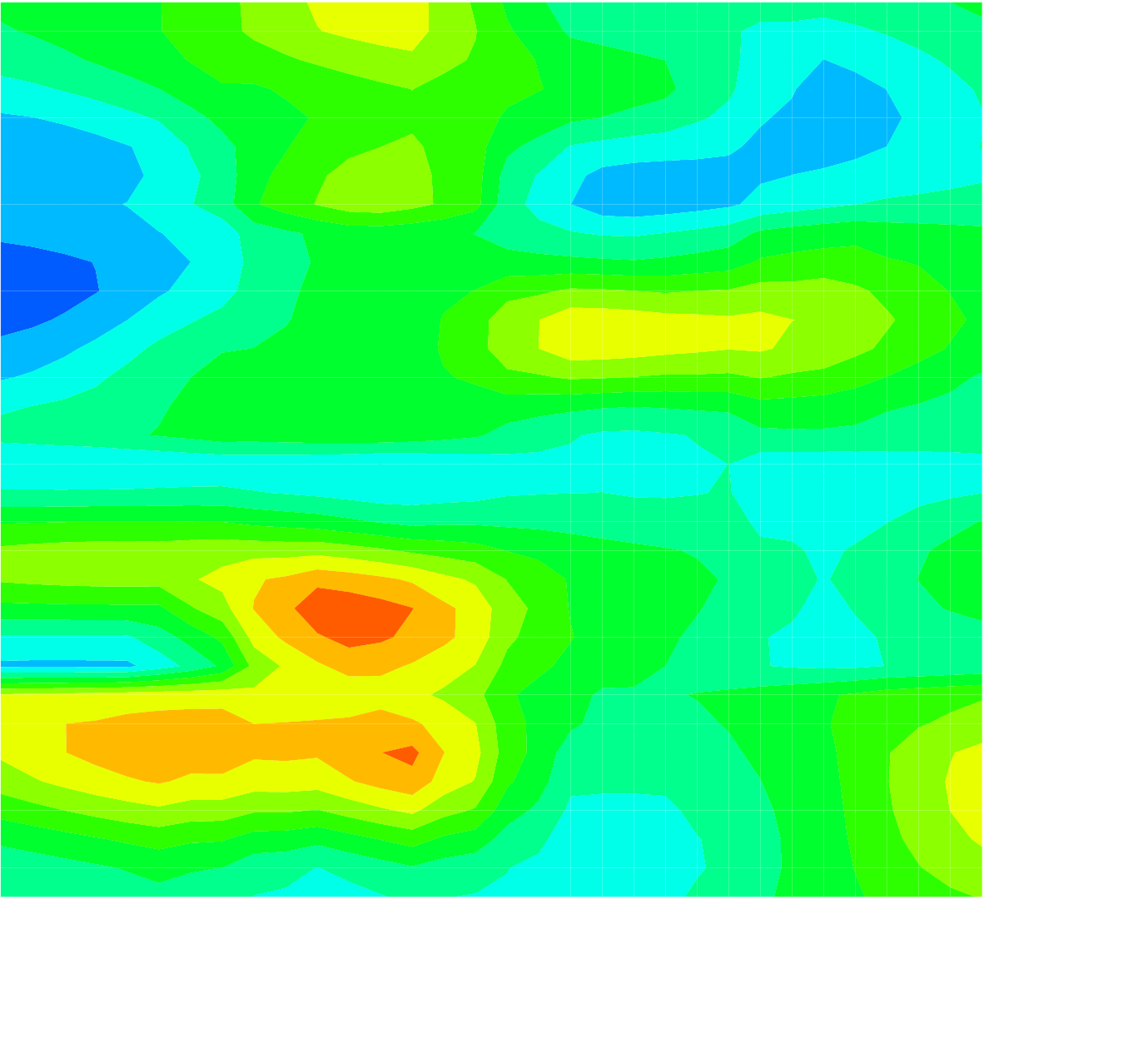}
	\includegraphics[width=1.5in]{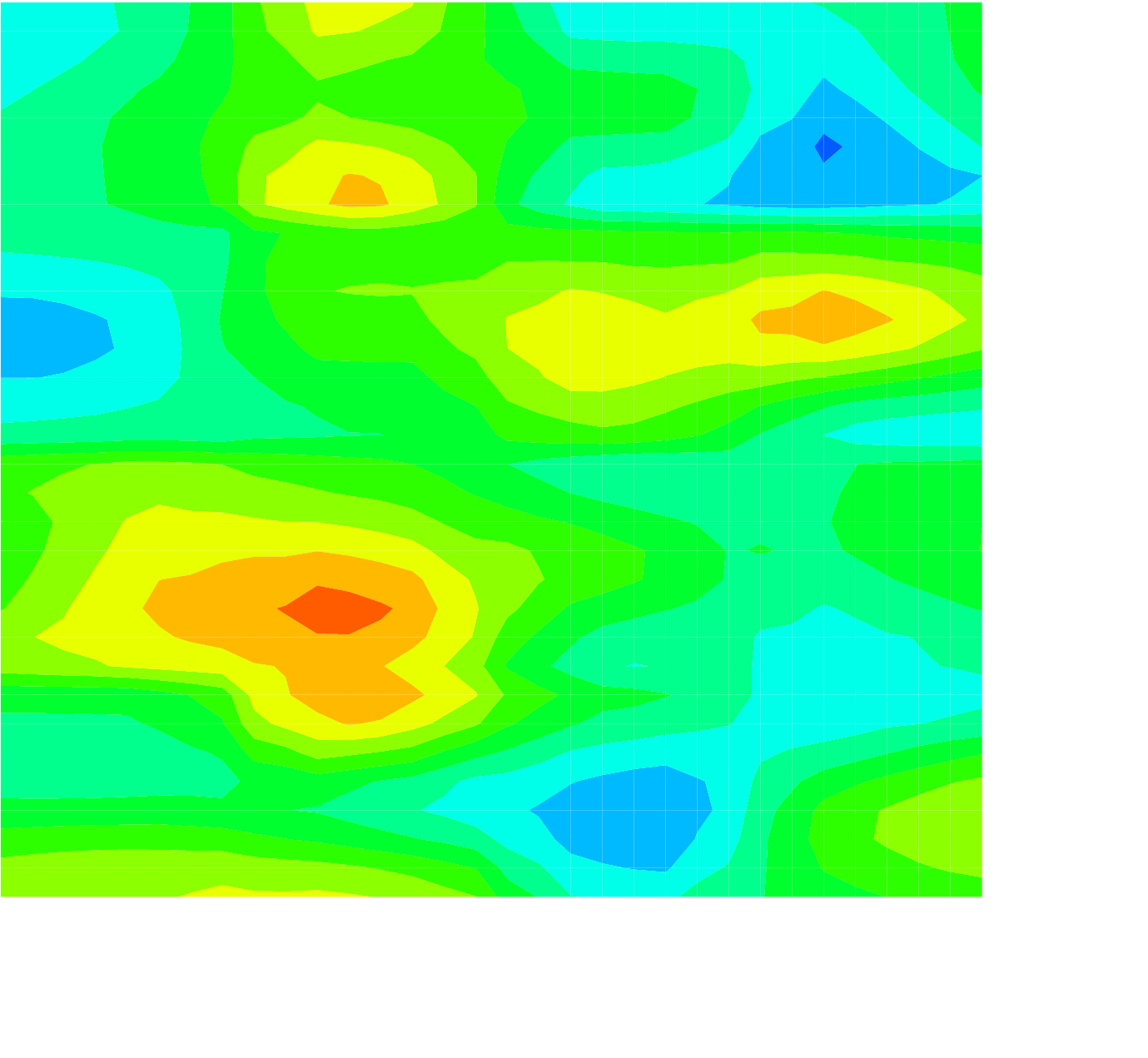}
				
	\caption{First row: Reference log permeability filed. Second row: Accepted permeability fields in the global sampling method. Third row: Accepted permeability fields in MSM $2\times 2$. Fourth row: Accepted permeability fields in MSM $4\times 4$.
	From left to right, log permeability fields at 20000, 50000 and 100000 iterations, respectively, from chain 1 in the second example.}
	\label{perm_2_chain1}
\end{figure}

\begin{figure}
	\centering
	
\includegraphics[width=1.5in]{figures/ex_2_perms/ref_32x32_crop.pdf}\\

\includegraphics[width=1.5in]{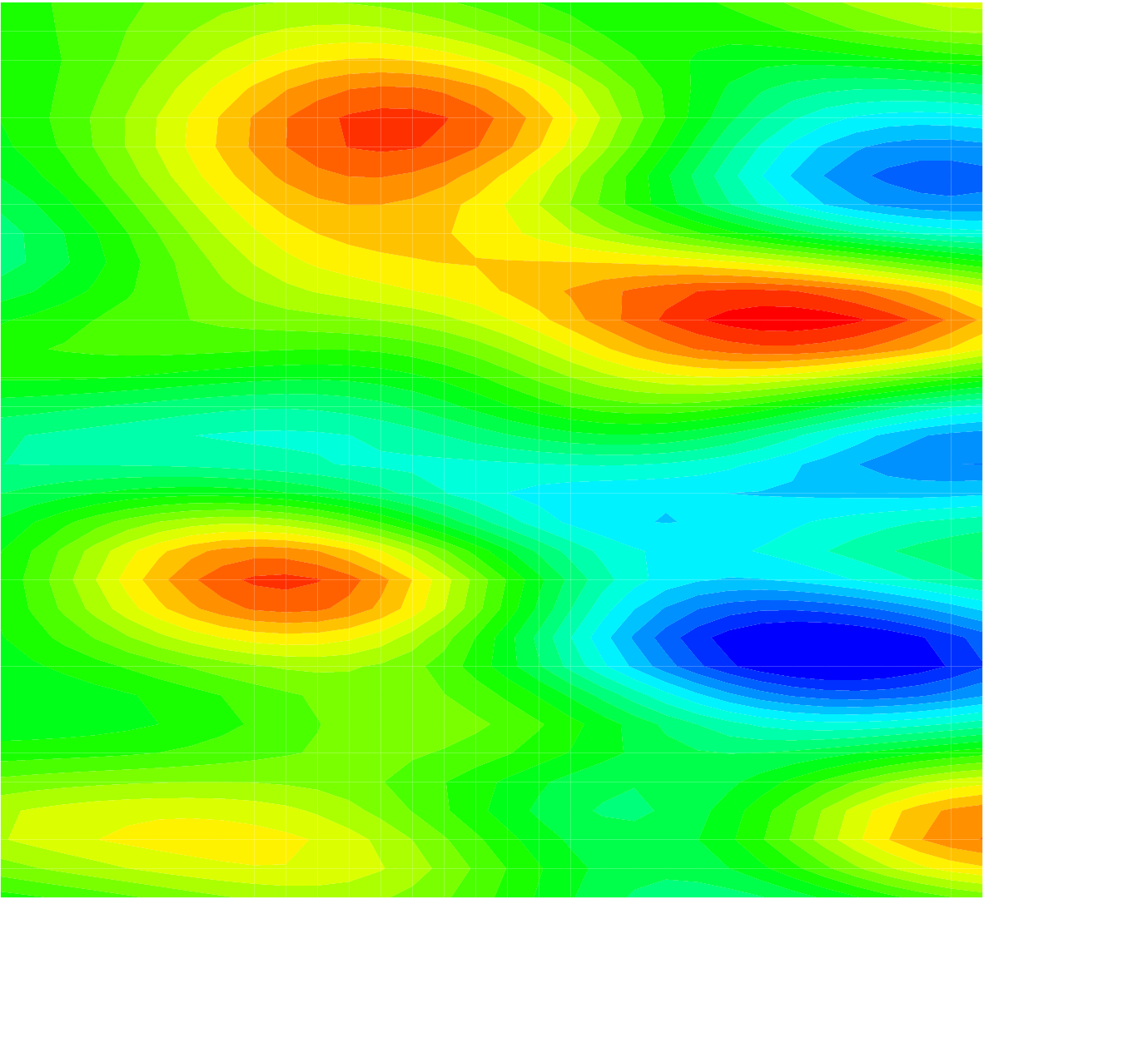}
\includegraphics[width=1.5in]{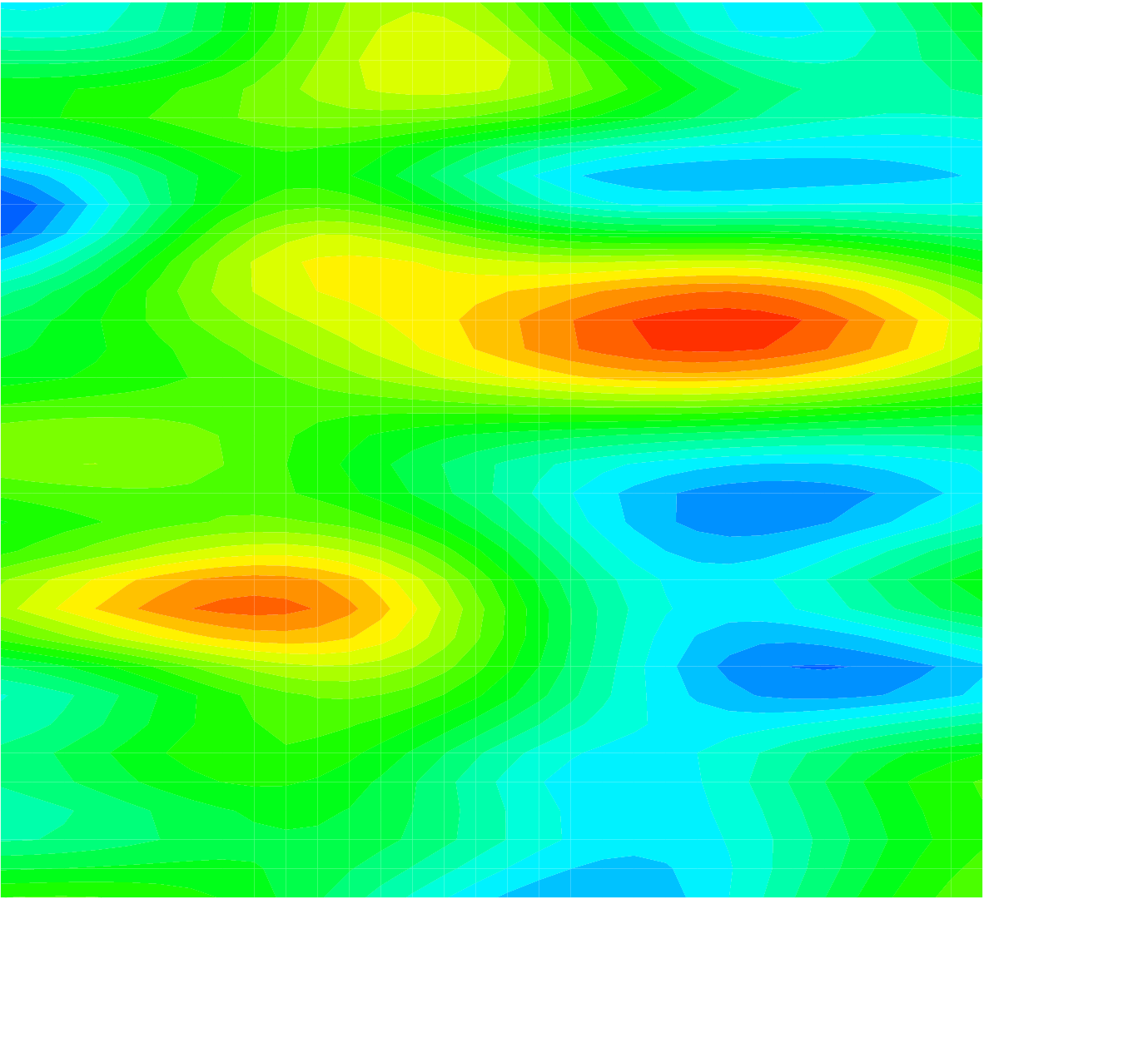}
\includegraphics[width=1.5in]{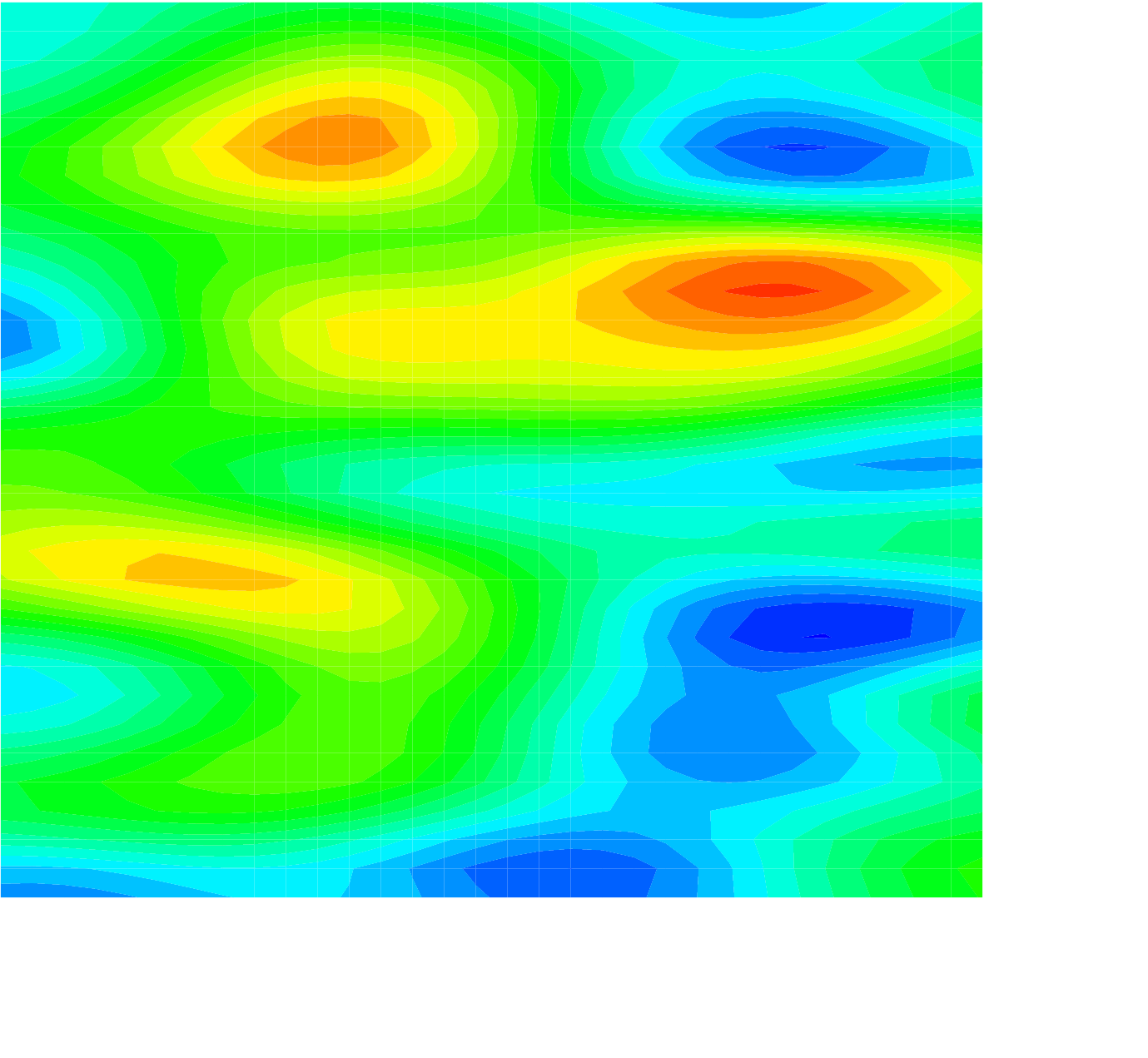}\\

\includegraphics[width=1.5in]{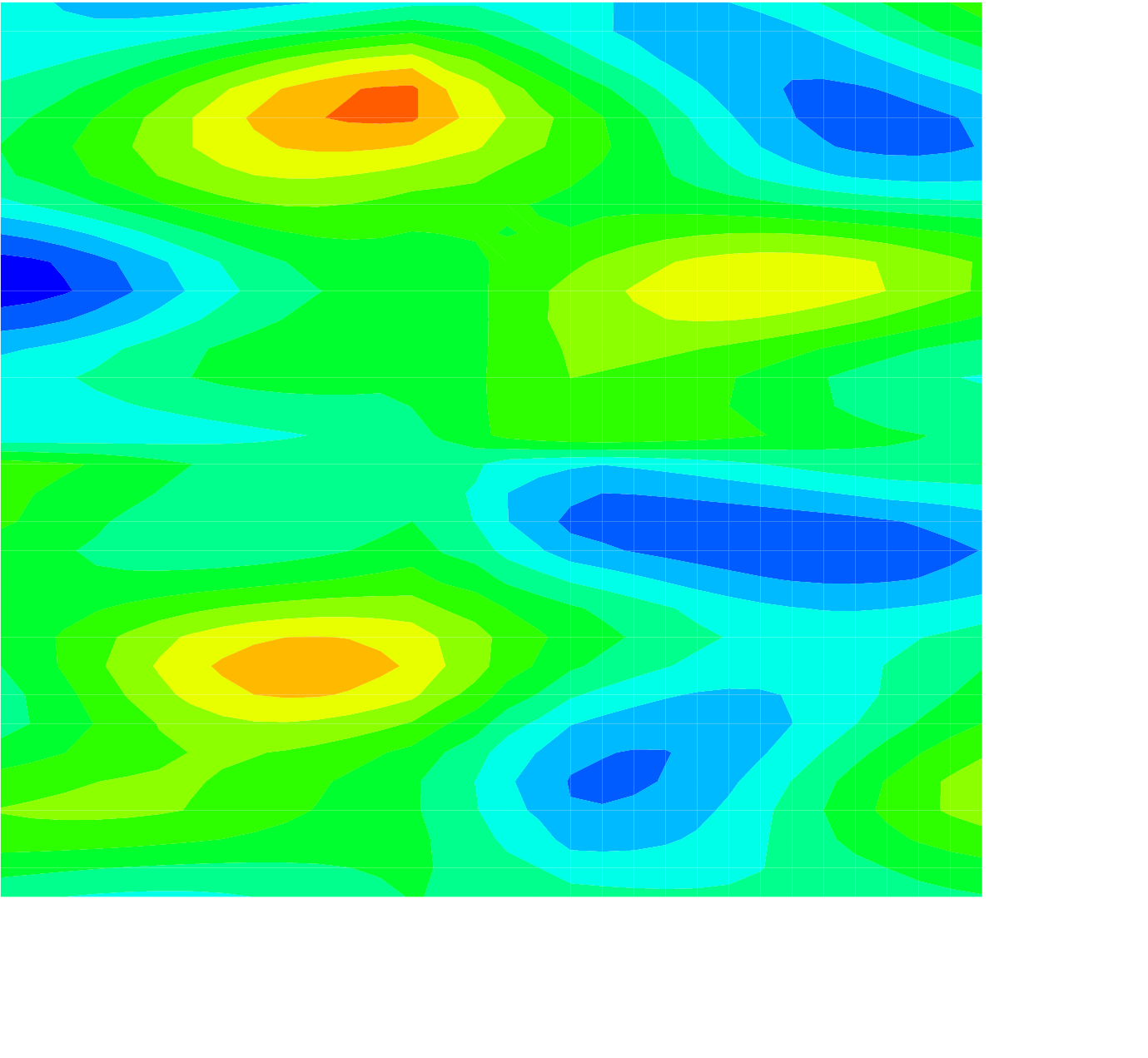}
\includegraphics[width=1.5in]{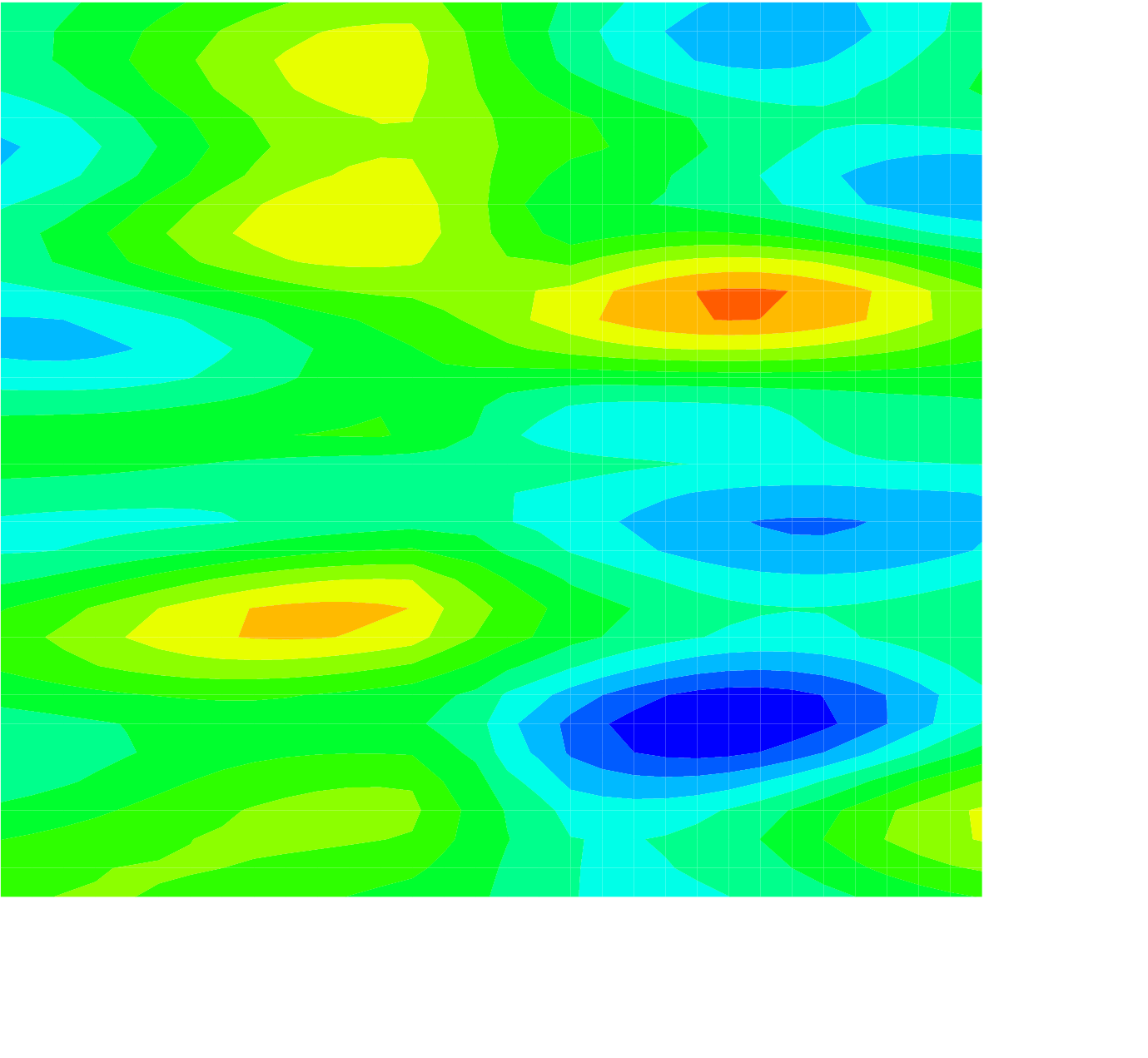}
\includegraphics[width=1.5in]{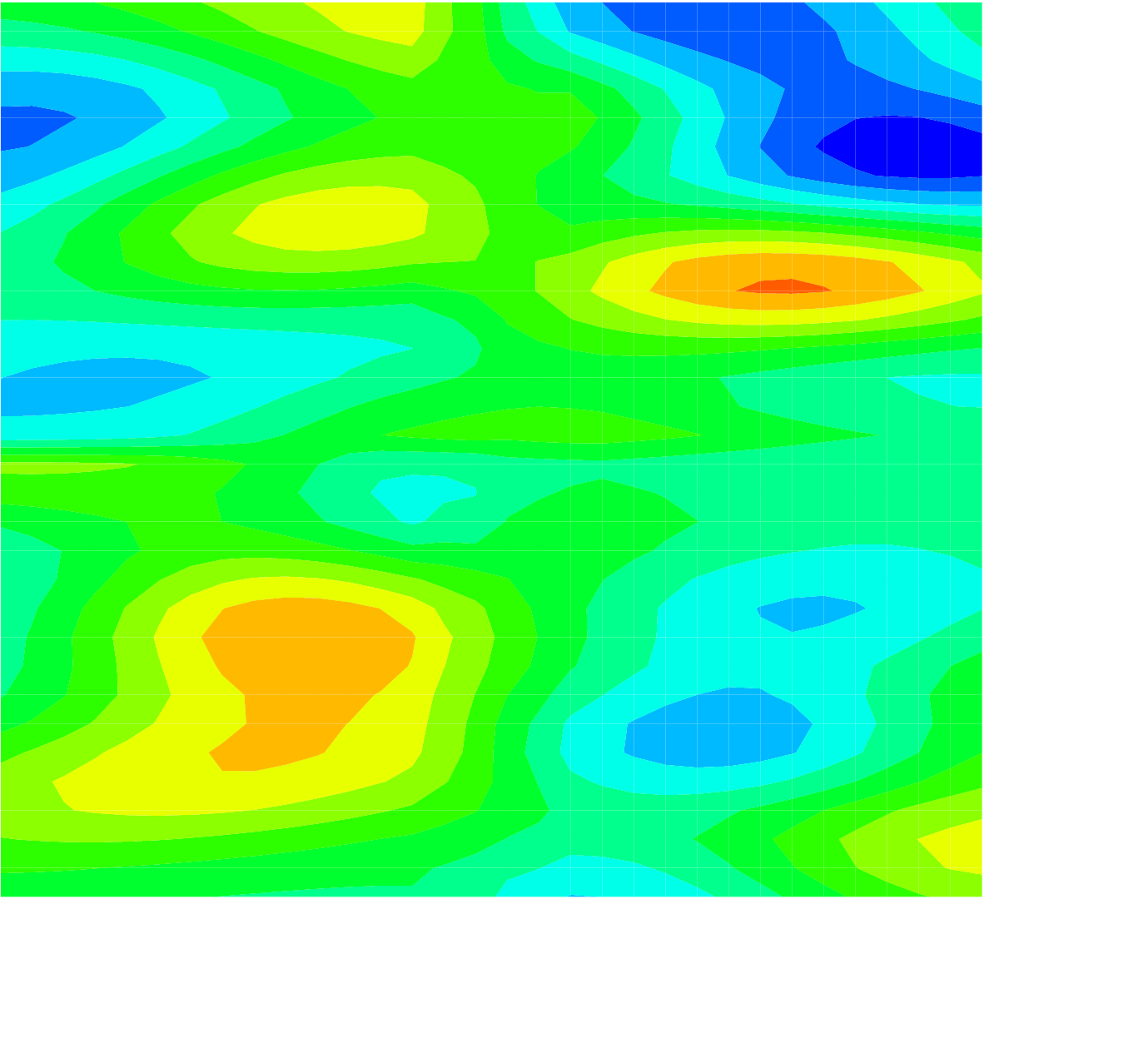}\\

\includegraphics[width=1.5in]{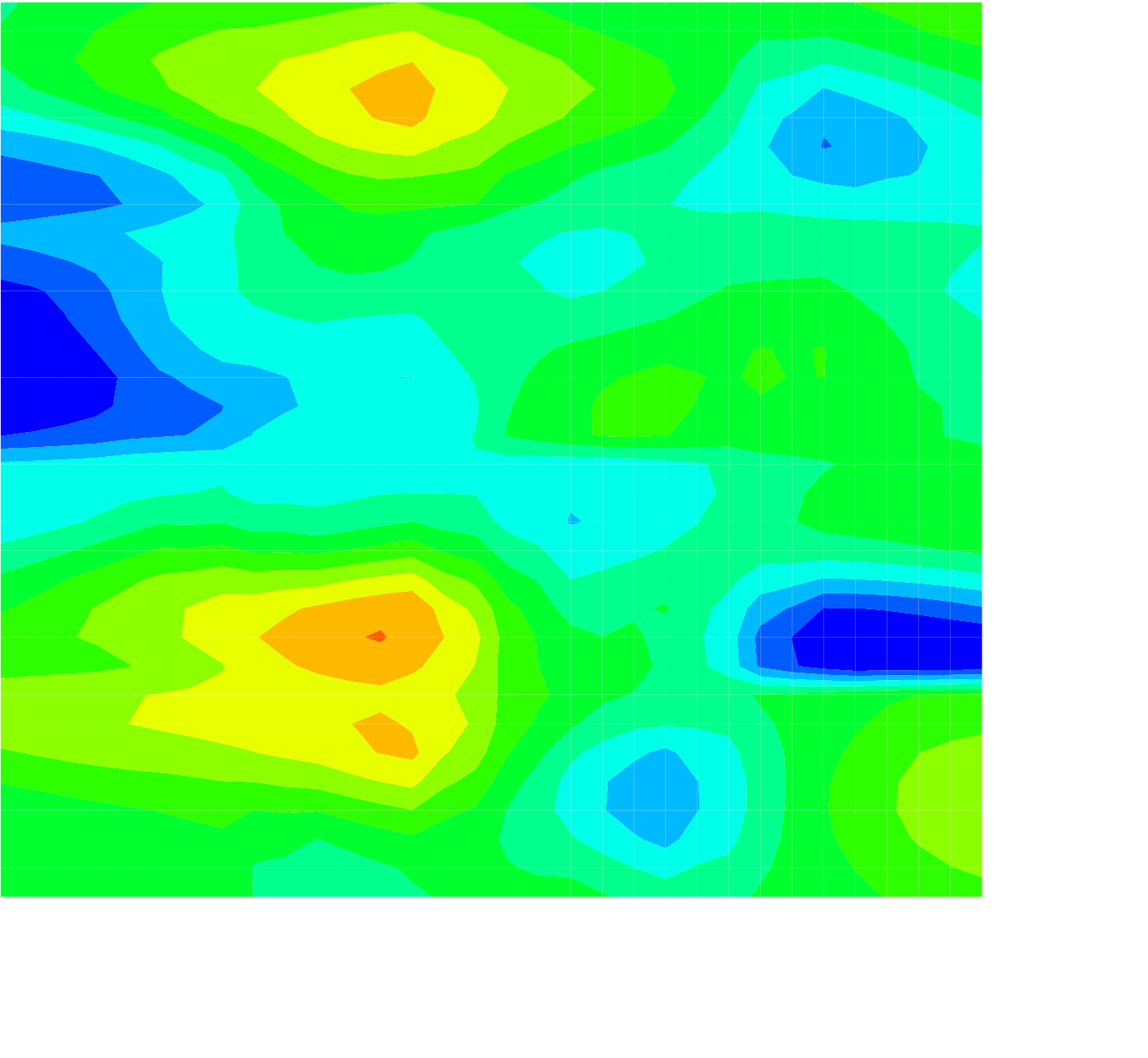}
\includegraphics[width=1.5in]{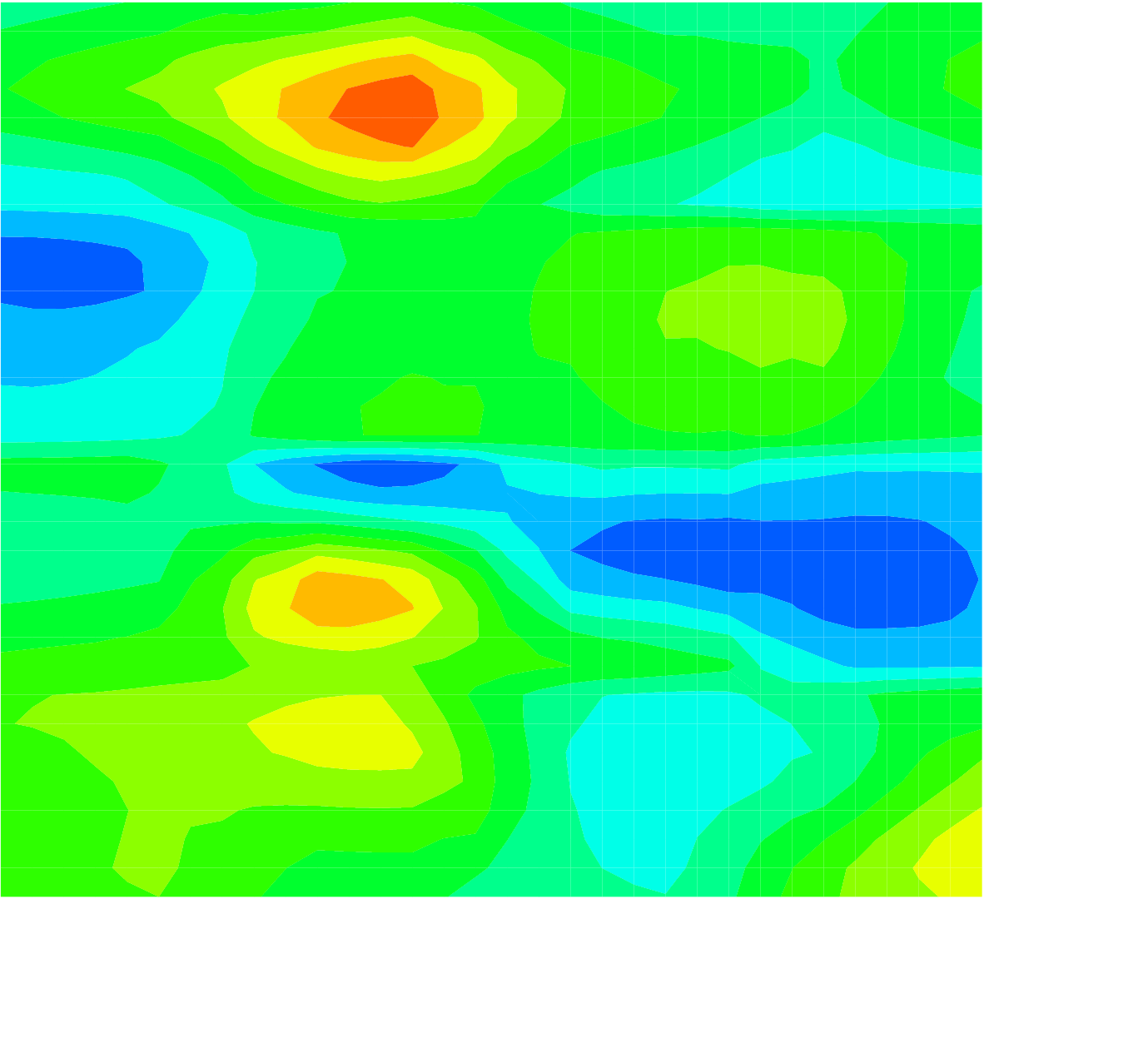}
\includegraphics[width=1.5in]{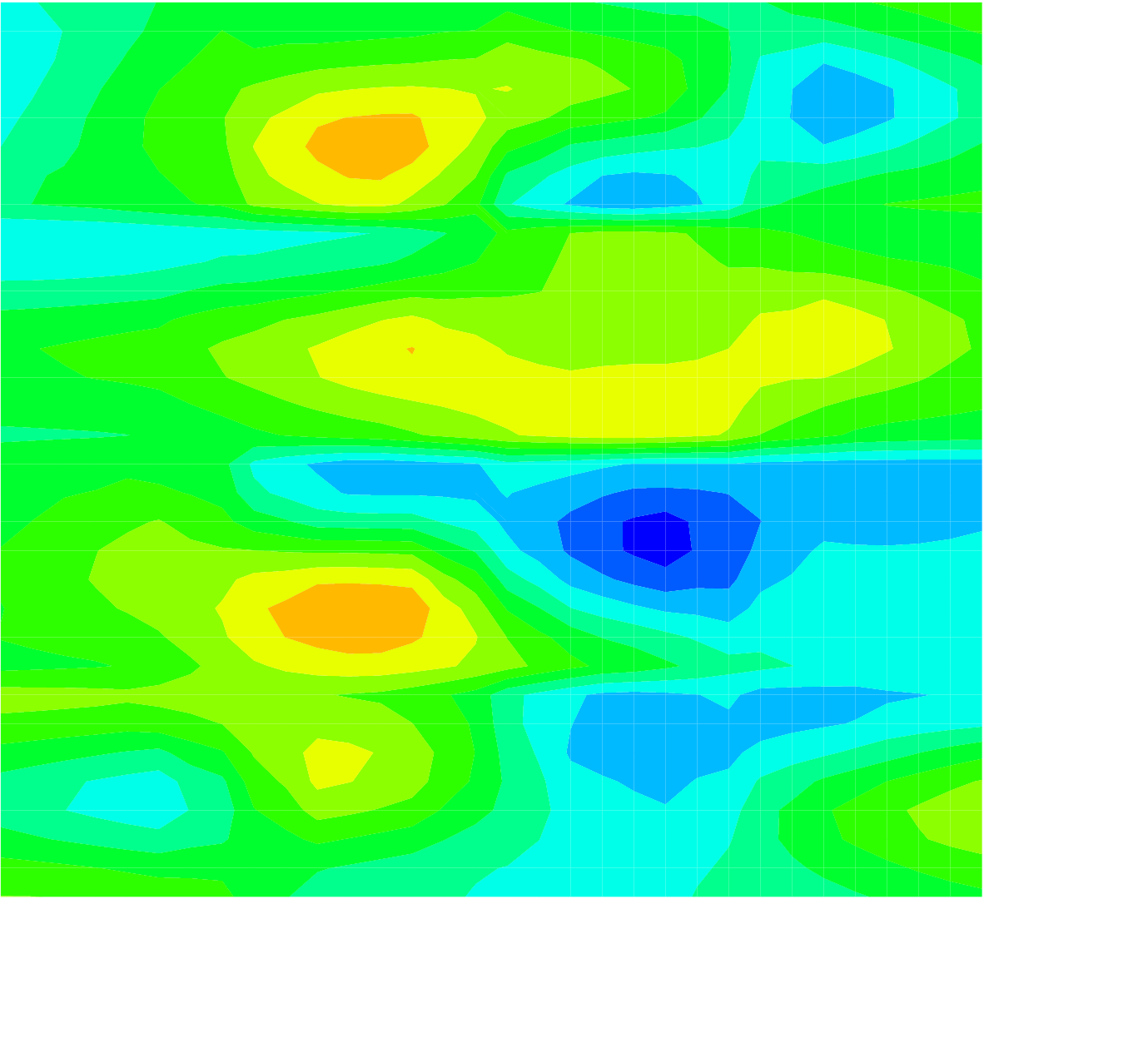}

	\caption{First row: Reference log permeability filed. Second row: Accepted permeability fields in the global sampling method. Third row: Accepted permeability fields in MSM $2\times 2$. Fourth row: Accepted permeability fields in MSM $4\times 4$.
	From left to right, log permeability fields at 20000, 50000 and 100000 iterations, respectively, from chain 2 in the second example.}
	\label{perm_2_chain2}
\end{figure}

Let us consider the convergence analysis of these methods. We take $170000$ proposals from each chain to compute the MPSRF and the maximum of PSRFs.  Figure  \ref{MPSRF_32x32_ch3} shows the maximum of PSRFs and MPSRF curves. At the tails of the maximum of PSRFs and MPSRF curves, we have the values $1.2$  and $1.4$, respectively, for MSM $4\times4$. These values are slightly higher for MSM $2\times2$. Thus, we can conclude that MSM $4\times 4$ converges to the stationary distribution faster than MSM $2\times 2$. On the other hand, the PSRF and MPSRF curves in the global sampling method do not show any sign of converging at the same number of iterations. Moreover, the acceptance rates are better for the multiscale sampling methods. See Table \ref{Lx_Ly_2}. The error curves for this study are shown in Figure  \ref{error_32x32_ch3}. They are very comparable. Figures  \ref{perm_2_chain1} and \ref{perm_2_chain2} present simulated permeability fields from two selected chains. From these figures, we observe that both MSM $2\times 2$ and MSM $4\times 4$ recover the permeability fields better than the global sampling method. 
\begin{table}[H] 
	\caption{A comparison of acceptance rates for the MCMC with and without MSM in the second example.}
	\center
	\begin{tabular}{|cccc|}
		\hline
		&  \quad  MCMC global & \quad MCMC with MSM $2\times 2$ & \quad MCMC with MSM $4\times 4$  \\
		\hline
		$\sigma_F^2$ &   $10^{-3}$  & $10^{-3}$  & $10^{-3}$   \\               
		$\sigma_C^2$   & $5\times10^{-3}$ & $5\times10^{-3}$&  $5\times10^{-3}$ \\  
		acc. rate & $50\%$   & $54\%$& $55\%$  \\
		\hline
	\end{tabular}
	\label{Lx_Ly_2}     
\end{table}


\subsection{Example 3}
In this example, we test the proposed method on a large grid size of $64\times 64$ with the correlation lengths, $L_x= L_y = 0.1 $ in Eq. \eqref{kle_3}. Figure  \ref{eigen_16x16_ch3} shows the decay of the eigenvalues for the methods, the global sampling, MSM $2\times 2$, and MSM $4\times 4$. We consider $64$ eigenvalues, which preserve $95.8\%$ of the total energy in KLE, in the global sampling. The numbers of eigenvalues for MSM $2\times 2$ and MSM $4\times 4$ are 16 and 4, respectively. The synthetic reference permeability field is generated on a computational mesh of size $64\times 64$. Then, the numerical simulator is used to generate the corresponding reference pressure field. Figure  \ref{ref_perm_3_1} shows the reference permeability field and the corresponding pressure distribution on the grid. We use a coarse mesh of size $16\times 16$ in the filtering step in the preconditioned MCMC. Let $N_{lb}=2$. We set $\beta = 0.2$ in Eq. \eqref{RW_sampler}. 
\begin{figure}[H]
	\centering
	\includegraphics[scale = 0.7]{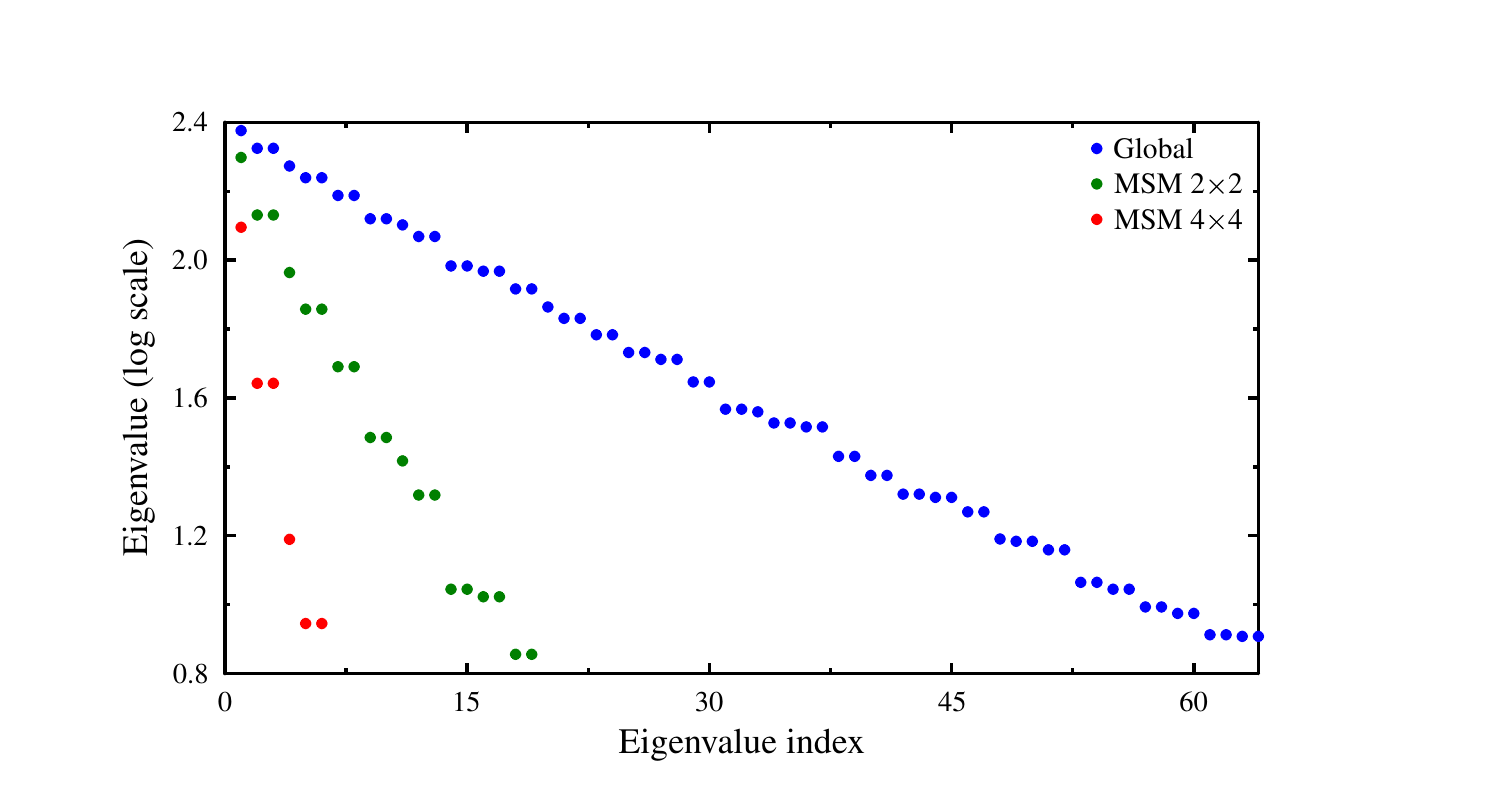}
	\caption{Decay of eigenvalues for the global and multiscale sampling for the third example.}
	\label{eigen_64x64_ch3}
\end{figure}

\begin{figure}[H]
	\centering
	\includegraphics[width= 2.2 in] {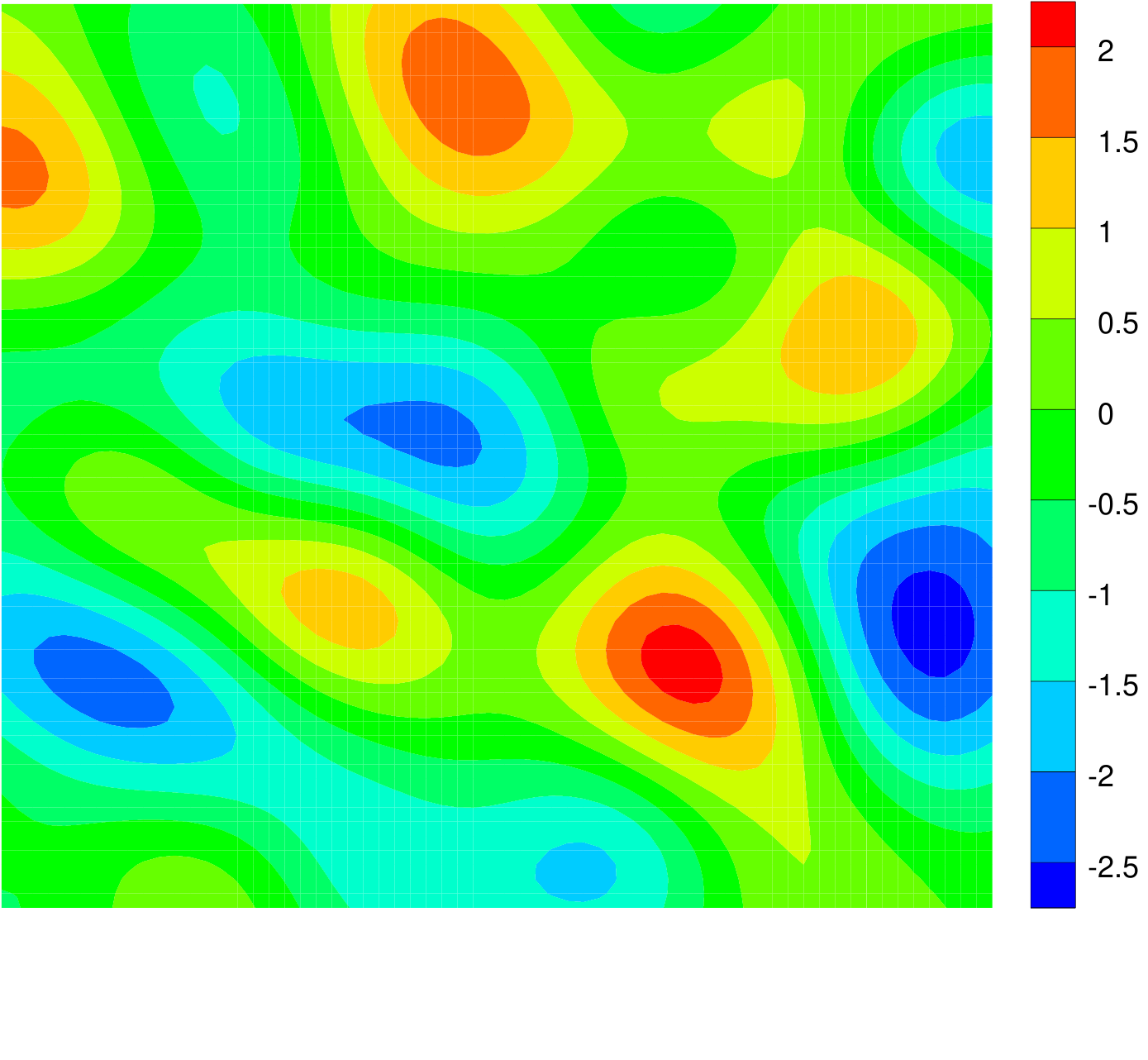}
	\hspace{5mm}
	\includegraphics[width= 2.2 in] {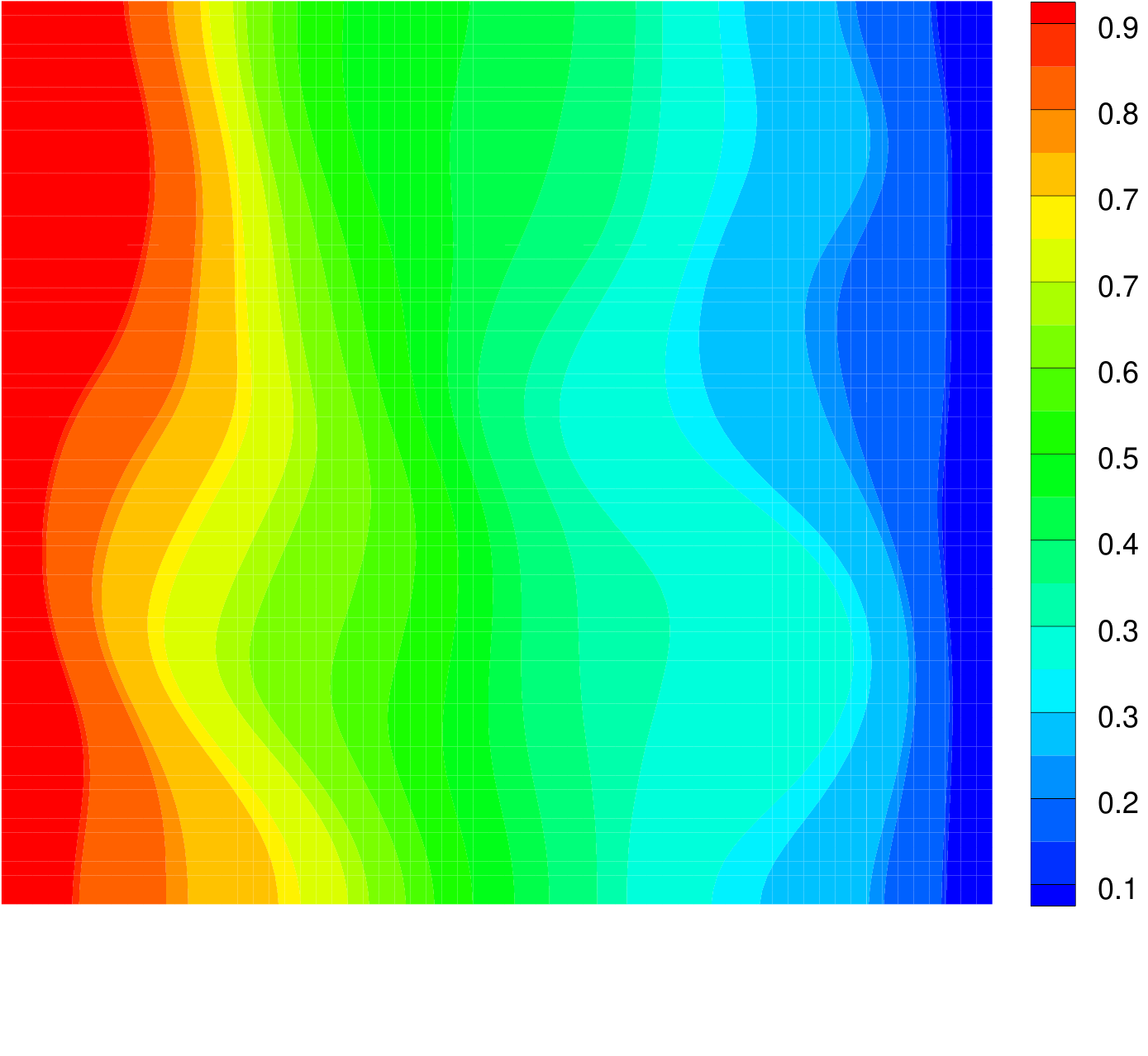}
	
	\caption{Reference log permeability field (left) and the corresponding reference pressure field (right) for the third example.}
	\label{ref_perm_3_1}
\end{figure}

\begin{figure}[H]
	\centering
	\includegraphics[scale = 0.55]{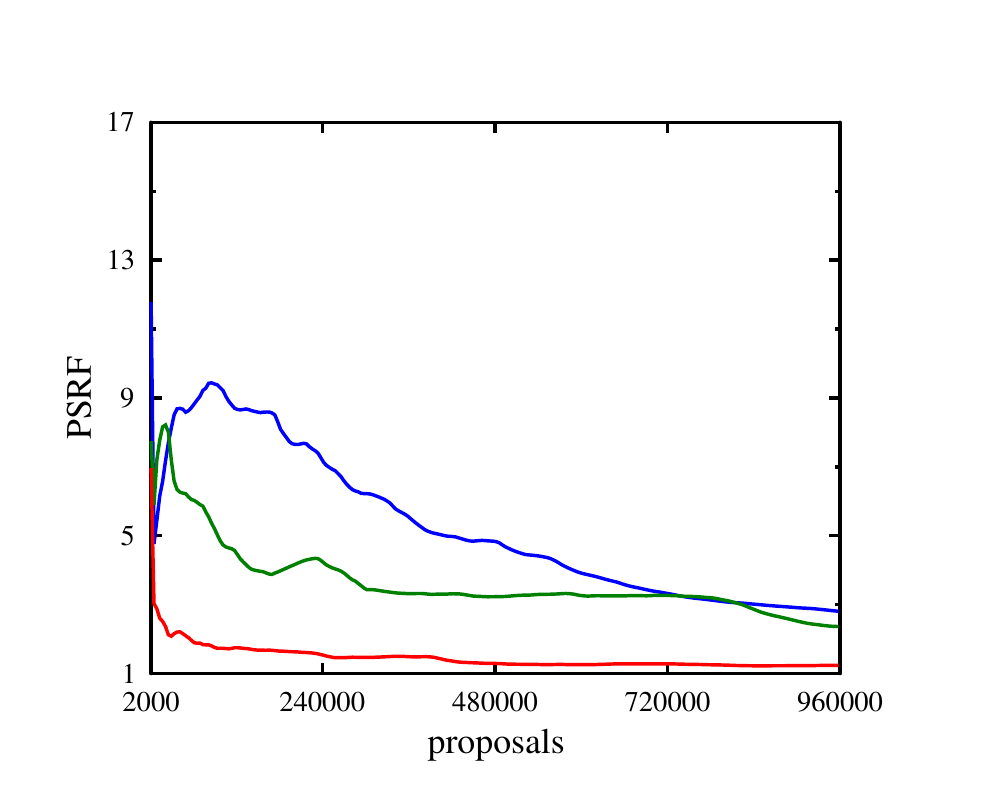}
	\includegraphics[scale= 0.55]{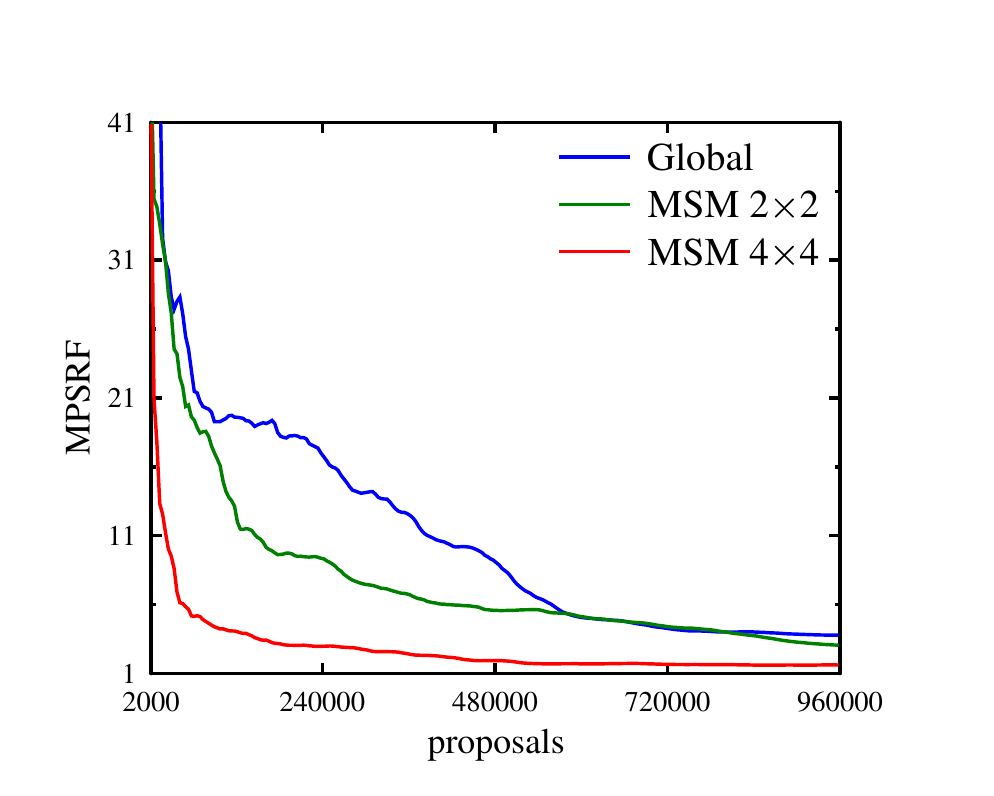}
	\caption{The maximum of PSRFs and MPSRF for the MCMC method with and without multiscale sampling for the third example.}
	\label{MPSRF_64x64_ch3}
\end{figure}

\begin{figure}[H]
	\centering
	\includegraphics[scale = 0.9]{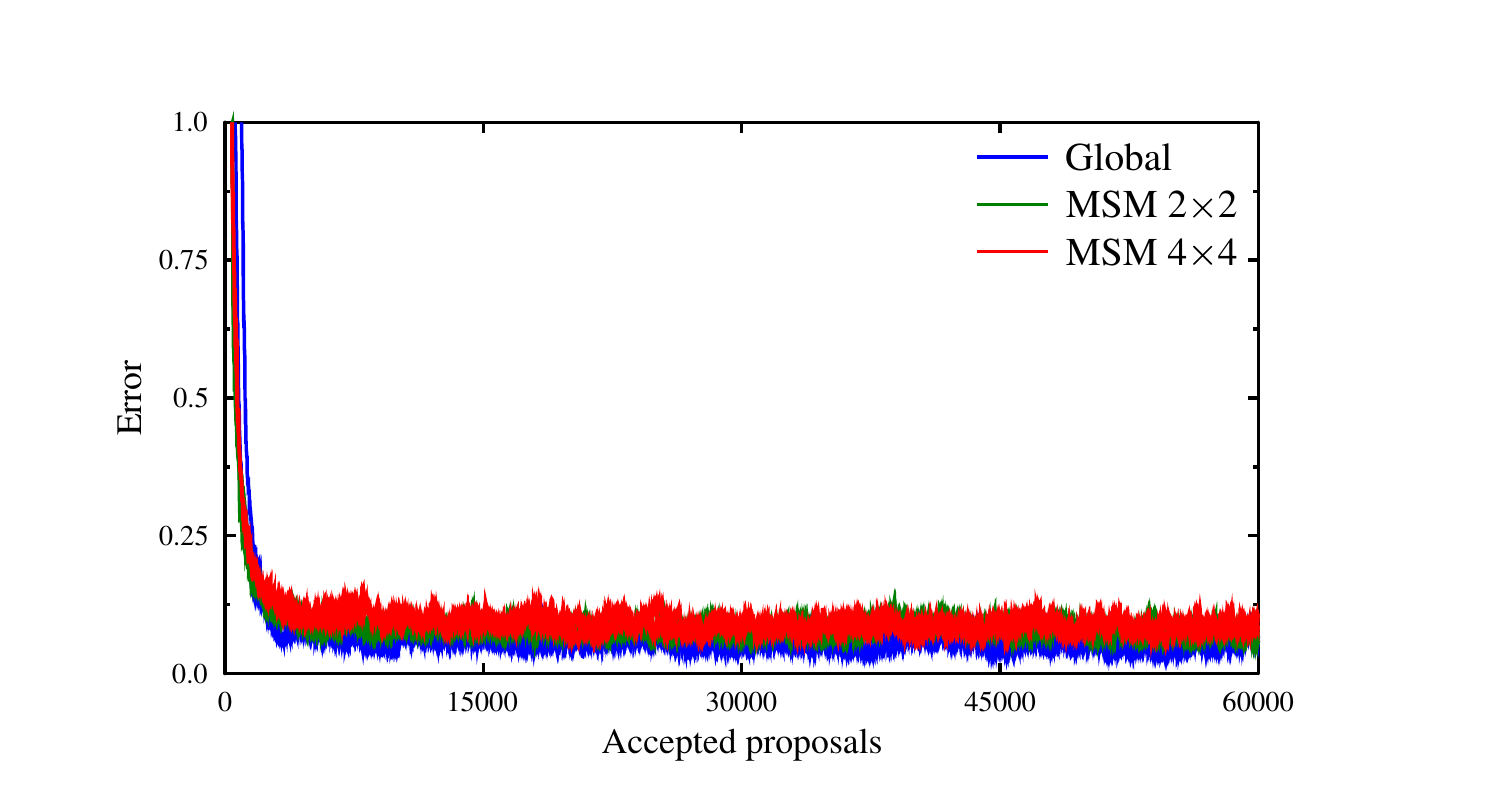}
	\caption{Error curves of the preconditioned MCMC with and without multiscale sampling for the third example.}
	\label{error_64x64_ch3}
\end{figure}

\begin{figure}[H]
	\centering
	\includegraphics[width= 1.5 in] {figures/ex_3_perms/Ref_ex_3_ref.pdf}\\
	\includegraphics[width= 1.5 in] {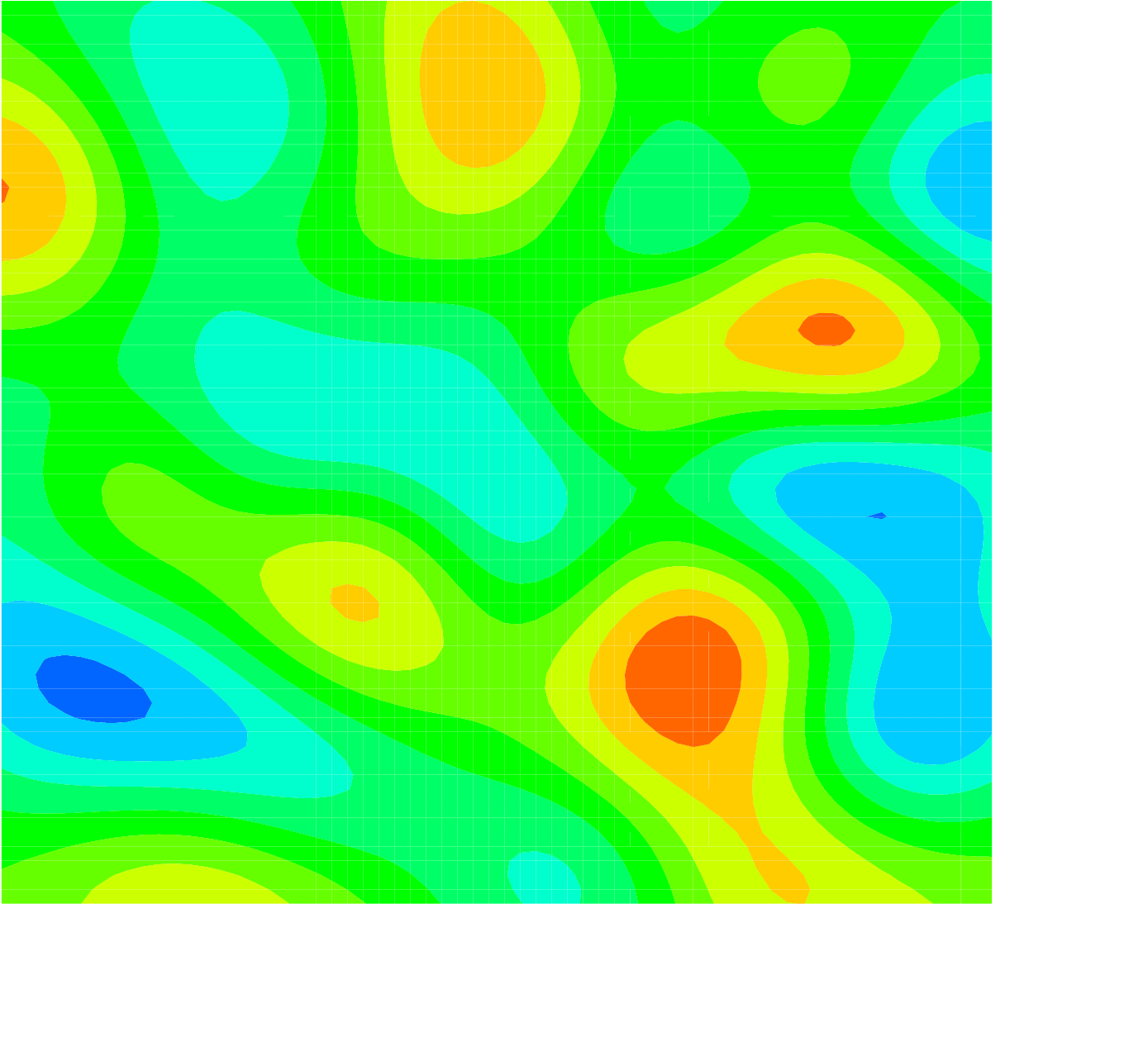}
	\includegraphics[width= 1.5 in] {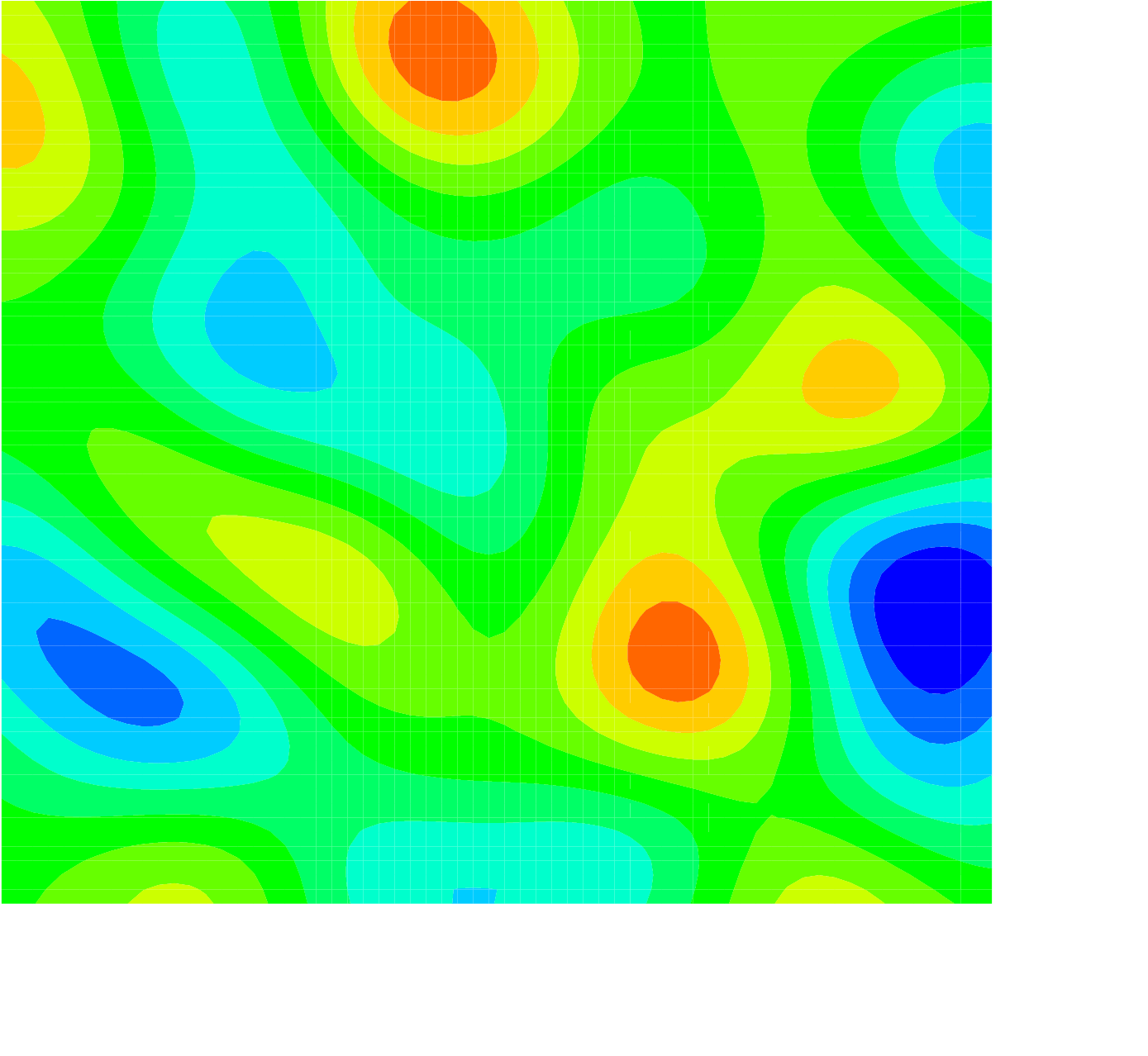}
	\includegraphics[width= 1.5 in] {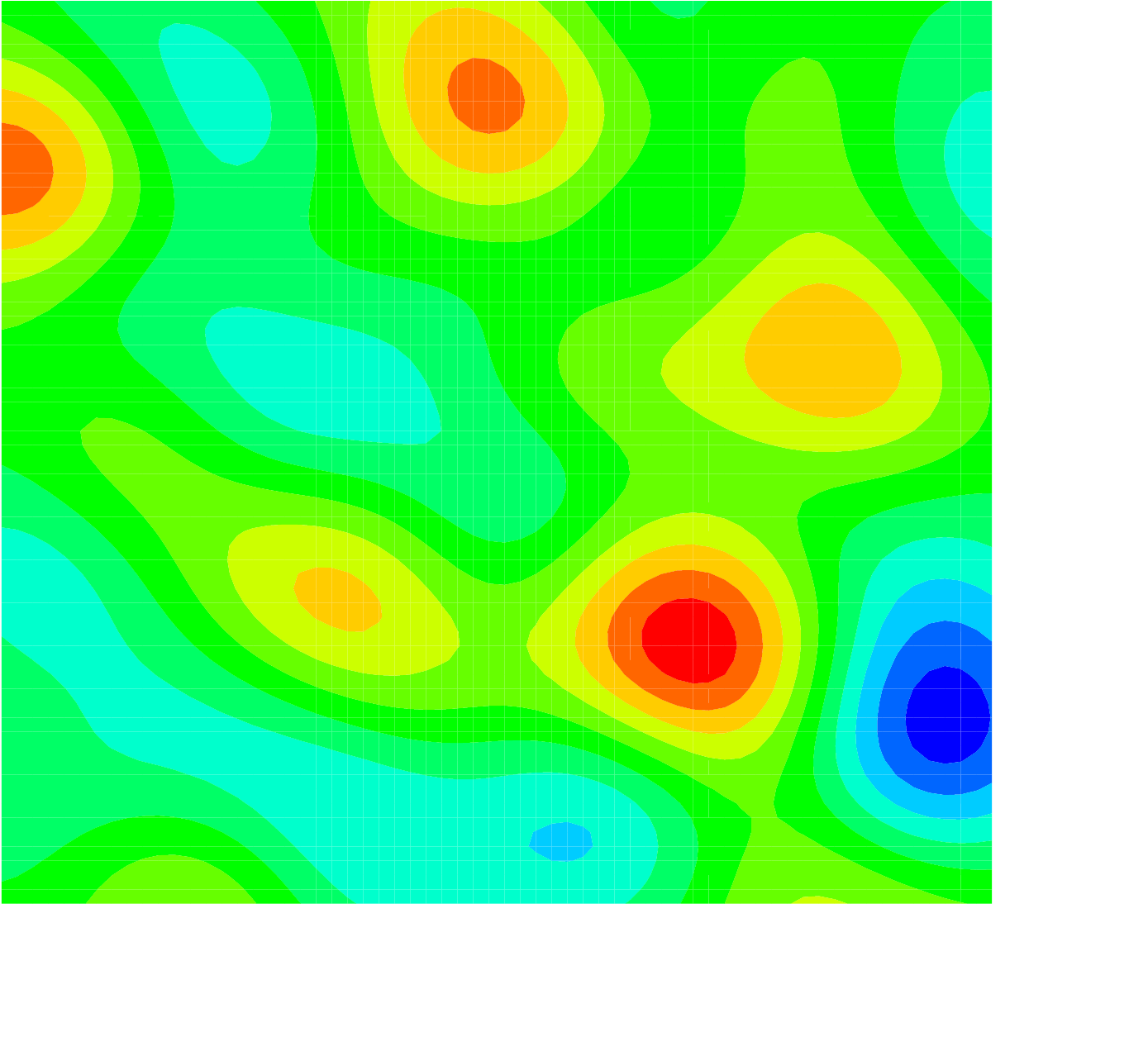}\\
	\includegraphics[width= 1.5 in] {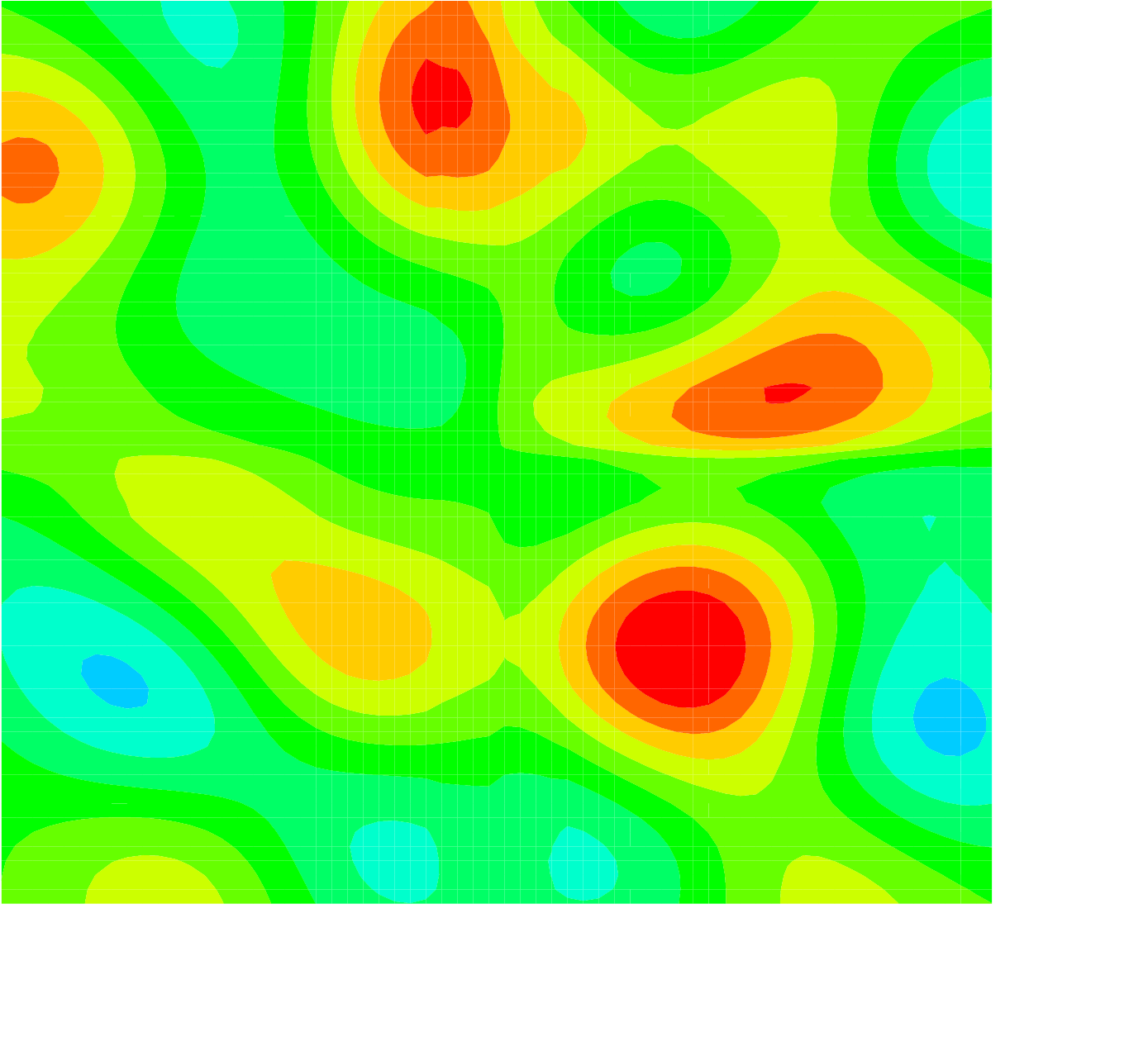}
	\includegraphics[width= 1.5 in] {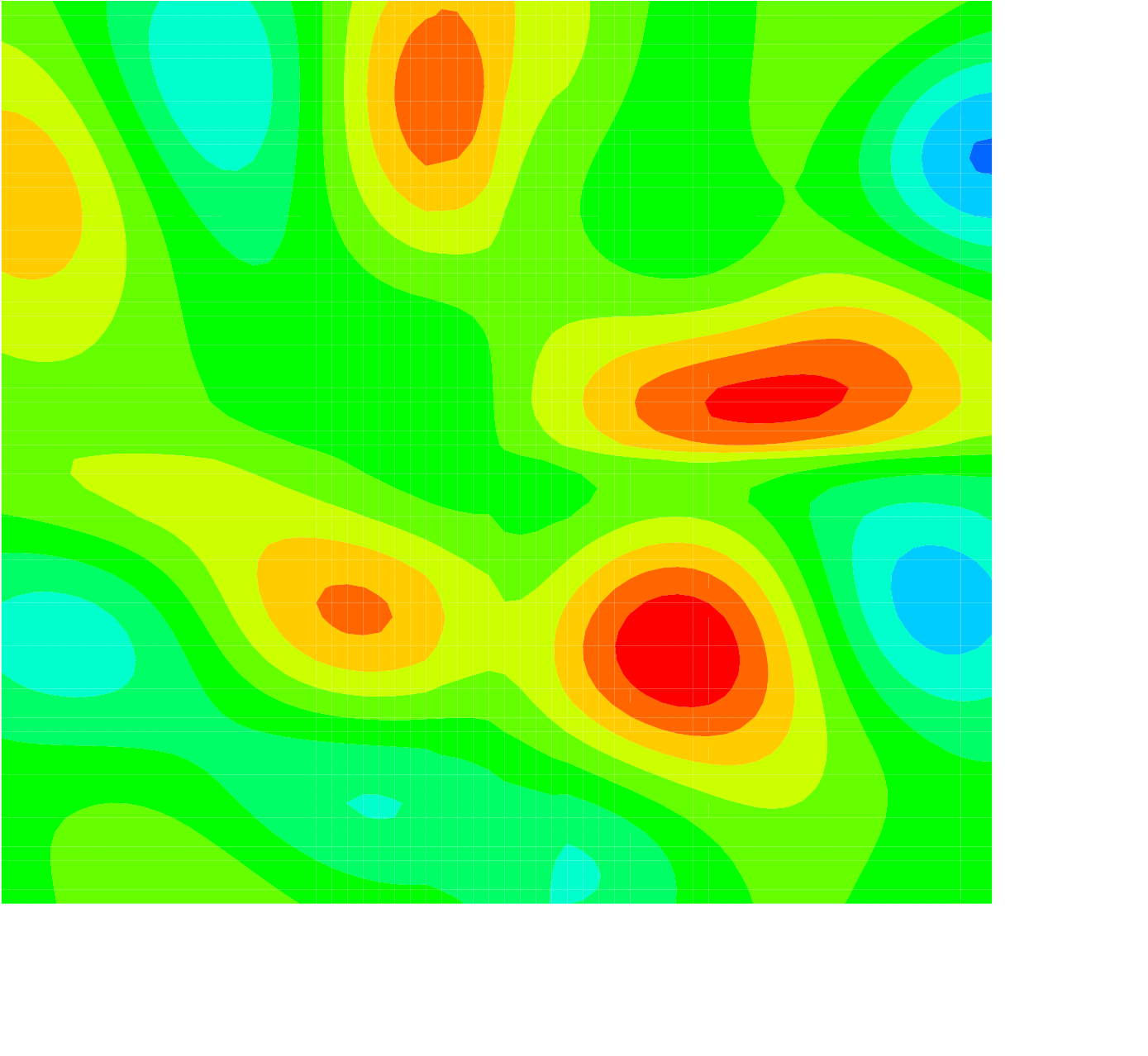}			\includegraphics[width= 1.5 in] {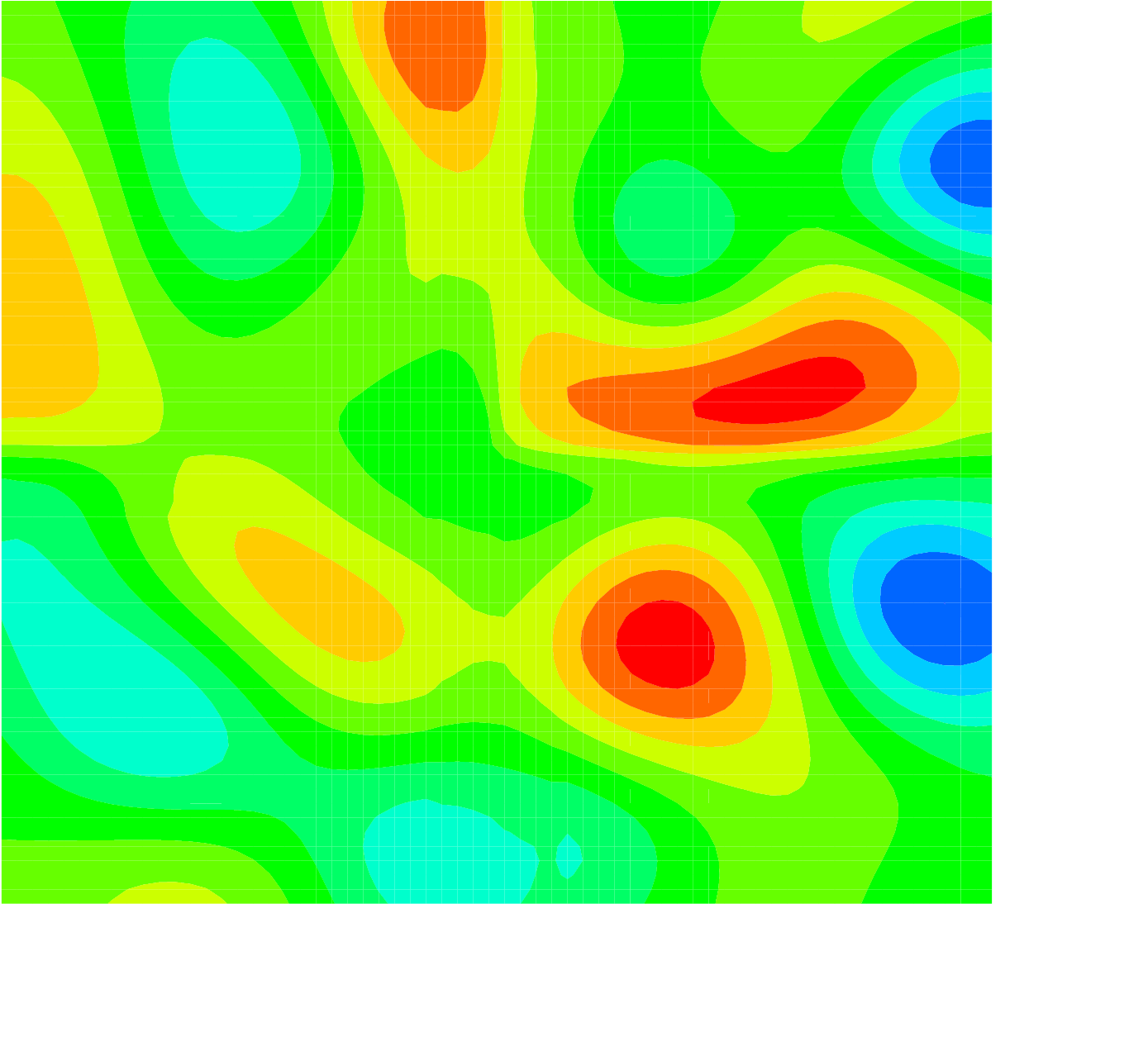}\\
	\includegraphics[width= 1.5 in] {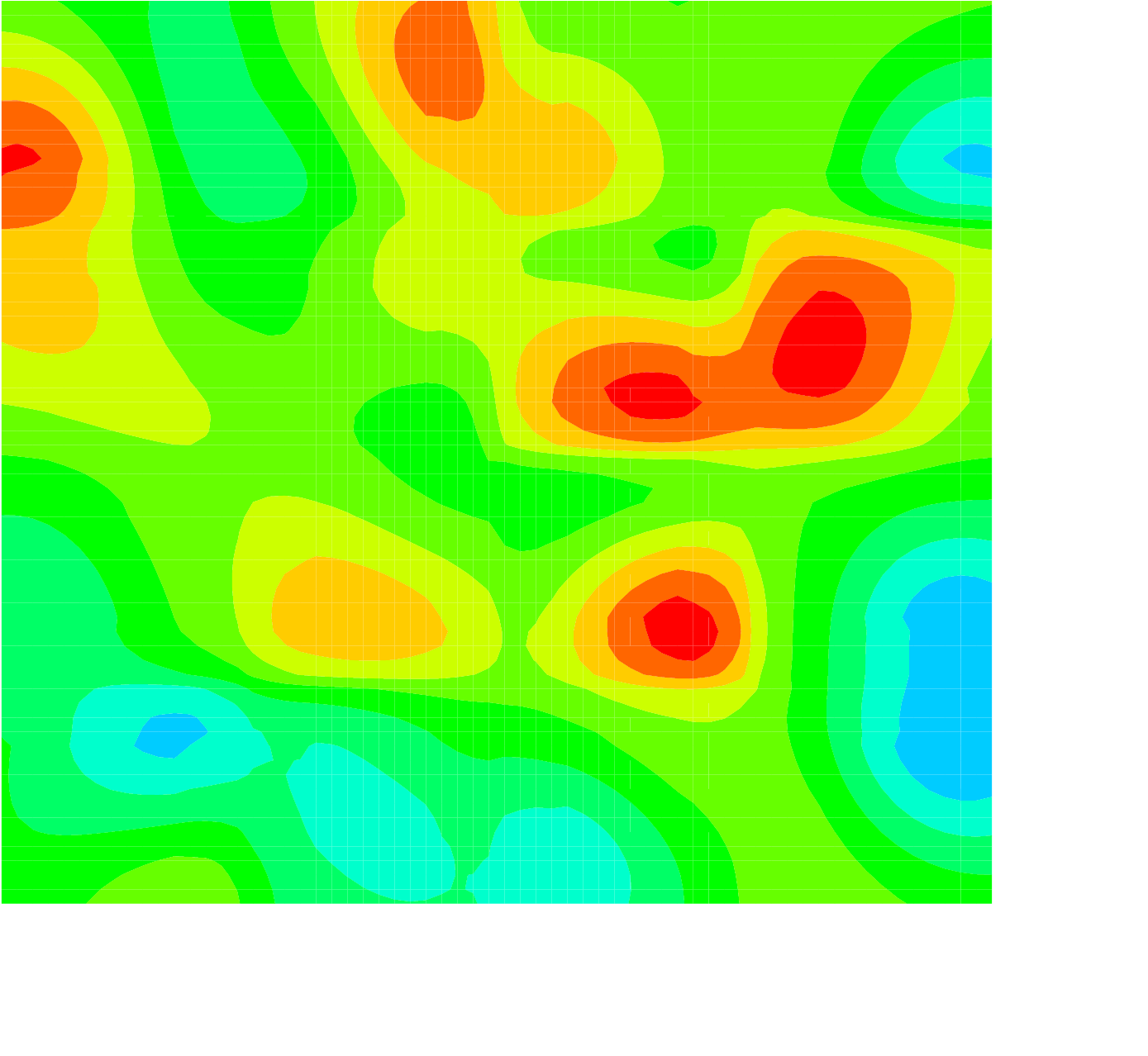}
	\includegraphics[width= 1.5 in] {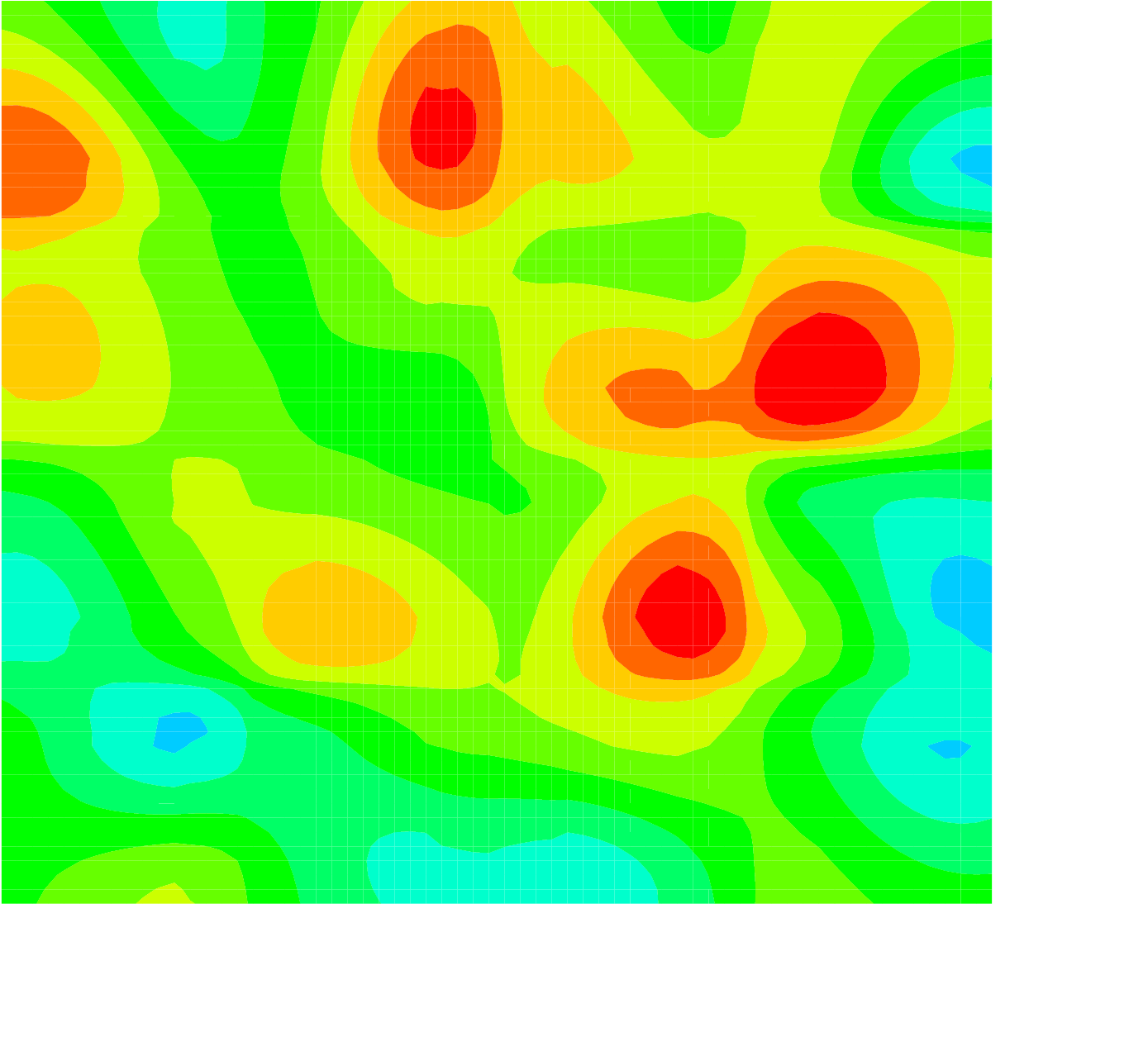}
	\includegraphics[width= 1.5 in] {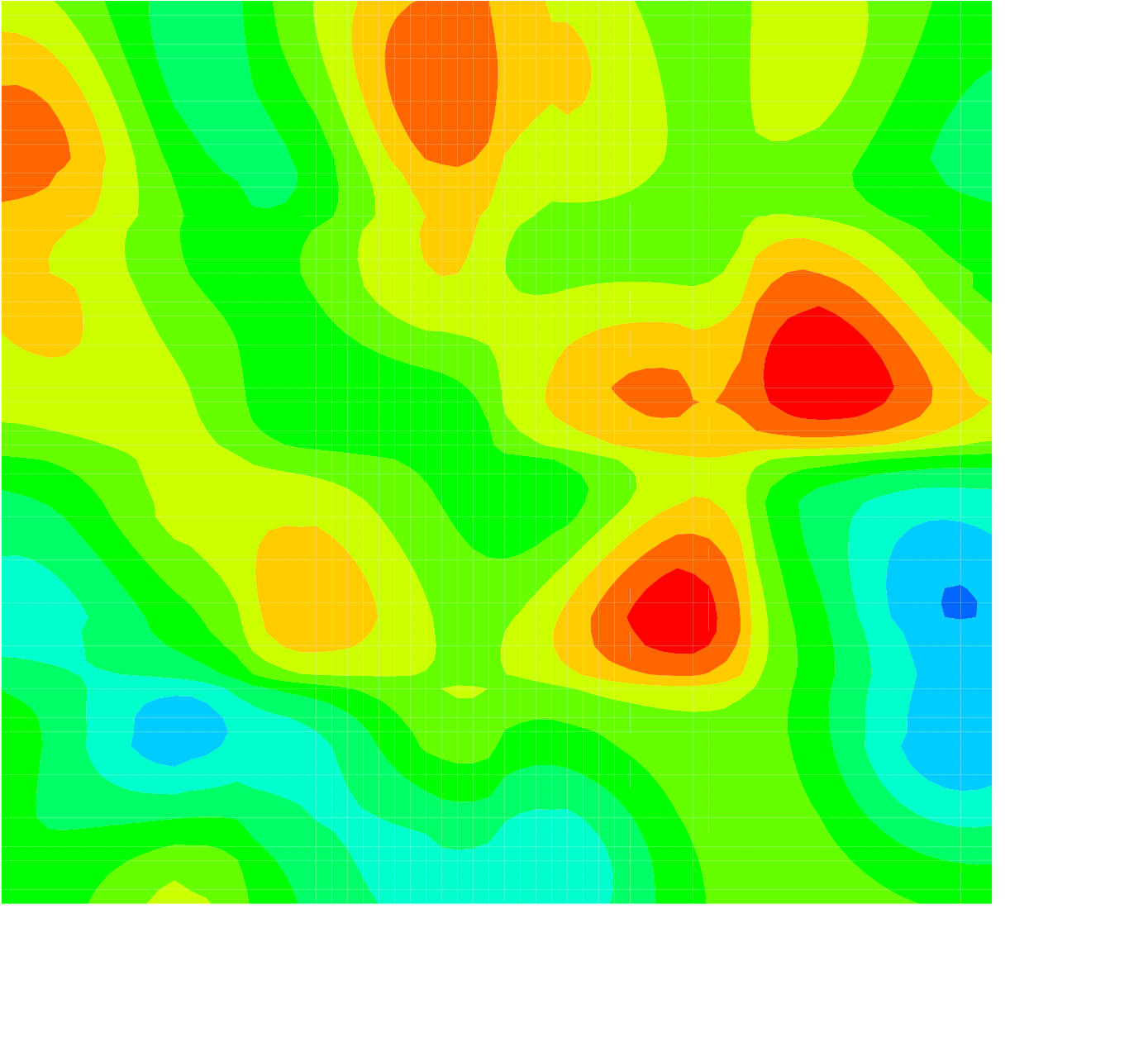}
	
	\caption{First row: Reference log permeability filed. Second row: Accepted permeability fields in the global sampling method. Third row: Accepted permeability fields in MSM $2\times 2$. Fourth row: Accepted permeability fields in MSM $4\times 4$.
	From left to right, log permeability fields at 80000, 160000 and 240000 iterations, respectively, from chain 1 in the third example.}
	\label{perm_3_chain1}
\end{figure}

\begin{figure}[H]
	\centering
	\includegraphics[width= 1.5 in] {figures/ex_3_perms/Ref_ex_3_ref.pdf}\\
	\includegraphics[width= 1.5 in] {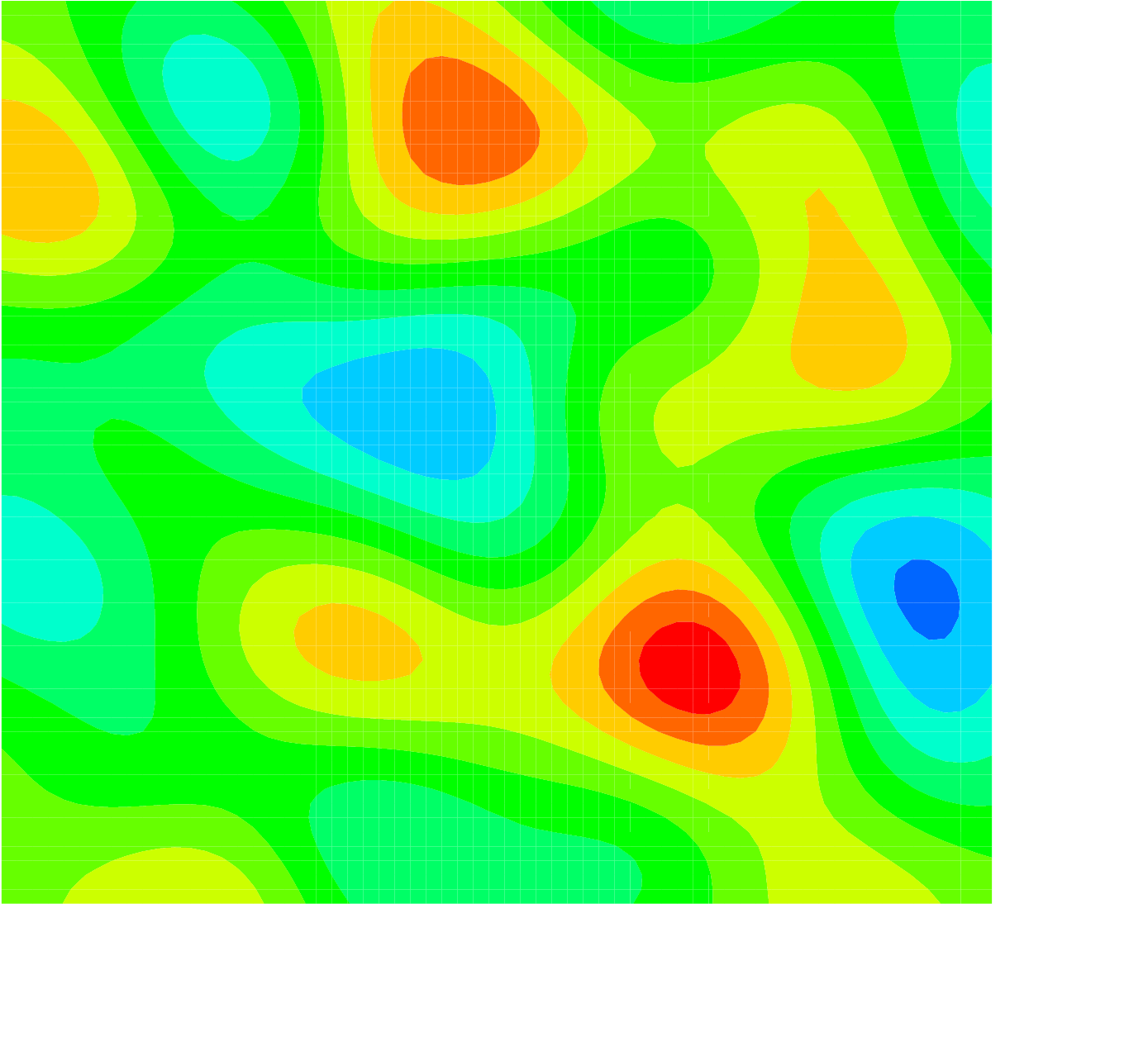}
	\includegraphics[width= 1.5 in] {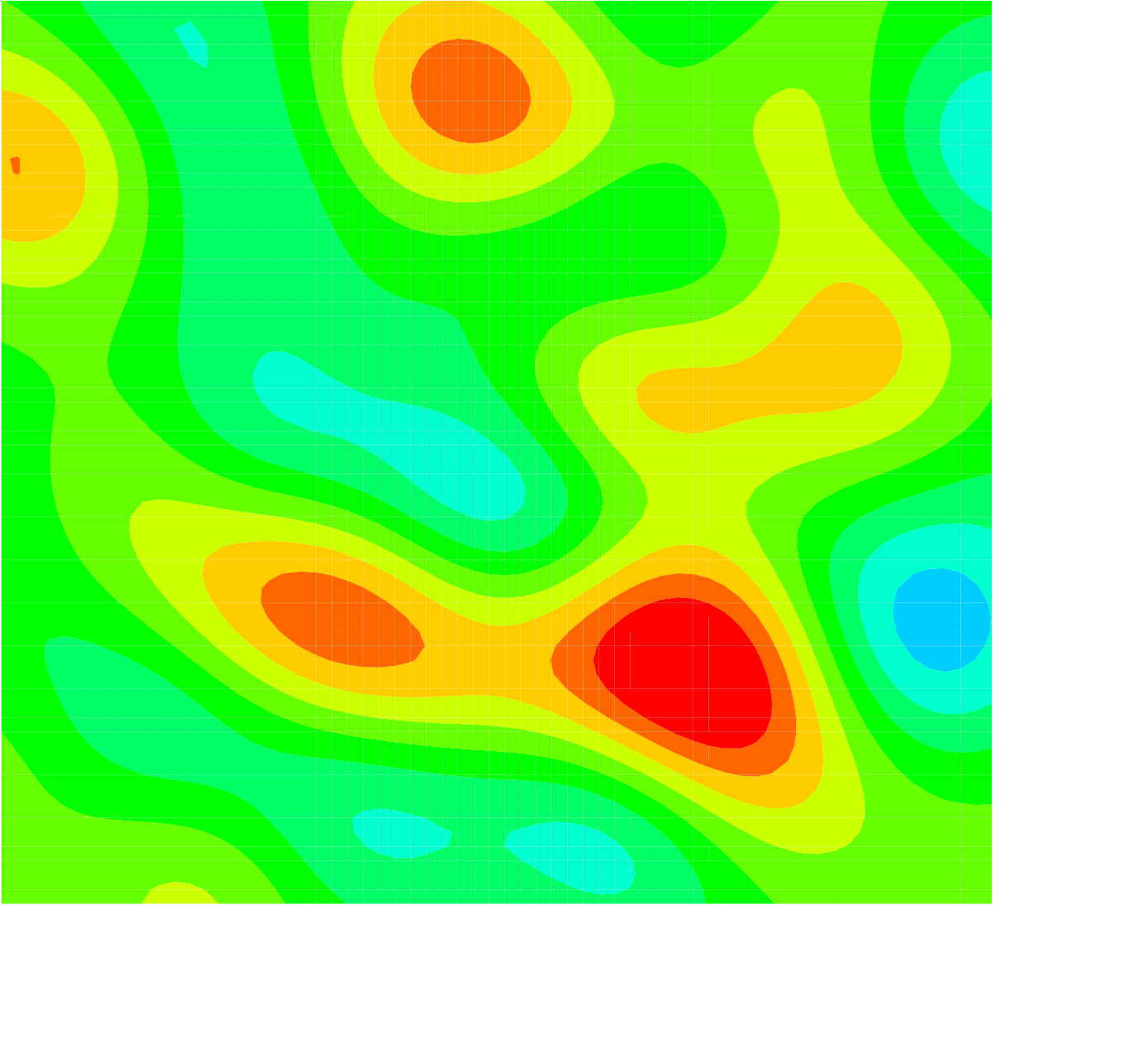}
	\includegraphics[width= 1.5 in] {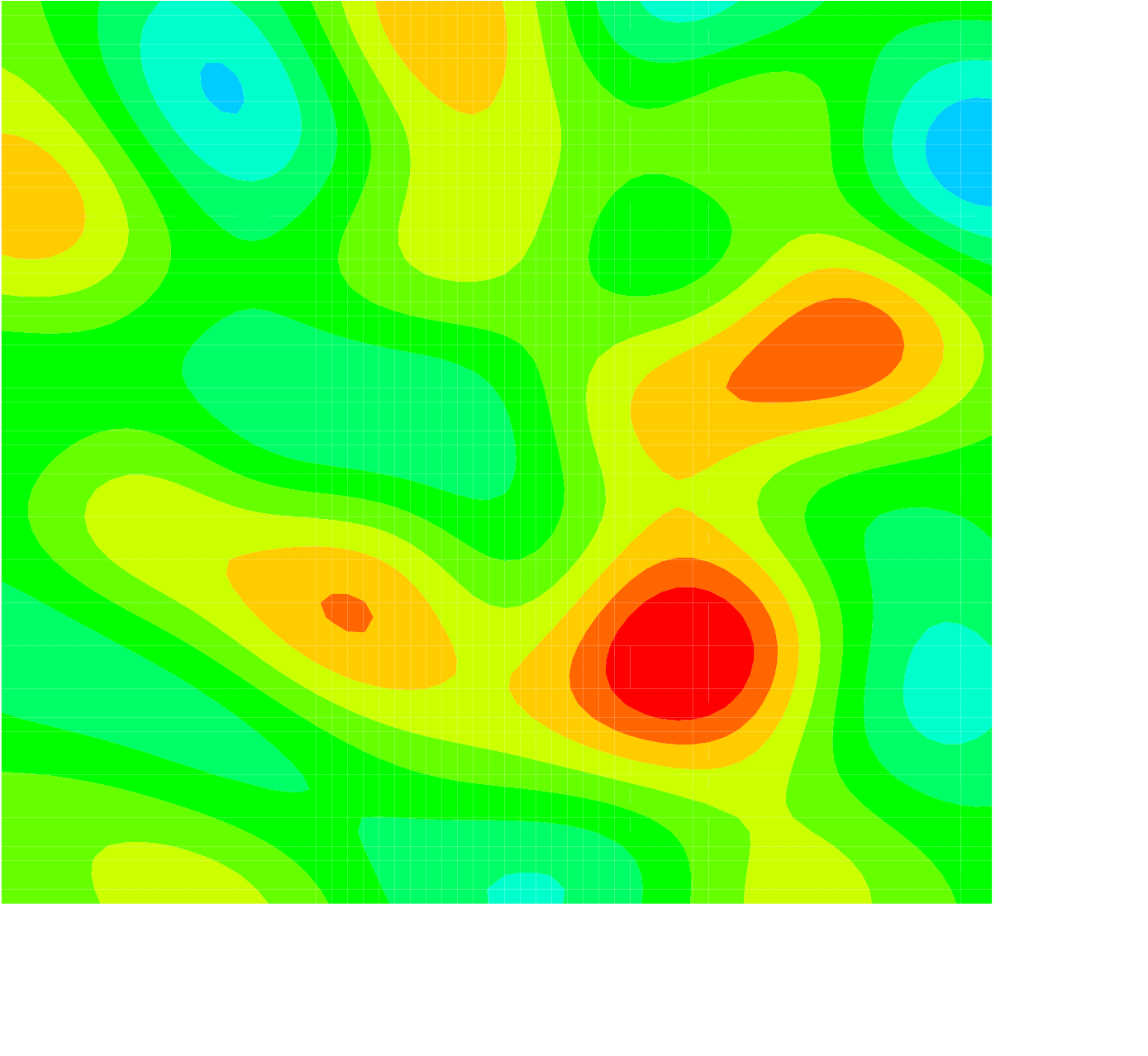}\\
	\includegraphics[width= 1.5 in] {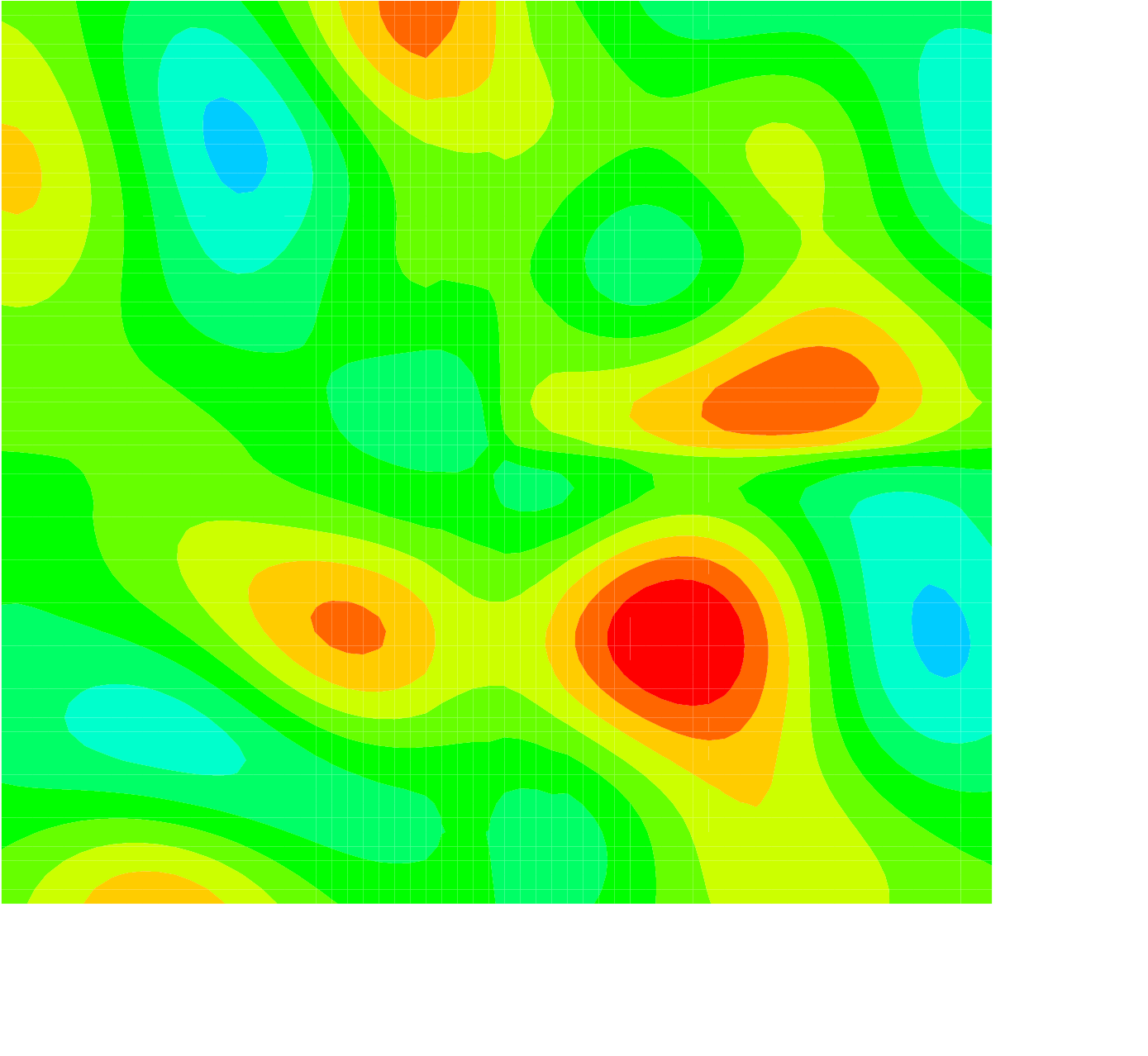}
	\includegraphics[width= 1.5 in] {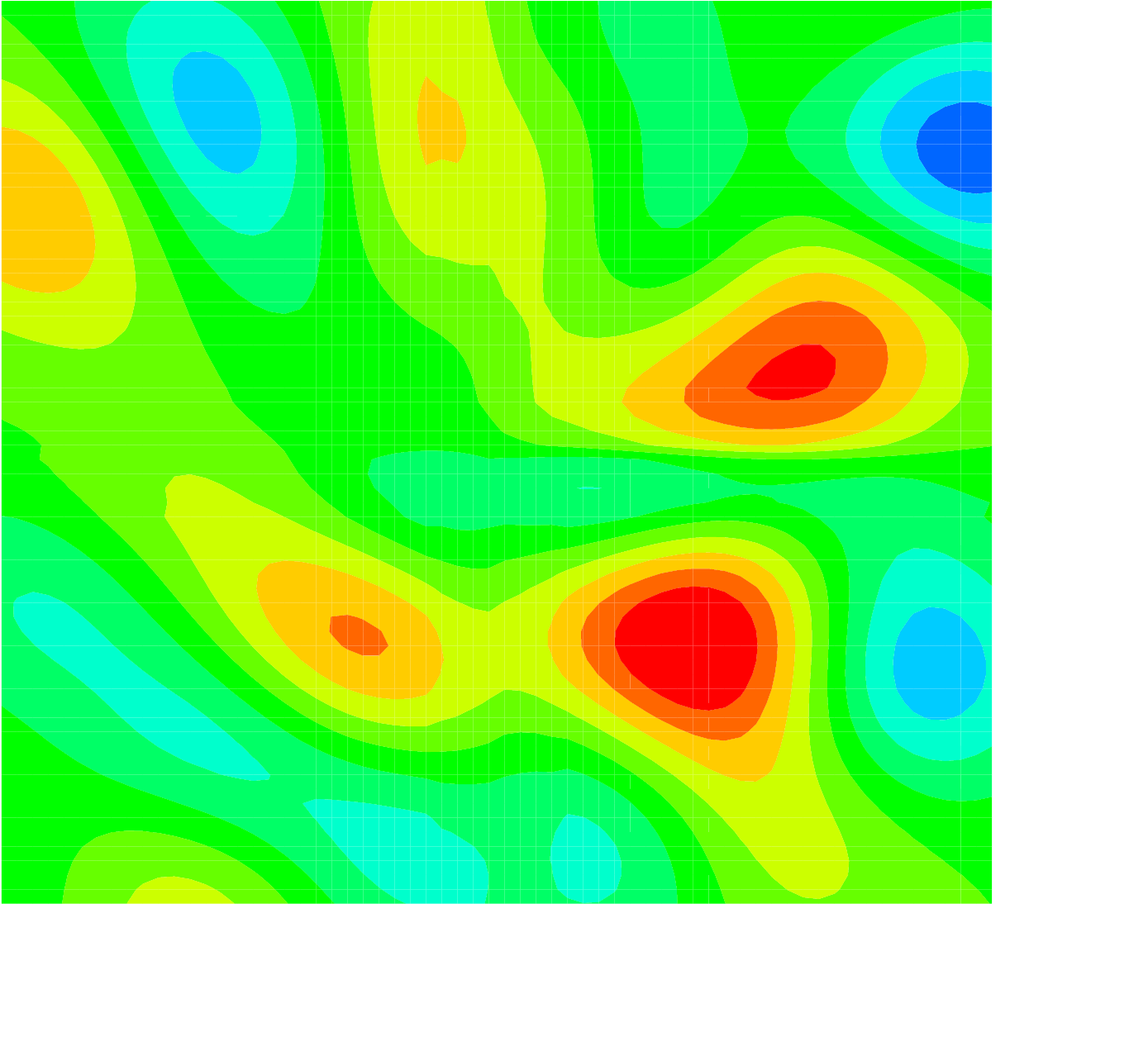}
	\includegraphics[width= 1.5 in] {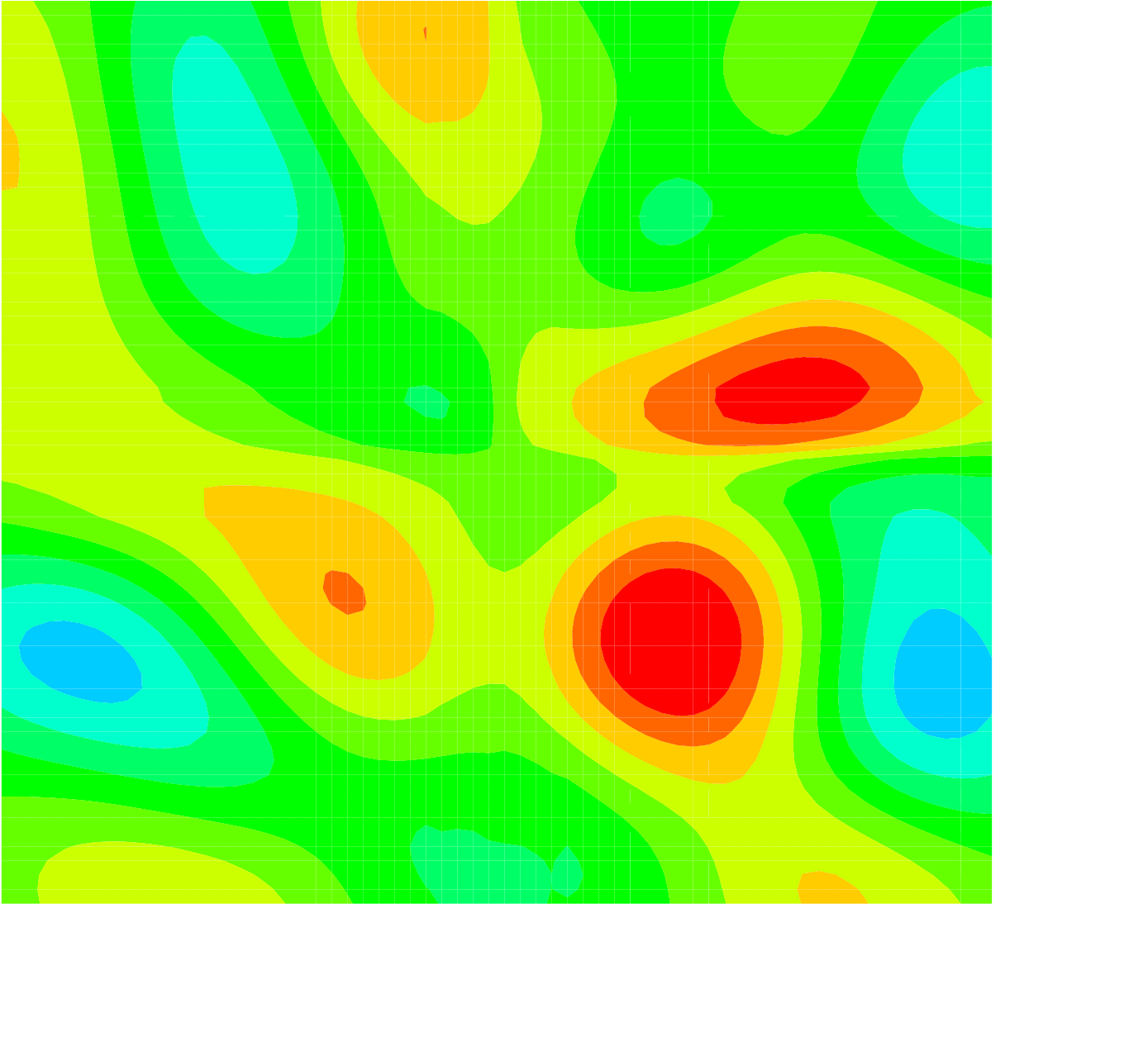}\\
	\includegraphics[width= 1.5 in] {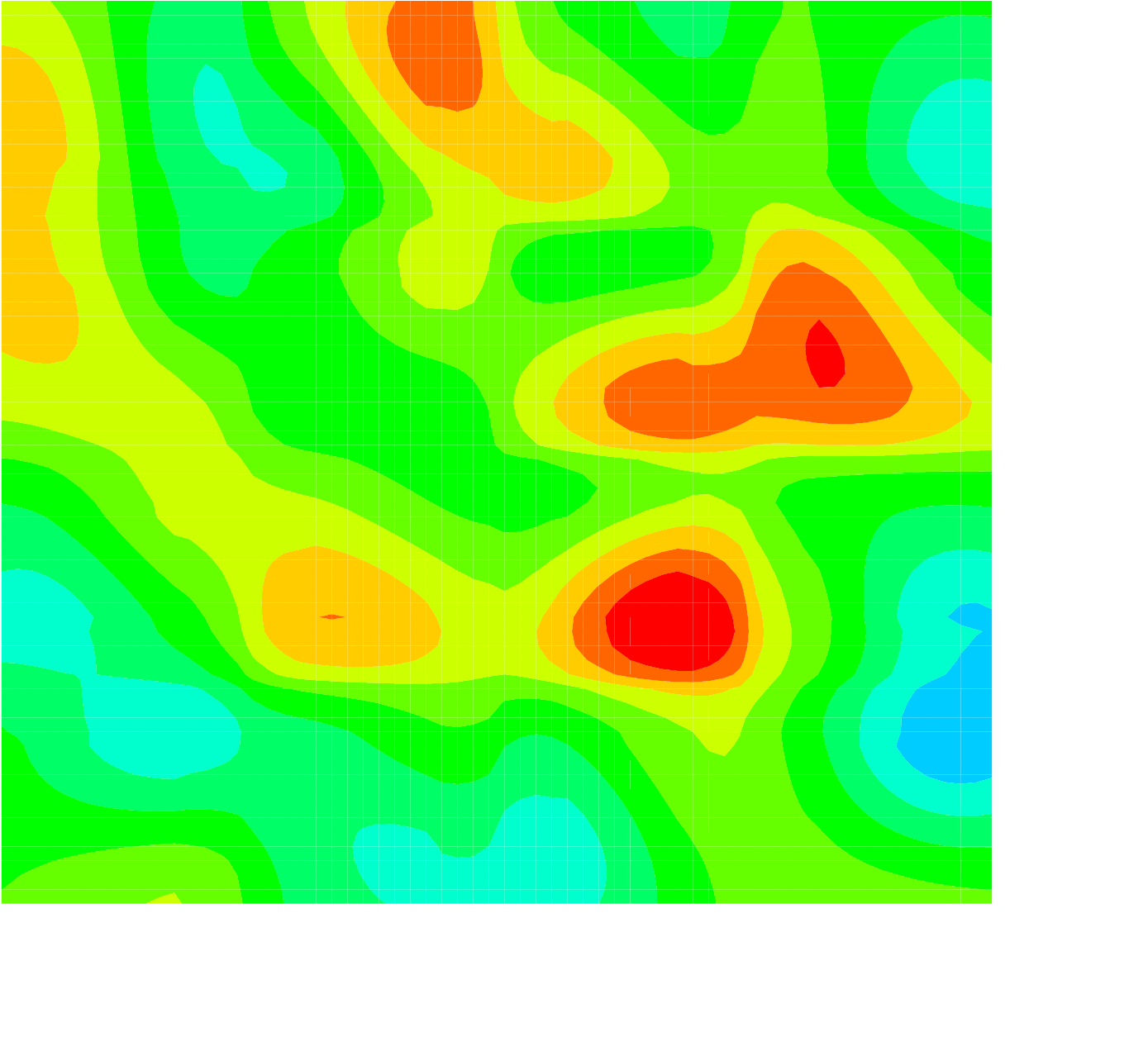}
	\includegraphics[width= 1.5 in] {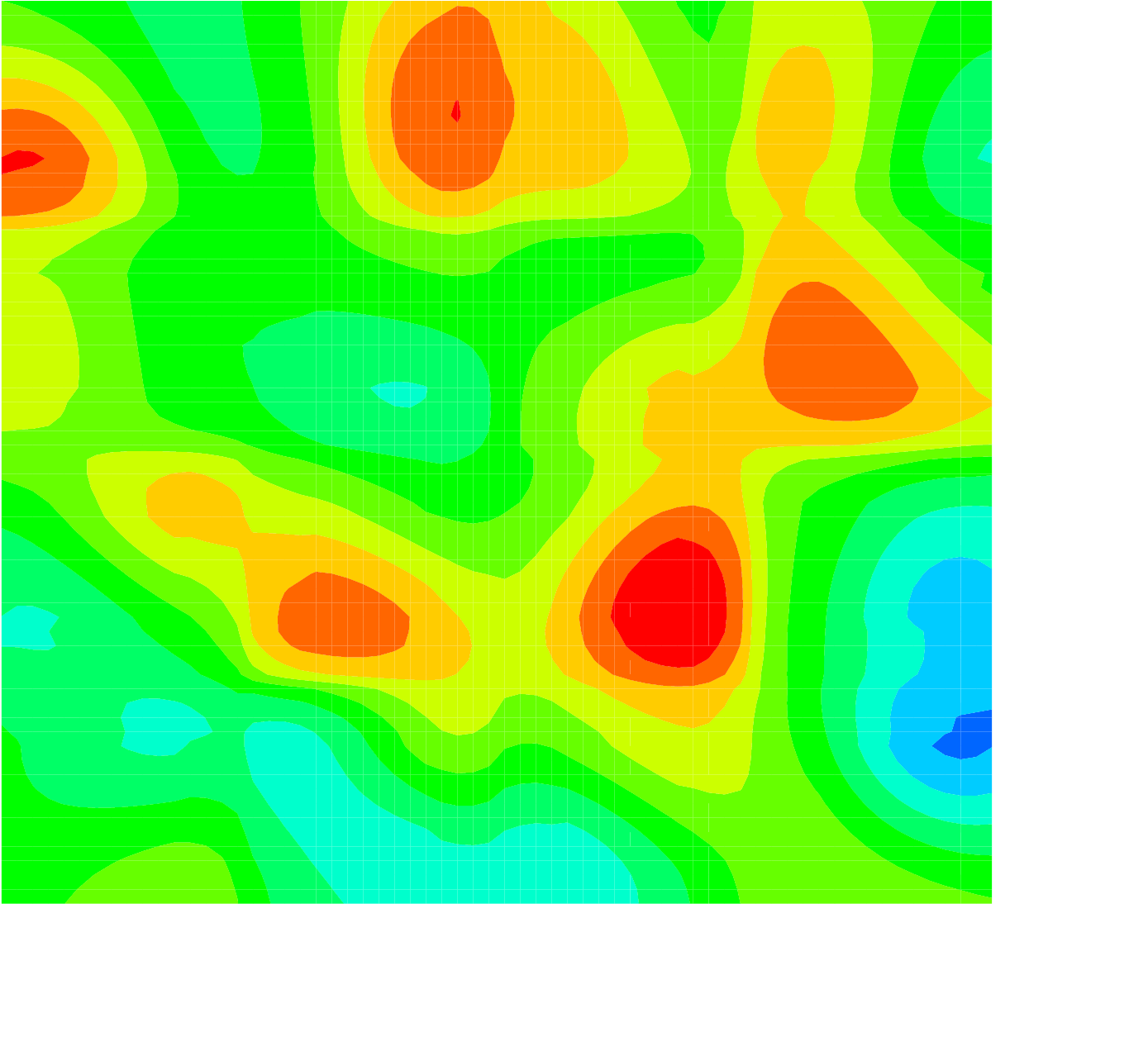}
   \includegraphics[width= 1.5 in] {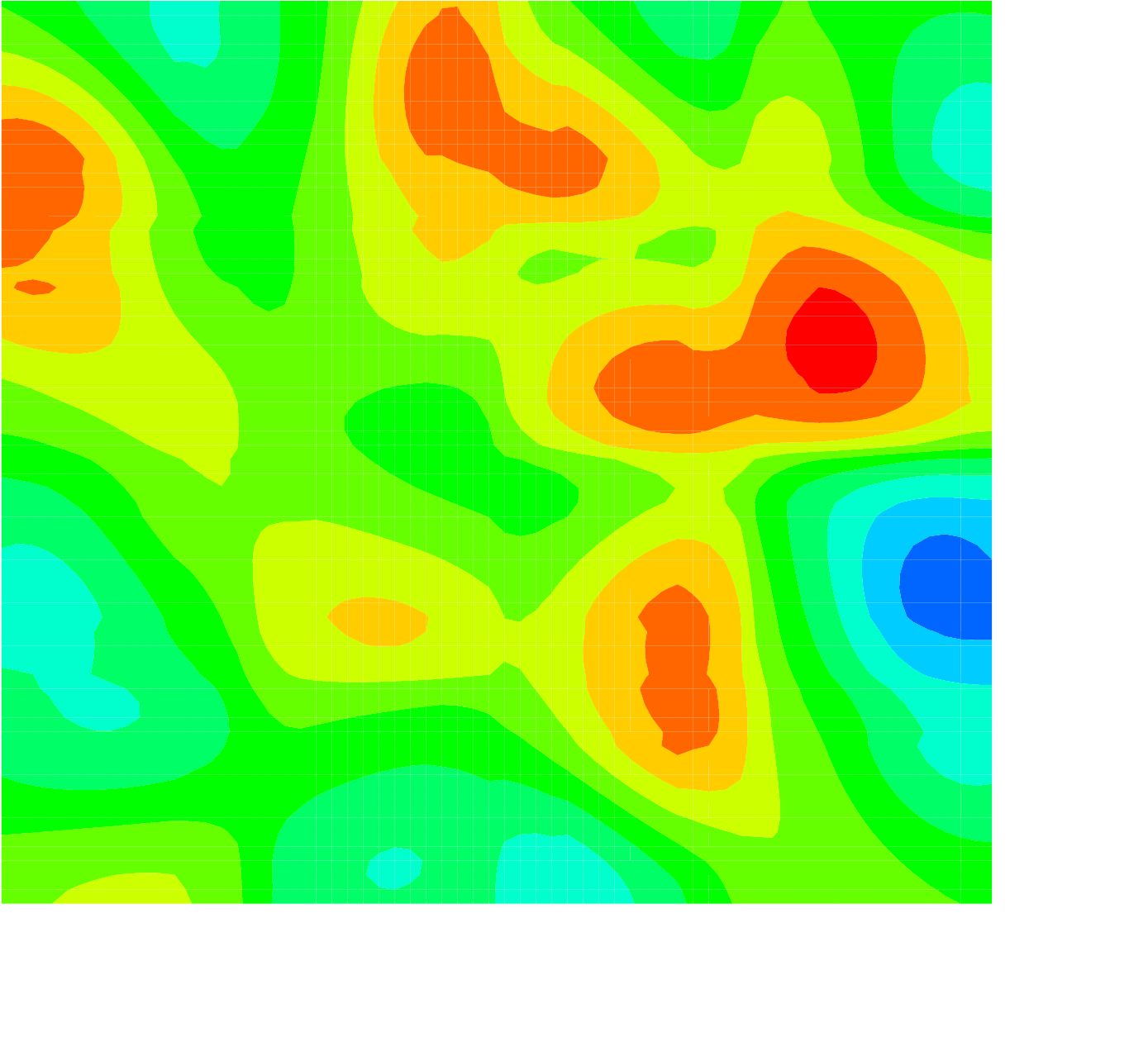}
	
	\caption{First row: Reference log permeability filed. Second row: Accepted permeability fields in the global sampling method. Third row: Accepted permeability fields in MSM $2\times 2$. Fourth row: Accepted permeability fields in MSM $4\times 4$.
	From left to right, log permeability fields at 80000, 160000 and 240000 iterations, respectively, from chain 2 in the third example.}
	\label{perm_3_chain2}
\end{figure}

Let us consider PSRFs and MPSRF curves for these methods. We take $240000$ proposals from each chain in constructing the PSRF and MPSRF curves. We show the maximum of PSRFs and MPSRF curves in Figure  \ref{MPSRF_64x64_ch3}. For MSM $4\times 4$, the values at the tails of the PSRF and MPSRF curves are $1.2$ and $1.6$, respectively. These values indicate that the curves in MSM $4\times 4$ are closer to the convergence. However, the curves in MSM $2\times 2$ and global sampling method are very far from reaching a convergence. Table \ref{Lx_Ly_3} shows that the acceptance rates for the multiscale sampling methods are also slightly better than that of the global sampling method. The error curves are comparable for these methods. See Figure  \ref{error_64x64_ch3}. Figures \ref{perm_3_chain1} and \ref{perm_3_chain2} compare the accepted permeability fields for two selected chains in the MCMC simulation. Both MSM $2\times 2$ and MSM $4\times 4$ recover the fields better than the global sampling method.

\begin{table}[H] 
	\caption{A comparison of acceptance rates for the MCMC with and without MSM for the third example.}
	\center
	\begin{tabular}{|cccc|}
		\hline
		&  \quad  MCMC global & \quad MCMC with MSM $2\times 2$ & \quad MCMC with MSM $4\times 4$  \\
		\hline
		$\sigma_F^2$ &   $10^{-3}$  & $10^{-3}$  & $10^{-3}$   \\               
		$\sigma_C^2$   & $5\times10^{-3}$ & $5\times10^{-3}$&  $5\times10^{-3}$ \\  
		acc. rate & $41\%$   & $43\%$& $43\%$  \\
		\hline
	\end{tabular}
	\label{Lx_Ly_3}     
\end{table}


\subsection{Example 4}

In this example we compare MSM $4\times4$ with and without conditioning for the problem in example 3. However, in the conditioning approach, we incorporate the permeability measurements at eight sparse locations in the field. 

\begin{figure}[H]
	\centering
	\includegraphics[scale = 0.55]{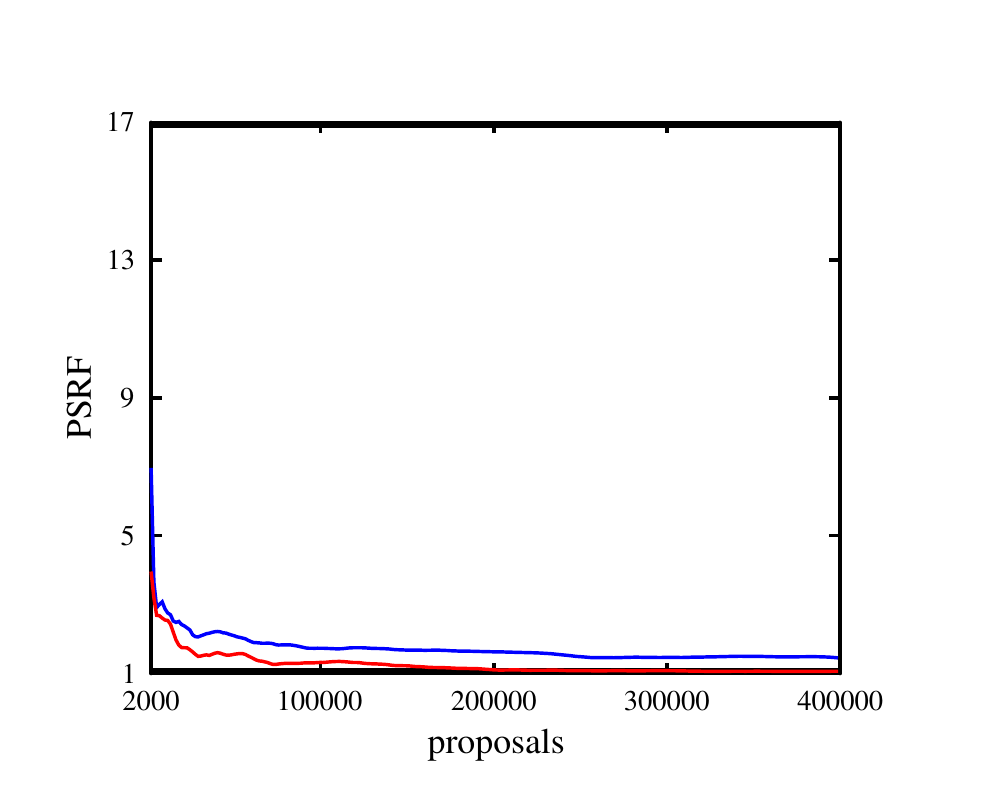}
	\includegraphics[scale= 0.55]{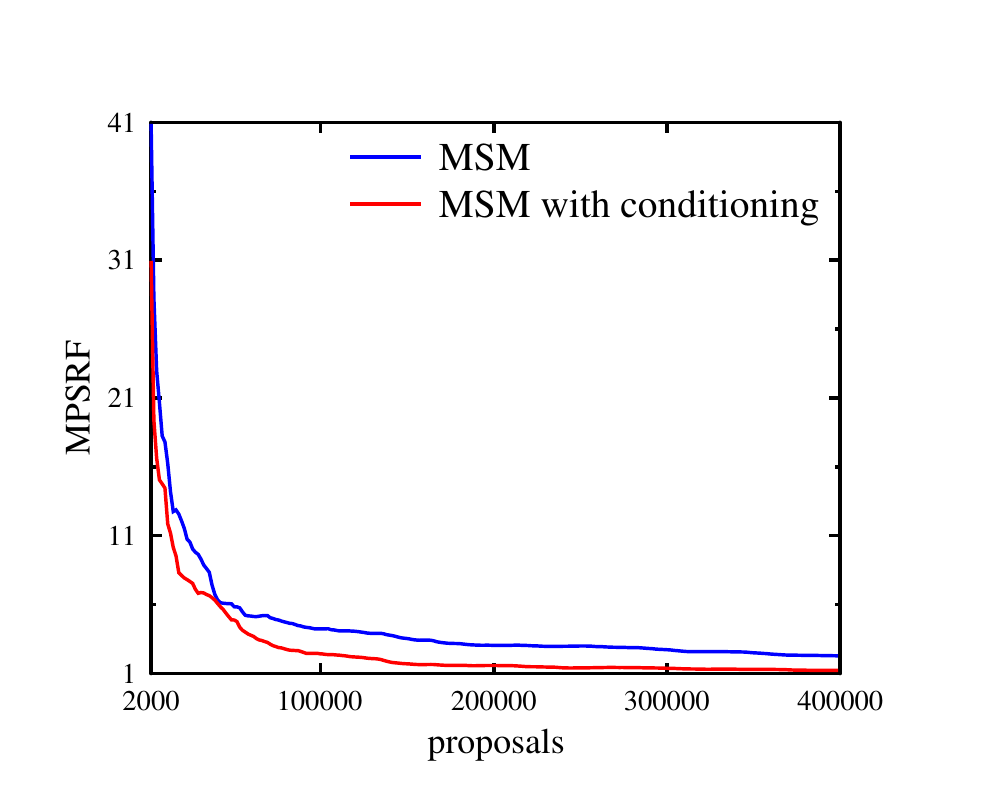}
	\caption{The maximum of PSRFs and MPSRF for the multiscale sampling MCMC method with and without conditioning for the fourth example.}
	\label{conver_cond}
\end{figure}

Four MCMCs are simulated in each case. We compute the maximum of the PSRFs and MPSRF by taking $100,000$ samples from each chain (total $400,000$ samples). In Figure  \ref{conver_cond} we present the maximum of PSRFs and MPSRF curves. The values in the tails of the maximum of the PSRFs and the MPSRF are $1.09$ and $1.17$, respectively, for the method with conditioning. For the method without conditioning, these values are $2.5$ and $3.5$, respectively. Therefore, we can say that MSM $4\times4$ with conditioning reaches convergence earlier than MSM without conditioning.

\begin{figure}[H]
	\centering
	\includegraphics[width= 1.5in] {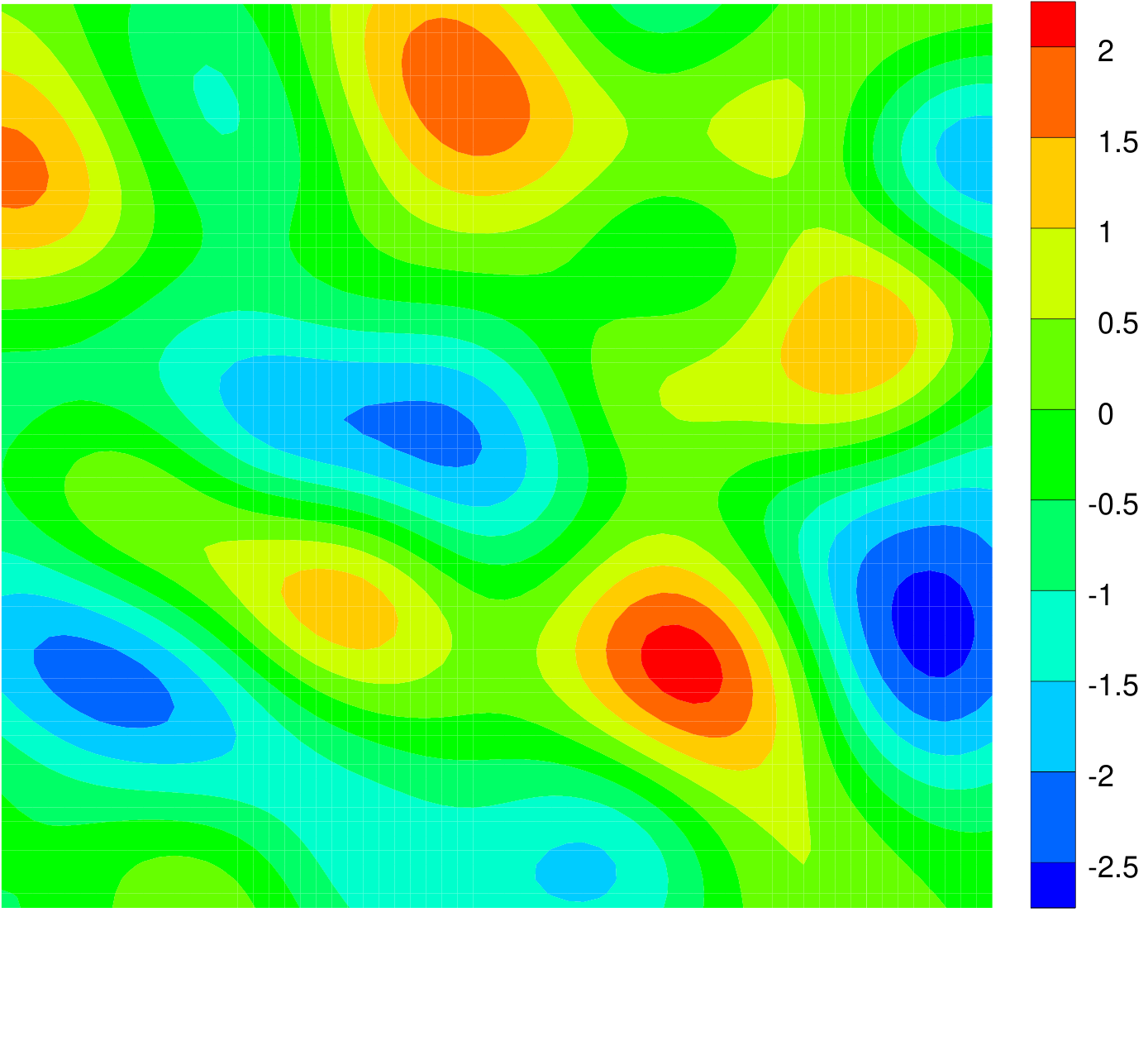}\\
		\includegraphics[width= 1.5 in] {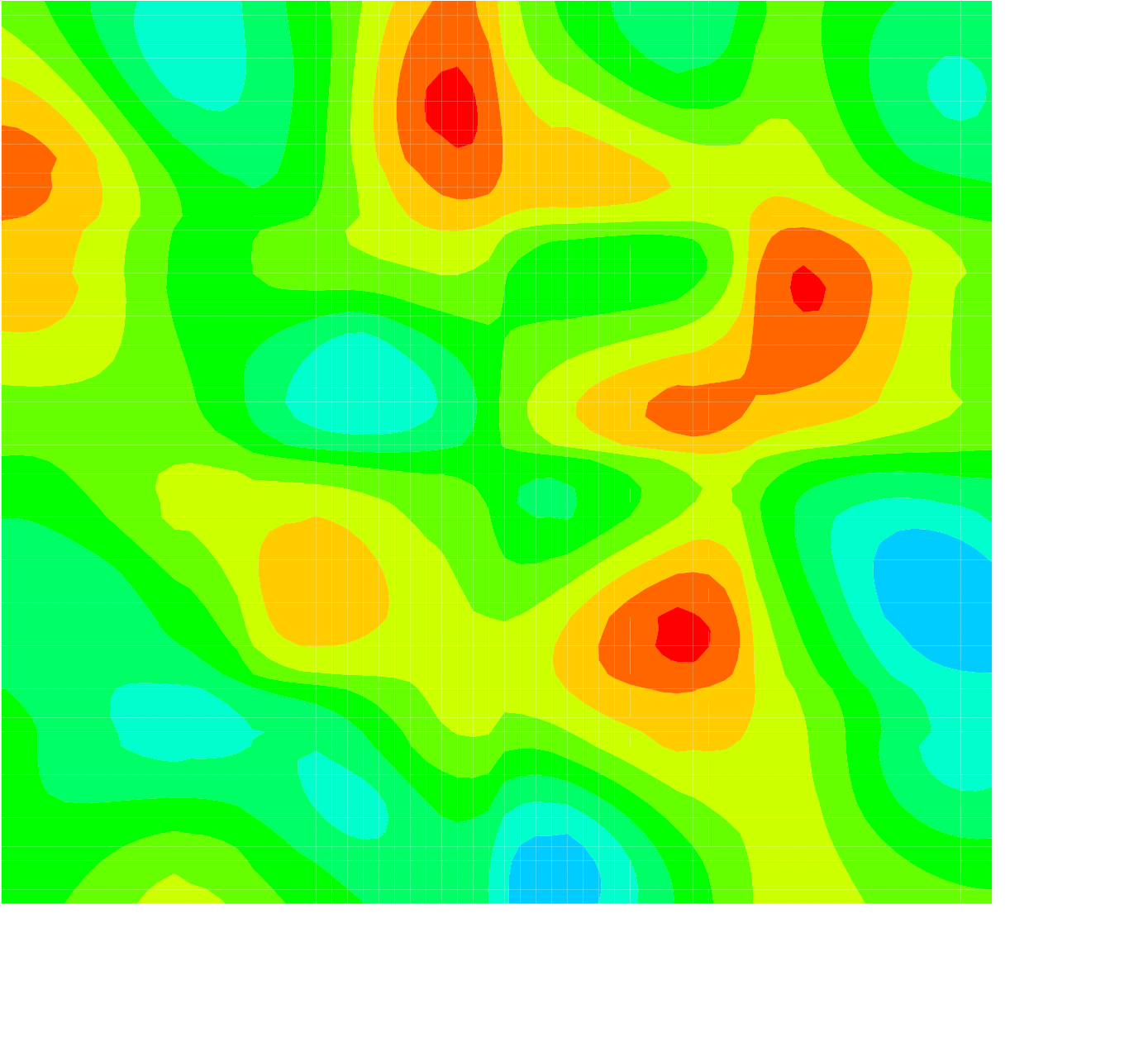}
		\includegraphics[width= 1.5 in] {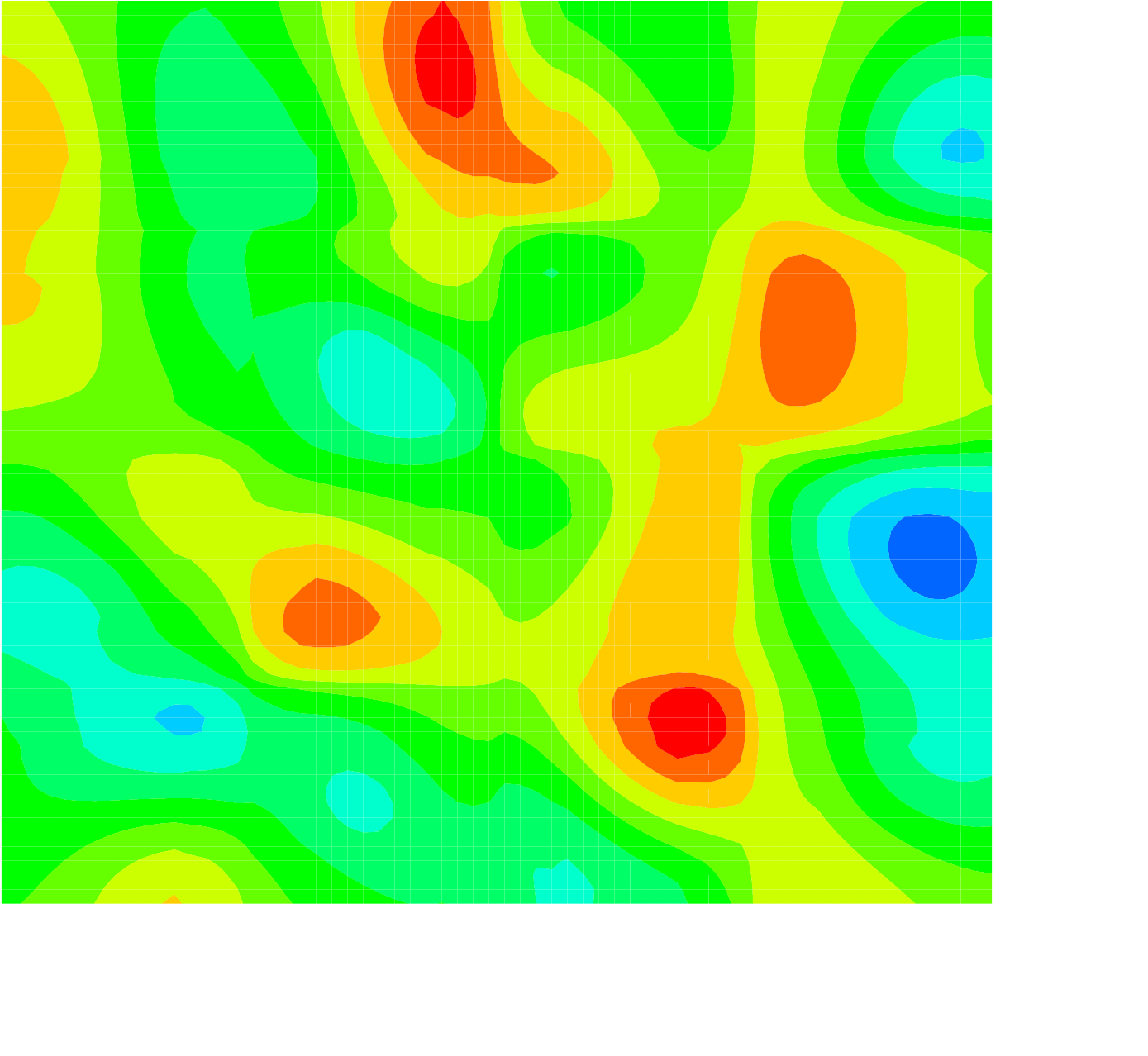}
		\includegraphics[width= 1.5 in] {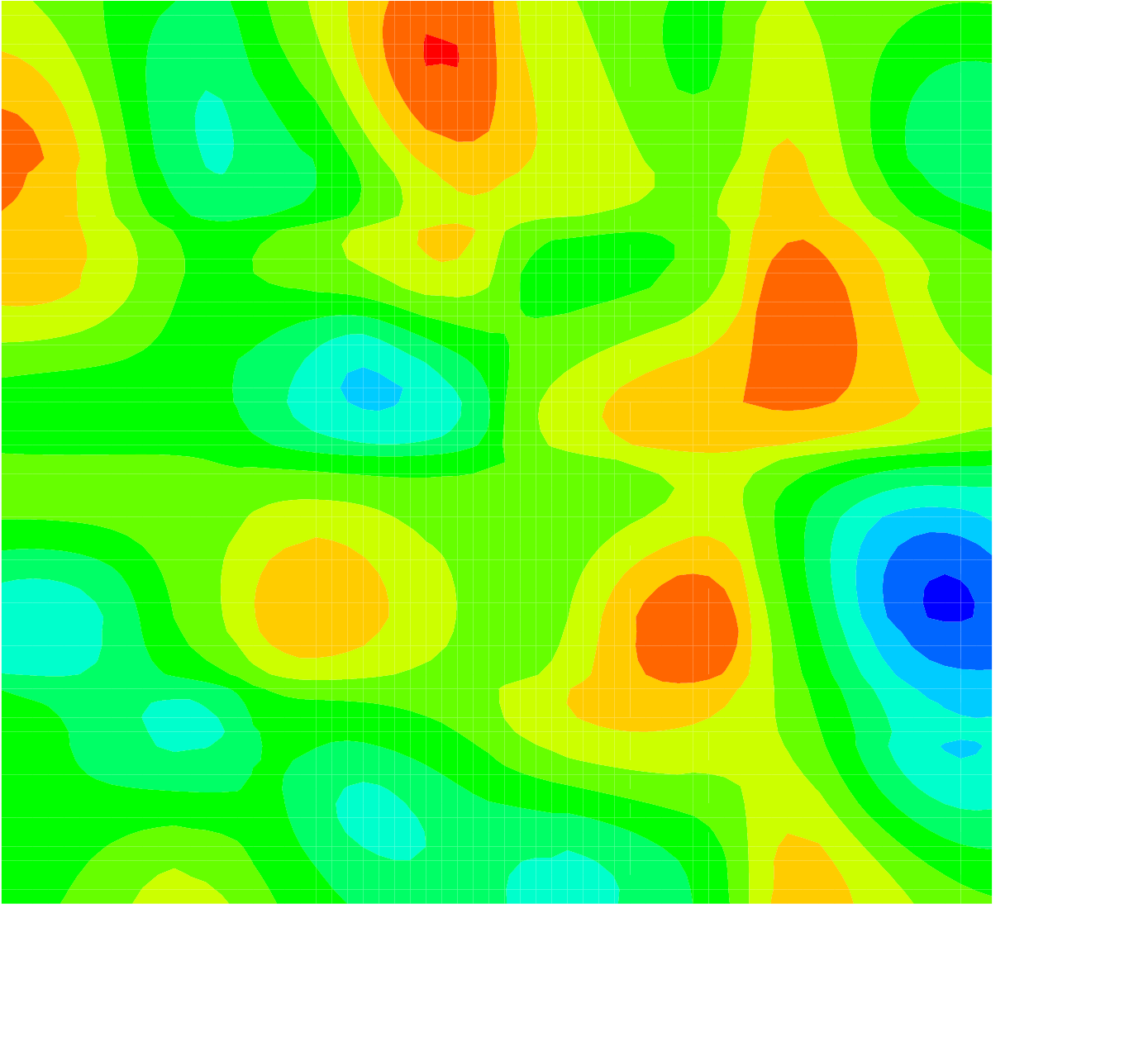}\\
		\includegraphics[width= 1.5 in] {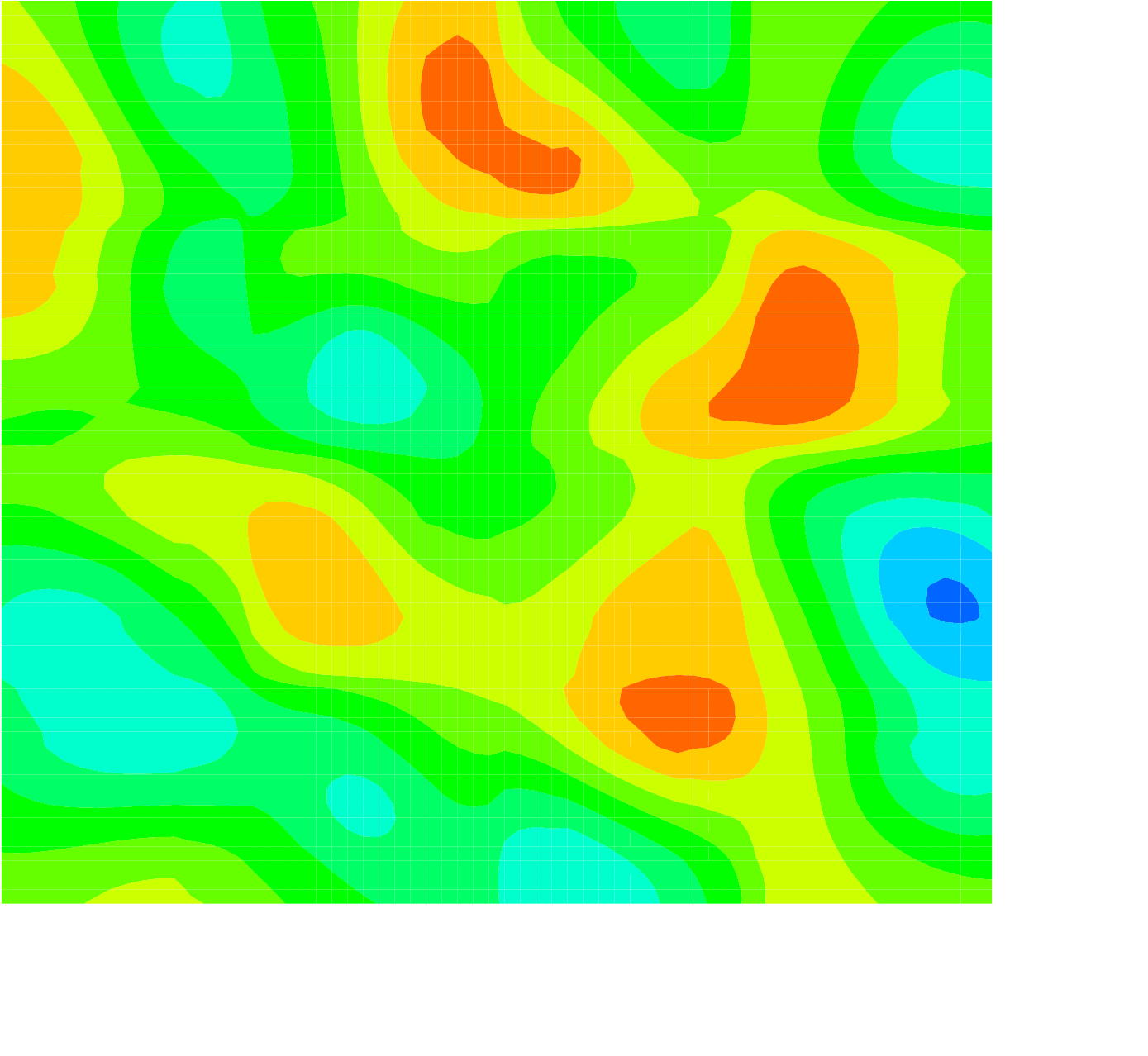}
		\includegraphics[width= 1.5 in] {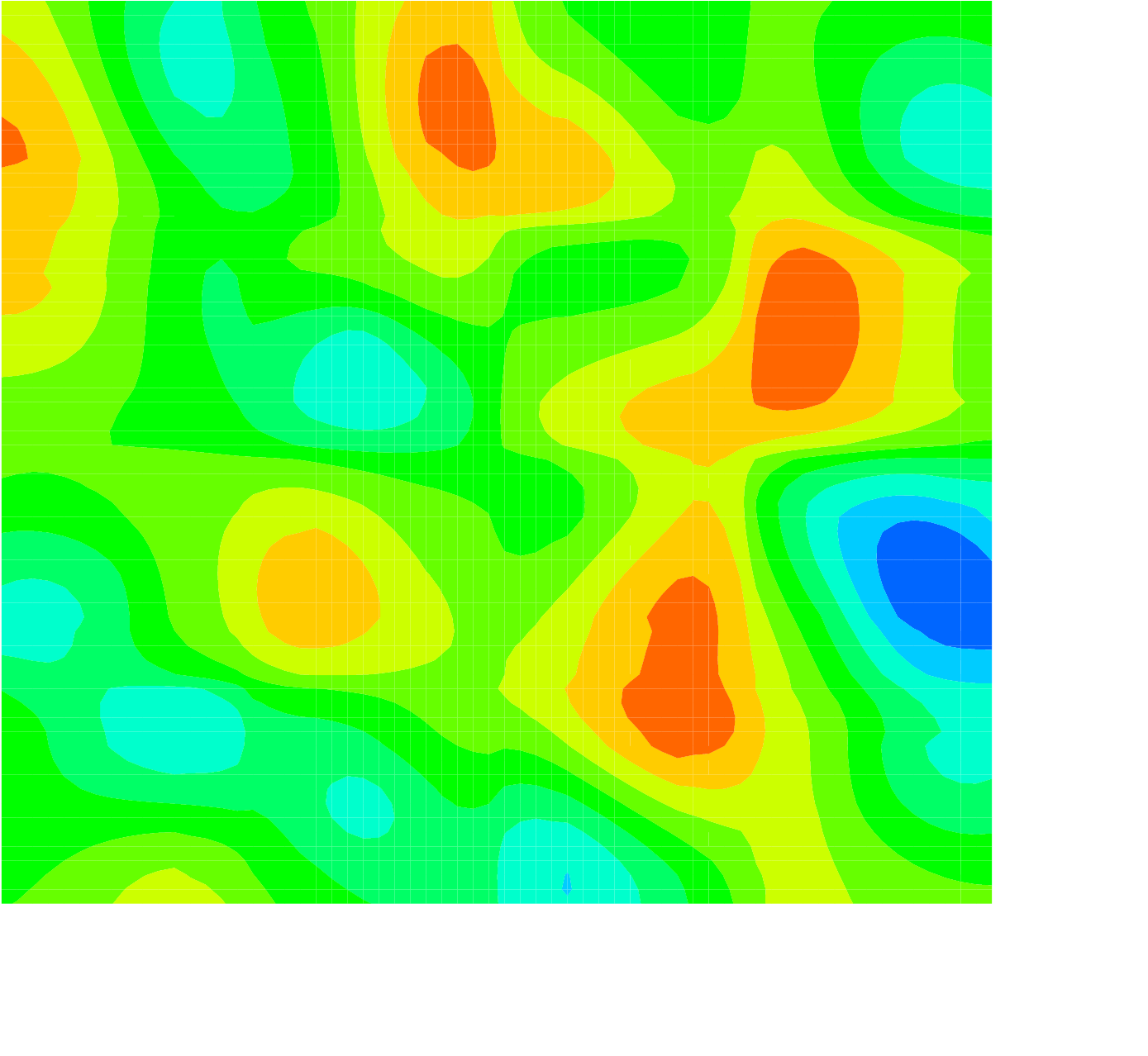}
  	  \includegraphics[width= 1.5 in] {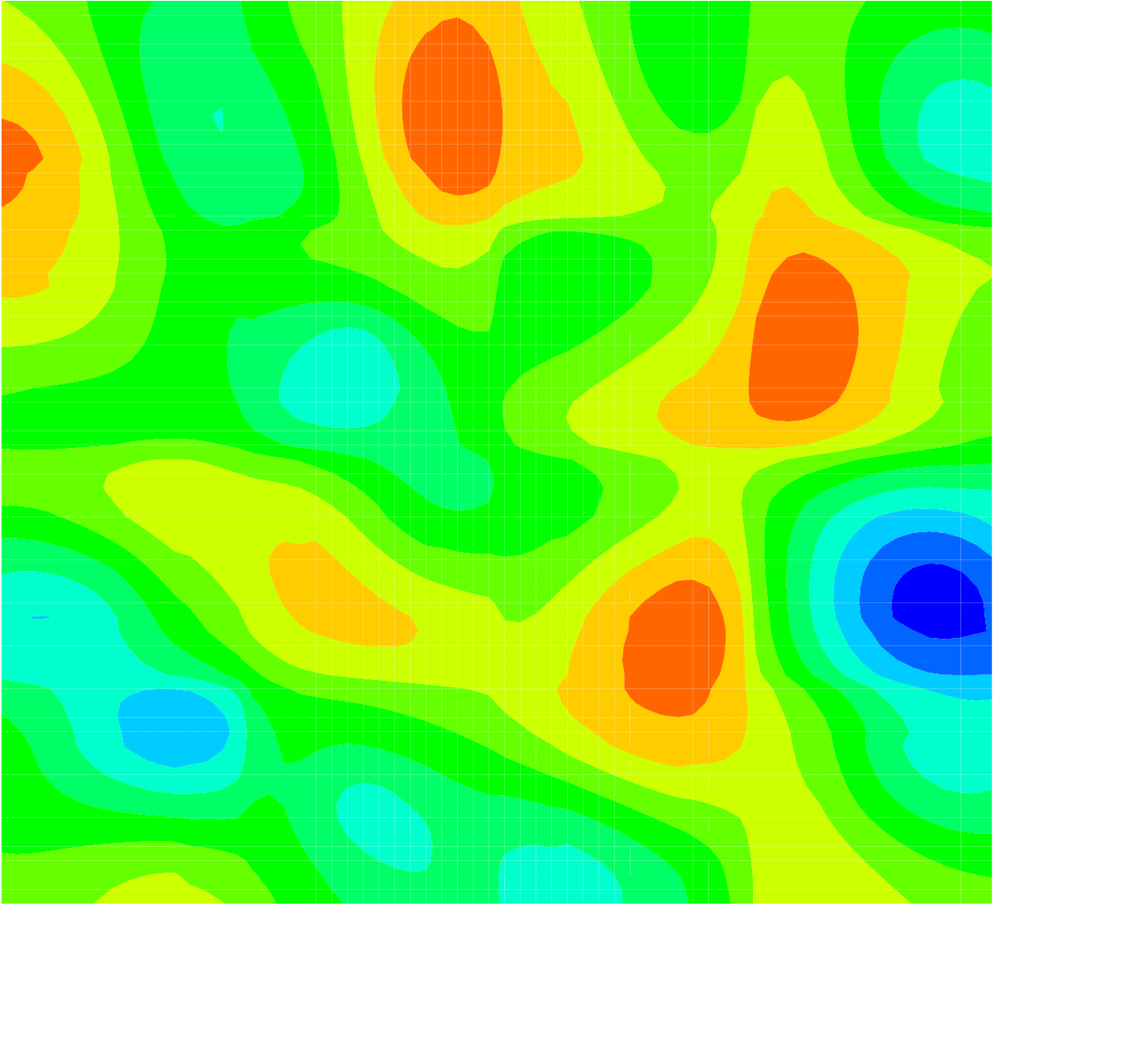}	
					
		\caption{ First row: Reference log permeability filed. Second row: Accepted permeability fields from chain 1. Third row: Accepted permeability fields from chain 2. 
	From left to right, log permeability fields at 40000, 60000 and 80000 iterations, respectively, in MSM $4\times 4$ with conditioning.
		}
		\label{perm_cond_64}
	\end{figure}

Fig. \ref{perm_cond_64} shows the accepted fields in the method with conditioning for two chains. The permeability fields were not fully recovered, however, we see a considerable improvement in the fields in comparison to the fields in Figures \ref{perm_3_chain1} and \ref{perm_3_chain2}, which were obtained using MSM $4\times4$ without conditioning. We thus conclude that the conditioning speeds-up the convergence and improves the characterization in this example

\section{Conclusions}
\label{concl_3}
We have presented a novel multiscale sampling method aiming at subsurface characterization. The proposed method is based on a non-overlapping partition of the domain of the governing partial differential equation that leads to the localization of the search in the underlying stochastic space. 
The novel method is implemented in the framework of a preconditioned Markov Chain Monte Carlo algorithm.

Through several multi-chain MCMC examples, motivated by subsurface flow problems, we compare the usual preconditioned Markov Chain Monte 
Carlo algorithm with the proposed procedure. Our results show that the new multiscale sampling method considerably
improves the convergence rate of the preconditioned Markov Chain Monte Carlo algorithm. We also incorporated sparse measurements of the permeability field in the multiscale sampling method, and showed that conditioning on this data further improves the 
convergence of the proposed method.

The authors and their collaborators are currently applying the method introduced here to solve the inverse problems associated with single and multiphase flows in porous media. 
In these studies the Multiscale Perturbation Method \cite{mpm_2020} will be used to speed-up the numerical solution of elliptic equations in the forward solution of the governing system of equations. Multiscale sampling procedures based on overlapping domain decompositions are also being considered.

\section*{Acknowledgments}
A. Rahunanthan research was supported by NIFA/USDA through Central State University Evans-Allen Research Program.

\bibliographystyle{model1-num-names}
\bibliography{alsadig_ref.bib}

\end{document}